%% file: article.tex
\documentclass[sn-mathphys-num]{sn-jnl}

\usepackage[utf8]{inputenc}
\usepackage{graphicx}%
\usepackage{multirow}%
\usepackage{amsmath,amssymb,amsfonts}%
\usepackage{amsthm}%
\usepackage{longtable}
\usepackage{mathrsfs}%
\usepackage []{csquotes}
\usepackage{enumitem}
\MakeOuterQuote{"}
\usepackage[title]{appendix}%
\usepackage[table]{xcolor}
\usepackage{textcomp}%
\usepackage{manyfoot}%
\usepackage{booktabs}%
\usepackage{algorithm}%
\usepackage{algorithmicx}%
\usepackage{algpseudocode}%
\usepackage{listings}%
\usepackage{svg}
\usepackage{rotating}
\usepackage{subcaption}
\usepackage{refcount}
\usepackage{pdflscape}
\usepackage{geometry}
\usepackage{xcolor}
\usepackage{longtable,tabularx}
\usepackage{tabularray}

\newcommand{\noun}[1]{\textsc{#1}}

\begin{document}

\title[Article Title]{System Architecture Optimization Strategies:
Dealing with Expensive Hierarchical Problems}

\author*[1]{\fnm{Jasper H.} \sur{Bussemaker}}\email{jasper.bussemaker@dlr.de}

\author[2]{\fnm{Paul} \sur{Saves}}\email{paul.saves@onera.fr}
\author[2]{\fnm{Nathalie} \sur{Bartoli}}\email{nathalie.bartoli@onera.fr}
\author[2]{\fnm{Thierry} \sur{Lefebvre}}\email{thierry.lefebvre@onera.fr}
\author[2]{\fnm{Rémi} \sur{Lafage}}\email{remi.lafage@onera.fr}

\affil*[1]{\orgdiv{Institute of System Architectures in Aeronautics}, \orgname{German Aerospace Center (DLR)}, \orgaddress{\street{Hein-Saß-Weg 22}, \city{Hamburg}, \postcode{21129}, \country{Germany}}}

\affil[2]{\orgdiv{DTIS}, \orgname{ONERA, Université de Toulouse}, \orgaddress{\street{2 Avenue Edouard Belin}, \city{Toulouse}, \postcode{31000}, \country{France}}}
\affil[2]{\orgdiv{Fédération ENAC ISAE-SUPAERO ONERA}, \orgname{Université de Toulouse}, \orgaddress{\street{2 Avenue Edouard Belin}, \city{Toulouse}, \postcode{31000}, \country{France}}}

\abstract{
Choosing the right system architecture for the problem at hand is challenging due to the large design space and high uncertainty in the early stage of the design process. Formulating the architecting process as an optimization problem may mitigate some of these challenges.
This work investigates strategies for solving System Architecture Optimization (SAO) problems: expensive, black-box, hierarchical, mixed-discrete, constrained, multi-objective problems that may be subject to hidden constraints.
Imputation ratio, correction ratio, correction fraction, and max rate diversity metrics are defined for characterizing hierarchical design spaces.
This work considers two classes of optimization algorithms for SAO: Multi-Objective Evolutionary Algorithms (MOEA) such as NSGA-II, and Bayesian Optimization (BO) algorithms.
A new Gaussian process kernel is presented that enables modeling hierarchical categorical variables, extending previous work on modeling continuous and integer hierarchical variables.
Next, a hierarchical sampling algorithm that uses design space hierarchy to group design vectors by active design variables is developed.
Then, it is demonstrated that integrating more hierarchy information in the optimization algorithms yields better optimization results for BO algorithms.
Several realistic single-objective and multi-objective test problems are used for investigations.
Finally, the BO algorithm is applied to a jet engine architecture optimization problem.
This work shows that the developed BO algorithm can effectively solve the problem with one order of magnitude less function evaluations than NSGA-II.
The algorithms and problems used in this work are implemented in the open-source Python library \noun{SBArchOpt}.
}


\maketitle

\section*{Nomenclature}
{\renewcommand\arraystretch{1.0}
\noindent\begin{longtable*}{@{}l @{\quad=\quad} l@{}}
$c_h$ & Hidden constraint \\
$c_{v,l}$ & Value constraint $l$ \\
$\mathrm{CR}$ & Correction ratio ($d$ subscript: discrete $\mathrm{CR}$; $c$ subscript: continuous $\mathrm{CR}$) \\
$\mathrm{CRF}$ & Correction fraction \\
DoE & Design of experiments\\
$f_{m}$ & Objective function $m$ \\
$g_k$ & Inequality constraint $k$ \\
HV & Hypervolume \\
$\mathrm{IR}$ & Imputation ratio ($d$ subscript: discrete $\mathrm{IR}$; $c$ subscript: continuous $\mathrm{IR}$) \\
$k_{\mathrm{doe}}$ & DoE multiplier \\
PLS & Partial least squares regression \\
KPLS & Kriging with PLS \\
$\mathrm{MRD}$ & Maximum rate diversity \\
$n_{c_v}$ & Number of value constraints \\
$n_f$ & Number of objectives \\
$n_g$ & Number of inequality constraints \\
$n_{\mathrm{act}}$ & Number of active design variables in a design vector \\
$n_{\mathrm{batch}}$ & Batch size for infill points \\
$n_{\mathrm{corr,discr}}$ & Number of correct discrete design vectors \\
$n_{\mathrm{doe}}$ & Design of experiments size \\
$n_{\mathrm{infill}}$ & Total number of infill points generated \\
$n_{\mathrm{kpls}}$ & PLS components for KPLS \\
$n_{\mathrm{parallel}}$ & Maximum number of parallel evaluations \\
$n_{\mathrm{valid,discr}}$ & Number of valid discrete design vectors \\
$n_x$ & Number of design variables \\
$n_{x_c}$ & Number of continuous design variables \\
$n_{x_d}$ & Number of discrete design variables \\
$n_{x,\mathrm{grp}}$ & Number of design variables in a sampling group \\
$N_{\mathrm{fe}}$ & Number of function evaluations \\
$N_j$ & Number of options for discrete variable $j$ \\
$j$ & Rate diversity of discrete variable $j$ \\
$\hat{s}$ & Uncertainty estimate by a surrogate model \\
TSFC & Thrust-specific fuel consumption \\
$w$ & Group weight for sampling \\
$\boldsymbol{x}$ & Design vector \\
$x_{\mathrm{act}}$ & Active design variables in a design vector \\
$x_{c,i}$ & Continuous design variable $i$ \\
$x_{\mathrm{corr}}$ & Design vector to be corrected \\
$x_{d,j}$ & Discrete design variable $j$ \\
$x_{\mathrm{valid,discr}}$ & Valid discrete design vectors \\
$\hat{y}$ & Function value estimate by a surrogate model \\
$\delta_i$ & Activeness function for design variable $i$ \\
\end{longtable*}}

\section{Introduction}\label{sec1}

System architecture is an important aspect of any engineered system. It conceptually defines how the system meets its stakeholder expectations by specifying the components the system consists of, how these components fulfill system functions, and how the components are related to each other~\cite{Crawley2015}. 
Designing a system architecture is anything but trivial, involving requirements analysis, interviewing stakeholders, defining the top-level functions the system should fulfill, identifying appropriate technologies for fulfilling these functions, managing complexity, and working with disciplinary experts to come up with a solution~\cite{chan2022aircraft}.
This process is challenged by the fact that due to the combinatorial nature of architecture decisions, it is infeasible to enumerate, analyze, and compare over every possible architecture alternatives~\cite{Iacobucci2012}.
Additionally, engineering teams might be biased towards known solutions, preventing them to consider the full range of possibilities~\cite{McDermott2020}, mainly due to the time it takes to analyze one architecture alternative. \textit{System Architecture Optimization} (SAO) is an emerging field, where these challenges are addressed by combining numerical optimization techniques with quantitative performance evaluation of architecture alternatives to automatically (but not exhaustively) search the design space in search of the best architecture(s) for the problem at hand~\cite{Judt2016}.
To turn an architecting process into an architecture optimization problem, several conditions must be satisfied~\cite{Bussemaker2021}:
\begin{itemize}
    \item Architectural decisions need to be formally specified in an architectural design space model, such that a computer program can reason about them, try out different combinations of decisions, and explore promising architecture candidates automatically.
    \item It should be possible to quantitatively evaluate the performance of each architecture so that architecture candidates can be compared to each other objectively, by implementing a computational framework that can handle all possible architectures and captures relevant (physical) phenomena with sufficient details. 
    \item Appropriate numerical optimization algorithms must be available to solve the architecture optimization problem efficiently.
\end{itemize}

Past work on SAO has explored various methods for modeling architectural design spaces, such as using feature models~\cite{Czarnecki2012}, function-means trees~\cite{Gedell2012}, morphological matrices~\cite{Mavris2008,Chakraborty2016}, and various graph-based models~\cite{Simmons2008,Herber2020a}.
The reader is referred to~\cite{Bussemaker2024adsg} for a more comprehensive overview of architecture design space modeling methods.
Quantitative evaluation of architecture candidates can be implemented using Multidisciplinary Design Analysis and Optimization (MDAO) techniques, which enables integrating the diverse and coupled engineering disciplines usually involved in system architecting problems~\cite{SobieszczanskiSobieski2015}.
Recent effort in coupling system architecting with MDAO includes~\cite{Chaudemar2021,Helle2022a,Bussemaker2022}, which is further extended with the capability to automatically modify the MDAO workflow for each architecture candidate generated during the SAO loop~\cite{Sonneveld2023,Bruggeman2024,Garg2024mdo}.

This work focuses on the last point in the list above: numerical optimization algorithms that can handle all challenges of architecture optimization problems, namely mixed-discrete, hierarchical, constrained, multi-objective design spaces with expensive, black-box evaluation functions.
Previous work in this area has mostly focused on the application of evolutionary algorithms~\cite{Frank2016,Judt2016,Apaza2021}, as these provide good performance in mixed-discrete, multi-objective, constrained design spaces.
Surrogate-Based Optimization (SBO) algorithms have been applied for SAO by \noun{Bussemaker et al.}~\cite{Bussemaker2021,Bussemaker2024hc} to support the expensive nature of physics-based multidisciplinary evaluation functions. 
SBO has also found application to optimization problems 
featuring mixed-discrete, hierarchical, multi-objective design spaces 
and
expensive evaluation functions in the field of hyperparameter optimization~\cite{Feurer2019,Bischl2023}.
%
This work contributes to the application of global optimization algorithms for solving SAO problems by improving several aspects of both evolutionary and BO algorithms:
\begin{itemize}
    \item The definition of metrics for quantifying the level of hierarchy in hierarchical design spaces; 

    \item Specifically for BO: new Gaussian Process (GP) kernels for modeling hierarchical categorical variables; 

    \item The definition of three classes of test problems for supporting research into and development of optimization algorithms for SAO; 

    \item A new algorithm for sampling hierarchical design spaces; 

    \item A comparison between various a-priori correction algorithms; 
    and

    \item An overview and comparison of various strategies for integrating information about the hierarchical design space into optimization algorithms. 

\end{itemize}


This paper continues with a more detailed introduction of the challenges posed by architecture optimization in Section~\ref{sec:ArchOpt}, and a discussion of appropriate optimization algorithm classes and their implementations in Section~\ref{sec:ArchOptAlgos}.
Kernels for modeling categorical variables in hierarchical GP models are presented in Section~\ref{sec:HierModel}.
Section~\ref{sec:HierSampling} develops and investigates new sampling and correction algorithms for hierarchical design spaces, and Section~\ref{sec:HierIntegration} discusses and compares different levels of hierarchy integration in optimization algorithms.
The developed strategies are applied to a jet engine architecture optimization problem in Section~\ref{sec:JetEngineOpt} and Section~\ref{sec:Conclusions} concludes this paper.

\section{Characteristics of Architecture Optimization Problems}\label{sec:ArchOpt}

Optimization is the automation of a design task: the goal is to find one or more design vectors $\boldsymbol{x}$, representing specific design points, that minimize (or maximize) one or more objective functions $f(\boldsymbol{x})$, while satisfying design constraints $g(\boldsymbol{x}) \leq 0$~\cite{Martins2022}.  The design vectors (or design points) $\boldsymbol{x}$ are composed of one value for each design variable $x$.
Note that maximization of an objective function is enabled by negating the objective function value, and therefore we treat minimization as the default in the rest of this work.
In the context of System Architecture Optimization (SAO), a design point $\boldsymbol{x}$ represents an architecture instance. Compared to for example the work by \noun{Frank \textit{et al.}}~\cite{Frank2016}, we treat the terms "architecture instance", "architecture alternative", and "design point" as synonymous meaning one possible design vector in the architecture design space.

\subsection{Design Space Characteristics}\label{subsec:DSchallenges}

SAO problems are \textit{mixed-discrete}: both continuous and discrete design variables can be part of the problem formulation~\cite{Frank2016,saves2021constrained}. Continuous design variables $x_c$  can take any value between some lower bound $\underline{x}_{c}$ and some upper bound $\overline{x}_{c}$ (inclusive). Discrete design variables $x_d$ can only take one from a finite set of values, encoded as an index between $0$ and $N_j-1$ (inclusive), with $N_j$ representing the number of possible discrete values for discrete design variable $x_{d,j}$.
Discrete variables can either be of \textit{integer}, \textit{ordinal} or of \textit{categorical} type: for integer and ordinal variables the order of the values is relevant, for categorical values no notion of value order exists~\cite{Saves2023a}. The difference between integer and ordinal variables is that integer values operate between two bounds and the variable can take any integer value between these bounds, whereas ordinal variables can take any value of a specified list of integer values.
Examples of continuous design variables are wing sweep and engine bypass ratio, an example of an integer variable is the number of seat rows in an aircraft cabin (\textit{e.g.} between 20 and 40), an example of an ordinal variable is the number of engines for an aircraft (\textit{e.g.} 1, 2 or 4; order is relevant), and an example of a categorical variable is aircraft power source (there is no order between kerosene, hydrogen and electricity). In the context of discrete architecture choices, function-to-component allocation and component connection are categorical choices, and component or system instantiations are integer or ordinal choices~\cite{Bussemaker2022c}.

SAO design spaces feature strong interaction between decisions, in the sense that some decisions might affect which other decisions remain to be taken and/or which options remain available for other decisions.
For example, consider the Apollo mission architecture problem analyzed by \noun{Simmons}~\cite{Simmons2008}: the decision whether or not a lunar-orbit-rendezvous or earth-orbit-rendezvous maneuver will be performed is only relevant if an architecture with a lunar module is selected. Another example from the same architecture problem is crew assignment: the lunar module can have 1, 2 or 3 crew members, except if there are only 2 crew members in the command module, in which case the lunar module can only contain 1 or 2 crew members.
In the launch vehicle design problem by \noun{Pelamatti \textit{et al.}}~\cite{Pelamatti2020a}, the number of stages is a decision (1 to 3), as well as several decisions such as fuel type and dimensions for each of the stages: it follows that if a one-stage architecture is selected, the decisions regarding the second and third stages do not have to be considered. Machine learning hyperparameters tuning problems also exhibit such interactions between design variables~\cite{Feurer2019}.
In general, decisions constraining available options of other decisions is a common theme in architecture optimization, for example defined using an incompatibility matrix~\cite{Armstrong2008}. These interactions between design variables create a \textit{hierarchical} structure in the design space~\cite{Zaefferer2018a}, with design variables higher up in the hierarchy determining which variables lower in the hierarchy are \textit{active} or \textit{inactive}, and/or which options are available for active design variables.
The property of being active or inactive is called \textit{activeness}~\cite{Bussemaker2021}, and is defined through the activeness function $\delta_{i}(\boldsymbol{x})$~\cite{Hutter2013}, which returns a $1$ indicating active or $0$ indicating inactive for design variable $x_{i}$. Variables that may be inactive are denoted as \textit{conditionally active}.
In architecture optimization problems often only discrete design variables determine the activeness of other design variables, however in the general case also continuous variables could determine activeness as in the work by \noun{Zaefferer \& Horn}~\cite{Zaefferer2018a}.
Hierarchical design spaces are also known as conditional design spaces~\cite{Feurer2019} as design variables can be conditionally active, design spaces with tree-structured dependencies~\cite{Bergstra2011,Jenatton2017,Apaza2021}, and variable-size design spaces~\cite{Abdelkhalik2012,Talbi2024} as inactive design variables can also be seen as artifacts of a locally-reduced size of the design space.

Inactive design variables do not influence the performance of the design and are therefore irrelevant for objective and constraint evaluations. If an optimizer is not aware of this, it can however still generate design vectors with different values for inactive design variables, resulting in the possibility of multiple design vectors representing the same design (\textit{i.e.} the design vector to design mapping is no longer one-to-one) and therefore the same objective and constraint values, and thereby wasting computational resources or preventing the optimizer from finding the optimum.
To mitigate this, each inactive design variable can be assigned a canonical value, for example $0$ for discrete variables and mid-bounds for continuous variables~\cite{Zaefferer2018a}. The process of replacing inactive design variables with canonical values is called design vector \textit{imputation}~\cite{Zaefferer2018a,Levesque2017}, and the resulting design vector is called an imputed or canonical design vector. 
Restricting option availability of active variables lower in the hierarchy is done using \textit{value constraints}, also known as enumeration constraints in the context of architecture optimization~\cite{Selva2012} or configuration constraints in product line engineering~\cite{Czarnecki2012,Weilkiens2015}.
According to the taxonomy by \noun{Le Digabel \& Wild}~\cite{LeDigabel2023} value constraints are non-quantifiable (it is only know whether they are violated, not by how much), unrelaxable (a design point with a failed value constraint does not represent anything that can be evaluated), \textit{a-priori} (they can be corrected without running an evaluation), and known (they are derived from the design space definition).
The process of ensuring value constraints are satisfied by modifying a design vector is called \textit{correction}.

\begin{figure}
\centering
\centering
\includesvg[inkscapelatex=false,width=\textwidth]{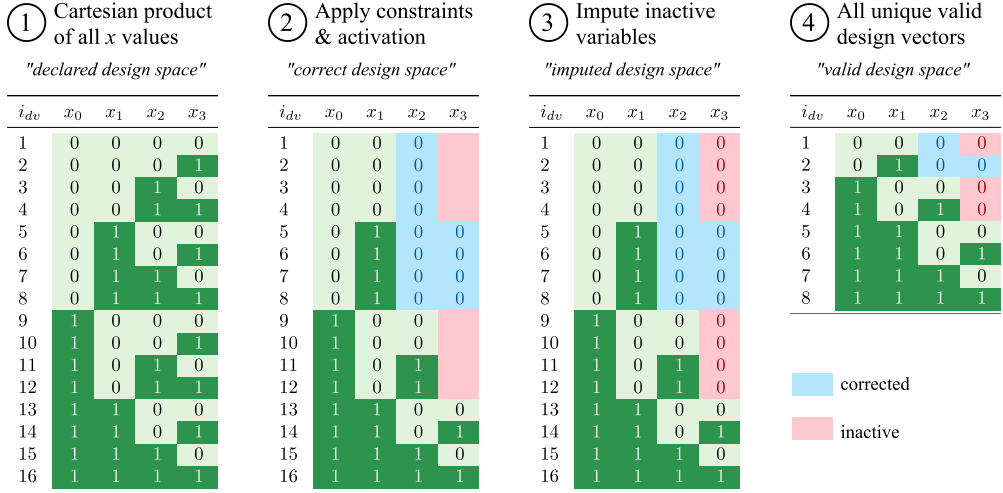}
\caption{Illustration of correction and imputation in hierarchical design spaces, showing how the Cartesian product of all discrete values relates to the set of correct, imputed, and valid design vectors. Correction (step 2) modifies design variables such that value constraints are satisfied. Imputation (step 3) assigns canonical values to inactive design variables.
$\delta_i$ represents the activeness function of design variable $x_i$.}\label{fig:corr_imp}
\end{figure}

A result of design variable hierarchy is that it is not possible to compute the number of discrete design options and continuous variables from the design variable definitions alone: for a non-hierarchical design space, the number of valid discrete design points $n_{\mathrm{valid,discr}}$ is equal to the Cartesian product of all discrete values $\prod N_{j}$. For hierarchical design spaces, however, the number of valid discrete design points can be much lower than that.
For example, for the suborbital vehicle design problem of \noun{Frank \textit{et al.}}~\cite{Frank2016} there are 2.8 million total combinations of discrete variables, however there are only 123 thousand (4.4\% of 2.8 million) valid designs. \noun{Ying \textit{et al.}}~\cite{Ying2019} report 423 thousand unique neural networks in a Cartesian design space of 510 million design vectors, a difference of a ratio 1206.
To highlight such differences, we distinguish between the \textit{declared} discrete design space, given by the aforementioned Cartesian product, and the \textit{valid} discrete design space, given by the number of valid discrete design vectors $n_{\mathrm{valid,discr}}$.
Correction combined with imputation moves a design point from the declared design space to the valid design space: correction ensures a design point is made \textit{correct}, imputation ensures that correct \textit{non-canonical} design points are made canonical and thereby \textit{valid} (correct and canonical).
The principles of correction and imputation are visualized in Figure~\ref{fig:corr_imp}. A hypothetical source-to-target assignment problem is shown where 1 or 2 energy consumers are assigned to 1 or 2 energy sources.
Two hierarchical interactions arise: if there is only one source, then both consumers have to be assigned to that source (a value constraint), and if there is only one consumer, the source choice for the second consumer is not active (hierarchical activation).
To move from the declared (Figure~\ref{fig:corr_imp} step 1) to valid design space, first design vectors are corrected (step 2) and then inactive design variables are imputed (step 3). The valid design space (step 4) contains the unique remaining design vectors.
For this example, the declared design space (step 1) consists of 16 design vectors, whereas the valid design space (step 4) only consists of 8 design vectors, demonstrating this discrepancy typical for hierarchical optimization problems.
The combined operations of correction and imputation can be implemented in an optimization problem using a repair operator~\cite{SalcedoSanz2009,Selva2012,Koch2015}.
Table~\ref{tab:x_cond} provides an overview of relevant design point statuses in the context of SAO.

\subsubsection{Imputation and Correction Ratio}

Due to hierarchy, the valid design space might be much smaller than the declared design space.
Since this might pose a challenge to optimization algorithms we now introduce a metric to quantify this discrepancy: the ratio between the declared and valid (see Table~\ref{tab:x_cond}) discrete design space sizes is defined as the discrete \textit{imputation ratio} $\mathrm{IR}_d$:
\begin{equation}\label{eq:IRd}
    \mathrm{IR}_d=\frac{\mathop{\prod_{j=1}^{n_{x_d}}N_{j}}}{n_{\mathrm{valid,discr}}}
\end{equation}
where $n_{\mathrm{valid,discr}}$ is the number of valid discrete design vectors, and $N_j$ is the number of options for discrete variable $j$. An imputation ratio of 1 indicates a non-hierarchical problem, values higher than 1 indicate hierarchy: for the suborbital vehicle design problem in~\cite{Frank2016}, the imputation ratio is \(2.8e6 / 123e3 = 22.8\).
The higher the value, the more invalid (non-canonical and non-corrected) vectors would be generated in a random search, and therefore the more difficult it is for an optimization algorithm to search the design space if this effect is not considered.
Variability factor~\cite{Benavides2010} as used in software product line engineering is the reciprocal of discrete imputation ratio.
The continuous imputation ratio $\mathrm{IR}_c$ is defined as follows:
\begin{equation}\label{eq:IRc}
    \mathrm{IR}_{c}=\frac{n_{\mathrm{valid,discr}}\cdot n_{x_c}}{\mathop{\sum_{l=1}^{n_{\mathrm{valid,discr}}}}\mathop{\sum_{i=1}^{n_{x_c}}\delta_{i}(\boldsymbol{x}_{d,l})}}
\end{equation}
where $\delta_i(\boldsymbol{x}_{d,l})$ is the activeness function for continuous variable $i$ for valid discrete design vector $x_{d,l}$, $n_{\mathrm{valid,discr}}$ the number of valid discrete design vectors, and $n_{x_c}$ the number of continuous design variables. A value of 1 indicates that all continuous variables are always active. A higher value indicates that for some discrete design vectors one or more continuous variables are inactive.
Note that this formulation assumes that only discrete design variables determine activeness of continuous design variables.
The overall imputation ratio for a given optimization problem is given by the product of the two:
\begin{equation}\label{eq:IR}
    \mathrm{IR} = \mathrm{IR}_d \cdot \mathrm{IR}_c
\end{equation}
To give an example, the simple engine architecture benchmark problem from~\cite{Bussemaker2021c} has $n_{\mathrm{valid,discr}} = 70$ valid discrete vectors, however the Cartesian product of its discrete design variables indicates 216 declared discrete vectors, resulting in a discrete imputation ratio of $\mathrm{IR}_d = 216/70 = 3.09$. The problem also contains 9 continuous design variables, of which 7.14 are active on average as seen over all discrete design vectors, which result in a continuous imputation ratio of $\mathrm{IR}_c = 9/7.14 = 1.26$. The overall imputation ratio therefore is $\mathrm{IR} = 3.09 \cdot 1.26 = 3.89$.
To compare, the example problem from Figure~\ref{fig:corr_imp} has an imputation ratio of $\mathrm{IR} = \mathrm{IR}_d = 16 / 8 = 2$.

Similarly to the discrepancy between the declared and valid design space sizes, we can also quantify the discrepancy between the declared and correct (see Table~\ref{tab:x_cond}) design spaces sizes. This can help determine how much of the design space hierarchy is due to value constraints that need correction, as opposed to design variable activeness.
We define the \textit{correction ratio} $\mathrm{CR}$ as:
\begin{align}
    \mathrm{CR}_d &= \frac{\mathop{\prod_{j=1}^{n_{x_d}}N_{j}}}{n_{\mathrm{corr,discr}}} \label{eq:CRd}\\
    \mathrm{CR}_c &= \frac{n_{\mathrm{corr,discr}}\cdot n_{x_c}}{\mathop{\sum_{l=1}^{n_{\mathrm{corr,discr}}}}\mathop{\sum_{i=1}^{n_{x_c}}\delta_{i}(\boldsymbol{x}_{d,l})}} \label{eq:CRc}\\
    \mathrm{CR} &= \mathrm{CR}_d \cdot \mathrm{CR}_c \label{eq:CR}
\end{align}
where $\mathrm{CR}_d$ is the discrete correction ratio, $n_{\mathrm{corr,discr}}$ the number of correct discrete design vectors, and $\mathrm{CR}_c$ the continuous correction ratio.
The impact of the need of correction to the design space hierarchy can be quantified by the \textit{correction fraction} $\mathrm{CRF}$:
\begin{equation}\label{eq:CRF}
    \mathrm{CRF} = \frac{\log \mathrm{CR}}{\log \mathrm{IR}}
\end{equation}
$\mathrm{CRF}$ varies between 0\% and 100\%, where 0\% indicates no hierarchy is due to correction ($\mathrm{CR} = 1$) and 100\% indicates all hierarchy is due to correction ($\mathrm{CR} = IR$, and there is no hierarchy due to activeness).
The valid design vectors shown in Table~\ref{tab:cr_example} represent a design space with a declared size of 12 ($N_0 \cdot N_1 = 4 \cdot 3 = 12$), however $n_{\mathrm{valid,discr}} = 6$, and therefore $\mathrm{IR} = \mathrm{IR}_d = 12 / 6 = 2.0$.
Additionally, the example has $n_{\mathrm{corr,discr}} = 10$, because the inactive design variables can take any declared value and still represent a correct (but not necessarily valid) design vector. Therefore, $\mathrm{CR} = \mathrm{CR}_d = 12 / 10 = 1.2$. From this, $\mathrm{CRF} = \log 1.2 / \log 2.0 = 26\%$, indicating that 26\% of the design space hierarchy is due to the need for correction, while the other 74\% are due to activeness and the need for imputation.

\begin{table}
\caption{Example set of valid design vectors, showing inactive variables in red.
$x_0$ has 4 possible values ($N_0 = 4$); similarly for $x_1$: $N_1 = 3$.}\label{tab:cr_example}
\begin{tabular}{c}
\hspace{3cm}\begin{tabular}{lll}
\toprule
$i_{dv}$ & $x_{0}$ & $x_{1}$ \\
\midrule
1 & {\cellcolor[HTML]{E1F3DC}} \color[HTML]{000000} 0 & {\cellcolor[HTML]{E1F3DC}} \color[HTML]{000000} 0 \\
2 & {\cellcolor[HTML]{E1F3DC}} \color[HTML]{000000} 0 & {\cellcolor[HTML]{AEDEA7}} \color[HTML]{000000} 1 \\
3 & {\cellcolor[HTML]{AEDEA7}} \color[HTML]{000000} 1 & {\cellcolor[HTML]{E1F3DC}} \color[HTML]{000000} 0 \\
4 & {\cellcolor[HTML]{AEDEA7}} \color[HTML]{000000} 1 & {\cellcolor[HTML]{6ABF71}} \color[HTML]{000000} 2 \\
5 & {\cellcolor[HTML]{6ABF71}} \color[HTML]{000000} 2 & {\cellcolor[HTML]{FFC7CE}} \color[HTML]{9C0006} \\
6 & {\cellcolor[HTML]{2C944C}} \color[HTML]{F1F1F1} 3 & {\cellcolor[HTML]{FFC7CE}} \color[HTML]{9C0006} \\
\bottomrule
\end{tabular}\hspace{3cm}
\end{tabular}
\end{table}

\subsubsection{Rate Diversity}

In addition to the potentially large gap between declared and valid design space sizes, there might also be a discrepancy in how often individual discrete values appear in all possible discrete design vectors, as also noted by \noun{Crawley \textit{et al.}}~\cite{Crawley2015}. They mention that for a partitioning problem of 10 elements, where the goal is to choose the best subset of elements from the available set, there are 115 thousand possibilities to choose from, however 88\% of those possibilities are made up of the size-4, -5, and -6 subsets. This means that all other sizes are represented much less in the total number of possibilities.
This observation can be extended to design variable values too: for the realistic engine architecting benchmark problem of \noun{Bussemaker \textit{et al.}}~\cite{Bussemaker2021c}, there are a little over 142 thousand valid architectures, however only 27 of those (about 0.02\%) represent a pure jet engine architecture; the rest are turbofan architectures. Therefore there is a large gap between how often the two possible values of this design variable appear in all discrete design vectors.
This gap can be quantified by the \textit{rate diversity} $\mathrm{RD}_{j}$, which is defined for each discrete design variable $x_{d,j}$, and the \textit{max rate diversity} $\mathrm{MRD}$, which is defined on the problem level:
\begin{align}
    \mathrm{Rates}_j &= \left\{ \mathrm{Rate}(j,\delta_j = 0), \mathrm{Rate}(j,0), .., \mathrm{Rate}(j,N_j-1) \right\} \\
    \mathrm{RD}_j &= \max \mathrm{Rates}_j - \min \mathrm{Rates}_j \\
    \mathrm{MRD} &= \max_{j \in 1,..,n_{x_d}} \mathrm{RD}_j \label{eq:MRD}
\end{align}
with $j$ the index of the discrete design variable, $\delta_j = 0$ denoting the case when design variable $j$ is inactive, and $\mathrm{Rate}(j,\mathrm{value})$ returning the relative occurrence rate of that value in all discrete design vectors. Rate diversity $\mathrm{RD}_j$ then represents the difference between the most and least occurring values for a given discrete design variable $j$, and max rate diversity $\mathrm{MRD}$ the maximum of all rate diversities.
Rate diversity and max rate diversity are normally defined without the inactive case $\mathrm{Rate}(j,\delta_j = 0)$ included. If the inactive case is included, the subscript "all" is added to denote all cases and values are considered.
Table~\ref{tab:tf_rates} shows occurrence rates, the rate diversity and maximum rate diversities of the simple turbofan benchmark problem from~\cite{Bussemaker2021c}. The rate diversity of the choice between pure jet and turbofan architectures $x_0$ is not as extreme as for the aforementioned realistic benchmark, however there is still a rate diversity of 60\% as only 20\% of possible discrete design vectors represent a jet engine architecture.
The example problem from Figure~\ref{fig:corr_imp} has two variables ($x_0$ and $x_3$) where one value occurs 2 of 8 times (25\%) and the other therefore 75\%, which results in an max rate diversity of $\mathrm{MRD} = 50\%$.

\begin{table}
\caption{Rate diversity $\mathrm{RD}$ of the simple turbofan architecting problem from~\cite{Bussemaker2021c}. The maximum rate diversity $\mathrm{MRD}$ values are underlined.
Only discrete variables are shown, as rate diversity does not apply to continuous variables.}\label{tab:tf_rates}
\begin{tabular}{llllllll}
\toprule
& $x_1$ & $x_4$ & $x_{11}$ & $x_{13}$ & $x_{14}$ & $x_{15}$ \\
\midrule
Inactive & $-$ & $-$ & 20\% & 20\% & 7.1\% & 7.1\% \\
$x_j = 0$ & 20\% & 7.1\% & 40\% & 40\% & 35.7\% & 35.7\% \\
$x_j = 1$ & 80\% & 28.6\% & 40\% & 40\% & 35.7\% & 35.7\% \\
$x_j = 2$ & $-$ & 64.3\% & $-$ & $-$ & 21.4\% & 21.4\% \\
$\mathrm{RD}_{all}$ & \underline{60\%} & 57.1\% & 20\% & 20\% & 28.6\% & 28.6\% \\
$\mathrm{RD}$ & \underline{60\%} & 57.1\% & 0\% & 0\% & 15.4\% & 15.4\% \\
 \bottomrule
\end{tabular}
\end{table}

Imputation ratio $\mathrm{IR}$ (Equation~\ref{eq:IR}), correction ratio $\mathrm{CR}$ (Equation~\ref{eq:CR}) and max rate diversity $\mathrm{MRD}$ (Equation~\ref{eq:MRD}) can be useful for quantifying how large the impact of hierarchy is for a given problem.

\subsection{Solution Space Characteristics}

The objective and design constraint functions $f(\boldsymbol{x})$ and $g(\boldsymbol{x})$, also known as the \textit{evaluation functions}, are \textit{black-box} functions: their behavior is not known in advance because the architecture optimization method can be applied to design any system. As a consequence, the evaluation functions exhibit non-linear, non-smooth (the function may contain discontinuities or gaps) and multi-modal behaviors (there may be multiple local minima; this also implies that the function is non-convex), and gradients cannot be assumed to be available (also due to the mixed-discrete nature of the problem, see below)~\cite{Crawley2015}.
Evaluation functions are assumed to be deterministic: repeated function calls with the same inputs yield the same outputs.
In many cases the evaluation performed to obtain the objective and design constraint values is \textit{expensive} in terms of time and/or computational resources, due to the use of physics-based simulations and Multidisciplinary Design Analysis and Optimization (MDAO~\cite{SobieszczanskiSobieski2015}) approaches. The consequence of this is that the time to perform one evaluation can be orders of magnitude more than the time for one iteration of the optimization algorithm, and therefore the focus should be on reducing the amount of evaluations, rather than speeding up the optimization algorithm itself~\cite{Jones1998}.

When designing a system architecture, needs and goals of system stakeholders often conflict with each other~\cite{Crawley2015} and the associated architecting problem thus becomes a trade-off between these conflicting goals.
In general, therefore, the SAO problem should be considered to be a \textit{multi-objective} optimization problem~\cite{Judt2016}: the problem will have multiple conflicting objective functions $f$ to minimize at the same time.
Due to this conflicting nature, however, it becomes impossible to define one design point as optimal; rather, multi-objective optimization results in a set of Pareto-optimal design points~\cite{Miettinen1998}. This so-called Pareto-set is comprised of design points that are not dominating each other: a design point $\boldsymbol{x}_a$ dominates a design point $\boldsymbol{x}_b$ if $f_m(\boldsymbol{x}_a) \leq f_m(\boldsymbol{x}_b)$ for all $m$ and $f_m(\boldsymbol{x}_a) < f_m(\boldsymbol{x}_b)$ for at least one $m$.
The consequence of this is that within the Pareto-set, no design point is better than any other design point, and that if you want to improve one or more of the objective values by moving from one point to another, you will also always make at least one other objective value worse. The Pareto-front contains the objective values of the points in the Pareto-set.
The challenge that multi-objective optimization algorithms face is to simultaneously progress towards the Pareto-front and maintain sufficient diversity along the Pareto-front~\cite{Rudenko2004}.

\textit{Design constraints} $g(\boldsymbol{x})$ can arise from physical (\textit{e.g.} maximum material stresses or temperatures) or operational limitations (\textit{e.g.} maximum take-off field length, maximum wing span), for example.
Solving constrained optimization problems means that a design point can only be considered optimal if all the design constraints are satisfied.
Equality constraints can be eliminated by design variable substitution, by non-linearly solving implicit residual equations within the evaluation function~\cite{Martins2022}, or by using dedicated infill functions and GP models for BO~\cite{SEGO-UTB}.
However, in this work, only \textit{inequality} constraints are considered.
We assume that the design constraints are QRSK in the taxonomy of \noun{Le Digabel \& Wild}~\cite{LeDigabel2023}: quantifiable (the degree of feasibility can be quantified), relaxable (a violated constraint is still meaningful to the optimizer), simulation-based (to know whether the constraint is satisfied, an evaluation must be run), and known (the constraint is defined in the problem formulation).
Design points where one or more design constraints are violated are called \textit{infeasible}, as opposed to \textit{feasible} if all constraints are satisfied.

Simulations used for evaluating architecture performance can fail, for example due to an unstable system of equations, infeasible underlying physics, or infeasible geometric parameterization~\cite{Forrester2006}. Any problem that employs simulation could include a so-called \textit{hidden constraint} (also known as unknown, unspecified, forgotten, virtual, and crash constraints~\cite{Mueller2019,LeDigabel2023}): a constraint that manifests itself through failed evaluations. The objective $f$ and design constraint $g$ function outputs are assigned NaN (Not a Number) values when a hidden constraint is violated.
We assume that hidden constraints are deterministic, resulting in the same status when repeatedly evaluating the same design point. Additionally, finding out about if the hidden constraint is violated takes at least as long as completing a successful evaluation takes, if not longer~\cite{Mueller2019}.
The region of the design space where the hidden constraint is violated is called the \textit{failed} region consisting of \textit{failed} design points, the region where this is not the case is called the \textit{viable} region consisting of \textit{viable} points.
The failed region can span a large part of the design space: \noun{Krengel \& Hepperle} report up to 60\% for a wing design problem~\cite{Krengel2022}, and for an airfoil design problem \noun{Forrester \textit{et al.}} report up to 82\% failed points~\cite{Forrester2006}.

\subsection{Problem Formulation}

Combining all relevant behavioral aspects, an SAO problem can be formulated as:
\begin{equation}\begin{array}{cclll}   
\mathrm{minimize} & & f_{m}(\boldsymbol{x},\delta(\boldsymbol{x})), &  & m=1,2,\ldots,n_f\\
\mathrm{w.r.t.} & & \underline{x}_{c,i}\leq x_{c,i}\leq\overline{x}_{c,i} &  & i=1,2,\ldots,n_{x_c}\\
 & & x_{d,j}\in\{0,..,N_{j}-1\} &  & j=1,2,\ldots,n_{x_d}\\
\mathrm{subject\ to} & & g_{k}(\boldsymbol{x},\delta(\boldsymbol{x}))\leq0, &  & k=1,2,\ldots,n_g\\
 & & c_{v,l}(\boldsymbol{x})=0 &  & l=1,2,\ldots,n_{c_v}\\
 & & c_{h}(\boldsymbol{x})=0
\end{array}\end{equation}
where $f_m(\boldsymbol{x},\delta(\boldsymbol{x}))$ represents the multiple objective functions to be minimized, $\boldsymbol{x}$ the design vector consisting of $n_{x_c}$ continuous and $n_{x_d}$ discrete design variables, and $\delta(\boldsymbol{x})$ the activeness function representing design space hierarchy. Continuous design variable $x_{c,i}$ bounds are represented by $\underline{x}_{c,i}$ and $\overline{x}_{c,i}$, and $N_j$ is the number of discrete options for the discrete variable $x_{d,j}$.
The inequality design constraint $k$ is given by $g_k(\boldsymbol{x},\delta(\boldsymbol{x}))$. The unrelaxable value and hidden constraints are defined by the functions $c_{v,l}(\boldsymbol{x})$ and $c_h(\boldsymbol{x})$, respectively, both returning $1$ if the constraint is violated and $0$ otherwise.
Considering all types of constraints, the design points and regions can have several different non-exclusive statuses as explained in the previous sections and summarized in Table~\ref{tab:x_cond}.

\begin{table}
\caption{Possible statuses of design points and regions. In order for any of the conditions to be met, conditions and operations above have to be met and performed as well.}\label{tab:x_cond}
\begin{tabular}{llll}
\toprule
& \multicolumn{2}{c}{Status if condition is ...} & Performed \\
Condition & met & not met & operations \\
\midrule
Any $\boldsymbol{x}$ in Cartesian product of $x_d$ values & Declared & & \\
All value constraints satisfied ($c_{v,l}(\boldsymbol{x})=0$) & Correct & Invalid & Correction \\
$\boldsymbol{x}$ is canonical & Valid & Non-canonical & Imputation \\
Hidden constraint satisfied ($c_{h}(\boldsymbol{x})=0$) & Viable & Failed & Evaluation \\
All design constraint satisfied ($g_{k}(\boldsymbol{x})\leq0$) & Feasible & Infeasible & Optimization \\
Optimality achieved & Optimal & Non-optimal & Optimization \\
\bottomrule
\end{tabular}
\end{table}

\section{Algorithms for System Architecture Optimization}\label{sec:ArchOptAlgos}

This section discusses optimization algorithms for solving SAO problems.
First, appropriate classes of algorithms are discussed in Section~\ref{subsec:Algos}: Multi-Objective Evolutionary Algorithms (MOEA) and Bayesian Optimization (BO) algorithms.
Section~\ref{subsec:NonHierBO} then provides details of the BO algorithm used in this work not related to design space hierarchy.
The implementation of optimization algorithms in the \noun{SBArchOpt} library is discussed in Section~\ref{subsec:SBArchOpt}.

\subsection{Appropriate Algorithm Classes}\label{subsec:Algos}

The  properties of SAO problems discussed in Section~\ref{sec:ArchOpt} have several consequences, also summarized in Table~\ref{tab:conseq}, on the types of optimization algorithms that can be used to solve these problems.
The fact that the objective and design constraint functions exhibit non-linear, non-smooth, and multi-modal behavior without gradient availability necessitates using a global, gradient-free optimization algorithm~\cite{Martins2022}.
Gradient-based optimization is additionally excluded from consideration due to the mixed-discrete, hierarchical nature of the design space.
Hierarchy induces the need to move design points from the declared to the valid design space by applying correction and imputation~\cite{Zaefferer2018a}.
The expensive nature of the evaluation functions drives the need to minimize the number of evaluations needed to find the optimum Pareto-set~\cite{Jones1998}.
Finally, constraint handling is needed for both the design constraints~\cite{Martins2022} and hidden constraints~\cite{Mueller2019}.

\begin{table}
\caption{Behavior of architecture optimization problems and needed algorithm capabilities.}\label{tab:conseq}
\begin{tabular}{llll}
\toprule
Space & Property & Capability needs \\
\midrule
\multirow{2}{*}{Design space} & Mixed-discrete & Gradient-free optimization~\cite{Martins2022} \\
 & Hierarchical & Correction and imputation~\cite{Zaefferer2018a} \\
\midrule
\multirow{5}{*}{Solution space} & Black-box & Global, gradient-free optimization~\cite{Martins2022} \\
 & Expensive & Minimize number of function evaluations~\cite{Jones1998} \\
 & Multi-objective & Find the Pareto-set~\cite{Miettinen1998} \\
 & Design constraints & Constraint handling~\cite{Martins2022} \\
 & Hidden constraints & Failed area avoidance~\cite{Mueller2019} \\
\bottomrule
\end{tabular}
\end{table}

Many classes of global optimization algorithms exist, designed for various purposes~\cite{Martins2022}. Exact optimization methods such as grid search~\cite{Yang2020} or DIRECT~\cite{Jones2020} are designed to yield reproducible results with provable convergence behavior, and are relevant for problems where finding the true optimum has a high priority~\cite{Locatelli2021}, however exact methods struggle to solve high-dimensional problems~\cite{Jones2020}.
On the other hand, heuristic optimization methods depend on strategies that work well in practice to search a design space more efficiently than exact optimization methods, without providing mathematical proof of convergence~\cite{Glover2003}.
One of the most powerful classes of heuristic global optimization methods are \textit{evolutionary algorithms} (EA). Evolutionary algorithms are population-based algorithms that mimic evolutionary processes found in nature: an initial population of individuals (design points) is evolved (optimized) for maximum fitness (objective value) by generating offspring (new design points) from selected parents whose genes (portions of the design vector) are crossed-over and mutated~\cite{Petrowski2017}.
Major variations between evolutionary algorithms lie in the way design points are encoded (\textit{i.e.} the encoding grammar~\cite{Selva2012}), how parents are selected for the basis of the new generation, types of cross-over and mutation operators used, and how the new generation is created from the current generation and offspring by a survival operator.
Evolutionary algorithms are robust in dealing with mixed-discrete design spaces: dealing with continuous variables was in fact a development that came later, for example through Differential Evolution~\cite{Petrowski2017} or CMA-ES~\cite{Hamano2022}.
Design constraints are handled either by applying a penalty on the objective functions or by integrating constraints in the selection and survival operators directly~\cite{Petrowski2017}.
Multi-Objective Evolutionary Algorithms (MOEA) have been developed to deal with multi-objective problems by modifying selection and survival operators for searching for a Pareto-front rather than a single optimal point~\cite{Petrowski2017}. One of the most popular MOEA is the Non-dominated Sorting Genetic Algorithm II (NSGA-II)~\cite{Deb2002} due to its low configuration effort and good performance.
Finally, EA can explore hierarchical design spaces using the hidden genes approach~\cite{Abdelkhalik2012,Nyew2015,Talbi2024} and deal with hidden constraints using the extreme barrier approach~\cite{Gratton2014}.
The appropriateness of evolutionary algorithms in general and NSGA-II in particular for SAO is noted by \noun{Crawley \textit{et al.}}~\cite{Crawley2015} and has been demonstrated in the past by various applications, see Table~\ref{tab:MOEA-app}.
Evolutionary algorithms have also found application for optimizing software product lines~\cite{LopezHerrejon2015}, which involves optimization problems characterized by strong hierarchy and choice dependencies, similarly to architecture optimization.

\begin{table}
\caption{Past applications of evolutionary algorithms to architecture optimization problems.
Nomenclature: number of discrete design points $n_{\mathrm{valid,discr}}$, number of continuous dimensions $n_{x_c}$, number of objectives $n_f$, number of function evaluations $N_{\mathrm{fe}}$, Ant Colony Optimization (ACO), Genetic Algorithm (GA), Non-dominated Sorting GA II (NSGA-II), Multi-Objective Evolutionary Algorithms (MOEA),  Guidance, Navigation and Control (GN\&C).}\label{tab:MOEA-app}
\begin{tabular}{llllllll}
\toprule
& & \multicolumn{2}{c}{Design space} & & & \\
Algorithm(s) & Application & $n_{\mathrm{valid,discr}}$ & $n_{x_c}$ & $n_f$& $n_g$ & $N_{\mathrm{fe}}$ & Reference \\
\midrule
GA & Supersonic aircraft & 576 & 100 & 10 & & 200 000 &~\cite{Buonanno2005} \\
NSGA-II & Aircraft engine & 1 163 & 30 & 3 & 15 & 4 000&~\cite{Bussemaker2021c} \\
ACO + GA & Aircraft subsystem & 9 600 & & 2 & & 96 000 &~\cite{Judt2016} \\
NSGA-II & Aircraft family & 21 875 & 4 & 4 & & 200 000 &~\cite{Pate2012} \\
NSGA-II & Launch vehicle & 123 000 & 23 & 3 & & 50 000 &~\cite{Frank2016} \\
NSGA-II & Space transport & 49e6 & 30 & 2 & 4 & 50 000 &~\cite{Frank2018} \\
MOEA & GN\&C system & 79e6 & & 2 & & 1 000&~\cite{Apaza2021} \\
(\footnotemark[1]) & Satellite instruments & 8.8e12 & & 4 & & 20 000 &~\cite{Selva2012} \\
\bottomrule
\end{tabular}
\footnotetext[1]{A special-purpose architecture optimization evolutionary algorithm was used.}
\end{table}

\begin{figure}
\centering
\includesvg[inkscapelatex=false,width=\textwidth]{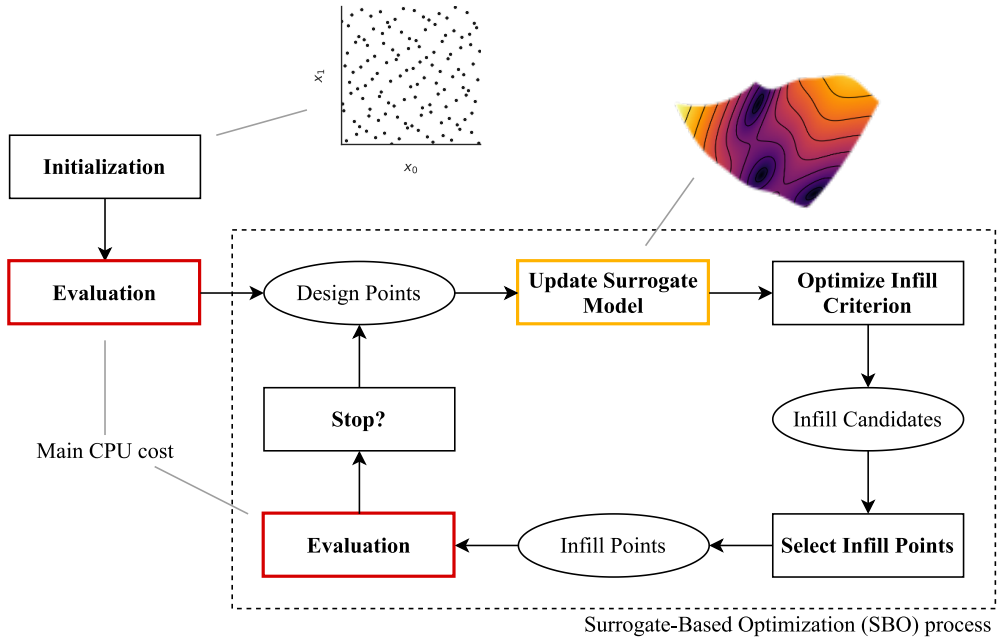}
\caption{Principle of Surrogate-Based Optimization (SBO), adopted from~\cite{Bussemaker2021}}\label{fig:sbo}
\end{figure}

As powerful as evolutionary algorithms are for dealing with the challenges of architecture optimization, there is one aspect that poses a problem if the objective and constraint functions are expensive to evaluate: the high number of needed function evaluations $N_{\mathrm{fe}}$~\cite{Chugh2019}. As can be seen in Table~\ref{tab:MOEA-app}, $N_{\mathrm{fe}}$ typically is in the order of thousands to hundreds of thousands for EA.
Each evaluation of the aircraft engine problem of \noun{Bussemaker \textit{et al.}}~\cite{Bussemaker2021c} took about two minutes, resulting in a little over 5 days to complete the 4000 evaluations. Depending on the context this could still be acceptable, however the time required to solve an optimization problem driven by an evolutionary algorithm quickly becomes intractable if more evaluations are needed, or if evaluation becomes more costly if more or higher fidelity analyses are used for architecture evaluation.
To solve such problems, \textit{Surrogate-based Optimization (SBO)} algorithms have been developed~\cite{Martins2022}. SBO algorithms rely on some surrogate model, also known as response surface model or metamodel, of the objective and constraint functions in order to efficiently determine interesting design point(s) to evaluate, the \textit{infill} point(s). Once the infill points have been evaluated using the expensive evaluation functions, the surrogate model is reconstructed and the process starts over, until some termination criterion has been reached.
Infill points are selected using infill criteria, also known as acquisition functions, that are normally defined to represent some kind of trade-off between exploration (finding new interesting regions) and exploitation (improving already interesting regions) of the design space. Figure~\ref{fig:sbo} visualizes the working principle of SBO.

The main choices involved in using SBO are the initial sampling method for creating the first surrogate model, the type of surrogate model, the infill criterion, and the termination criterion.
For example, a Radial Basis Function (RBF) surrogate has been combined with a Coordinate Perturbation infill criterion, striking a balance between the best predicted objective value (exploitation) and staying away from previously evaluated points (exploration)~\cite{Regis2013,Bagheri2017}.
SBO has also been extended with strategies for dealing with hidden constraints~\cite{Mueller2019}.
Especially powerful types of SBO algorithms are \textit{Bayesian Optimization (BO)} algorithms~\cite{Garnett2023}. BO algorithms use surrogate models that provide error estimates $\hat{s}(\boldsymbol{x})$ in addition to function estimates $\hat{y}(\boldsymbol{x})$, and therefore can use infill criteria that leverage this additional information in order to better predict where interesting points lie. 
Most applications of BO use Gaussian Process (GP) models~\cite{Chugh2019}, also known as Kriging models.
Extensive research has been performed on handling design constraints~\cite{Schonlau1998,Sasena2002} and solving multi-objective~\cite{Knowles2006,RojasGonzalez2019} problems using BO.
GP models were originally developed for continuous variables, however research has been performed into methods for supporting mixed-discrete~\cite{GarridoMerchan2020,Daulton2022,Pelamatti2019,Munoz2020,Saves2023a,Dreczkowski2023} and hierarchical~\cite{Pelamatti2020a,Audet2022,Zaefferer2018a,Horn2019,Lu2018,Jenatton2017,Saves2023SMT} design spaces.
Enabling the use of GP models for high-dimensional design spaces is an active area of research~\cite{Bouhlel2018,Garnett2023,Priem2023,Saves2024}.
For the especially non-smooth and hierarchical design spaces encountered in hyperparameter optimization~\cite{Bischl2023} problems Random Forest (RF)~\cite{Hutter2011,Lindauer2022} and Tree Parzen Estimator (TPE)~\cite{Bergstra2011,Ozaki2022} models have been applied. These two models naturally represent the tree structure of hierarchical design spaces~\cite{Eggensperger2015}, however recent evidence shows the superior performance of mixed-discrete GP models over RF and TPE models~\cite{GarridoMerchan2020}, and therefore in this work we will concentrate on the use of BO with GP models only.

In the rest of this work, both Multi-Objective Evolutionary Algorithms (MOEA) and Bayesian Optimization (BO) algorithms will be further investigated for use in architecture optimization problems. One is not better than the other: MOEA should be used if evaluation is not expensive (\textit{e.g.} one evaluation takes at most in the order of seconds), and BO should be used if evaluation is expensive (\textit{e.g.} in the order of minutes or more)~\cite{gamot2023hidden}. Bayesian Optimization should not be used for inexpensive problems, as then time for model fitting and searching for infill points becomes limiting, rather than function evaluation time~\cite{Chugh2019}.

\subsection{Non-hierarchical Foundation of the Bayesian Optimization Algorithm}\label{subsec:NonHierBO}

An advantage of BO is the high level of composability: many specializations of BO (\textit{e.g.} mixed-discrete, multi-objective) can be combined without negative interactions~\cite{Greenhill2020}. In this work we follow the same approach: specializations, \textit{e.g.} hierarchy and hidden constraints, build on prior specializations, \textit{e.g.} mixed-discrete and multi-objective.
The GP model used by the BO algorithm further developed in this work was already able to deal with mixed-discrete variables, including categorical variables, for non-hierarchical design spaces.
Categorical variables pose a particular challenge to GP models, since there is no inherent ordering in the possible values and therefore conventional distance measures used for modeling correlation for continuous and integer variables might not be able to model the variable accurately.
The approach for modeling categorical variables is based on the work by \noun{Saves \textit{et al.}}~\cite{Saves2023a}, which uses dedicated kernels to model categorical variables, for example the Gower or Exponential Homoscedastic Hypersphere (EHH) kernel. We use the implementation in the Surrogate Modeling Toolbox (SMT)\footnote{\label{note:smt}\href{https://smt.readthedocs.io/}{https://smt.readthedocs.io/}}~\cite{Saves2023SMT}.
By default the Gower distance kernel is used.

One of the main limitations of BO is that GP models do not handle high-dimensional design spaces well. Due to the number of hyperparameters to tune, training and sampling time quickly increases with the number of input dimensions and training points~\cite{Garnett2023}. However, a feature space can be defined that reduces the dimension of the original input space~\cite{Calandra2016}, supported by the observation that in some cases the optimization problem has a lower intrinsic dimensionality.
The BO algorithm uses Kriging with Partial Least Squares (KPLS) to construct such feature spaces as described by \noun{Bouhlel \textit{et al.}}~\cite{Bouhlel2016}. Recently KPLS has been extended to work with discrete variables too~\cite{Saves2022b,Charayron2023,Saves2024} and these methods are implemented in SMT~\cite{Bouhlel2019,Saves2023SMT}.
By default, KPLS is applied with $n_{\mathrm{kpls}} = 10$ if the design space contains more than 10 design variables.

Infill points are selected by formulating an infill ensemble: an approach first presented by \noun{Lyu \textit{et al.}}~\cite{Lyu2018} based on the observation that different infill criteria can suggest very different infill points. This can be especially true in the earlier phases of the optimization when the GP is less accurate.
There are two advantages to this ensemble-infill strategy~\cite{CowenRivers2020}: the selected infill points represent the best trade-off between the infill criteria, and it is easy to select multiple infill points per iteration for batch optimization~\cite{Garnett2023} without needing to retrain the underlying GP models as is needed for qEI~\cite{Ginsbourger2010}.
For single-objective optimization, an ensemble of Lower Confidence Bound (LCB)~\cite{Cox1992}, Expected Improvement (EI)~\cite{Jones1998} and Probability of Improvement (PoI)~\cite{Hawe2007} infills is used.
For multi-objective, the ensemble consists of the Minimum Probability of Improvement (MPoI)~\cite{Rahat2017} and Minimum Euclidean PoI (MEPoI)~\cite{Bussemaker2021} infills.
The infill batch size $n_{\mathrm{batch}}$ is set to the maximum amount of designs that can be evaluated in parallel: $n_{\mathrm{batch}} = n_{\mathrm{parallel}}$.

Design constraints are handled by constraint function mean prediction $\hat{g}$~\cite{Sasena2002} by default, or by Probability of Feasibility (PoF)~\cite{Sohst2021} if a more or less conservative (achieved by PoF above or below 50\%, respectively) criterion is requested.
Hidden constraints are handled by training an additional model to predict the Probability of Viability (PoV), and constraining it to be at least 25\% (by default) during infill search~\cite{Bussemaker2024hc}. The PoV prediction model is the same mixed-discrete GP as used for predicting $f$ and $g$.

\subsection{Implementation in SBArchOpt}\label{subsec:SBArchOpt}

The optimization algorithms presented in this work are implemented in the \noun{SBArchOpt}\footnote{\href{https://sbarchopt.readthedocs.io/}{https://sbarchopt.readthedocs.io/}} library: an open-source Python library for solving SAO problems~\cite{Bussemaker2023a}.
\noun{SBArchOpt} features a problem definition class with interfaces for executing hierarchical optimization problems.
The problem definition class is built on top of \noun{pymoo}'s\footnote{\href{https://pymoo.org/}{https://pymoo.org/}} \noun{Problem} class~\cite{Blank2020}, and adds functions for performing imputation and correction (standalone, not part of design point evaluation), getting information about which design variables are conditionally active, generating all valid discrete design vectors $x_{\mathrm{valid,discr}}$, and for getting various statistics (see Section~\ref{subsec:DSchallenges}).
The hierarchical structure of the problem can either be provided implicitly by implementing the correction and imputation function, or it can be modeled explicitly using \noun{ConfigSpace}\footnote{\href{https://automl.github.io/ConfigSpace/}{https://automl.github.io/ConfigSpace/}}~\cite{Lindauer2019}.
\noun{ConfigSpace} enables declaring design variables, specifying activation conditions between them, and specifying value constraints. For more discussion about the role \noun{ConfigSpace} can play in hierarchical optimization, refer to Section~\ref{sec:HierIntegration}.

\noun{SBArchOpt} implements various optimization algorithms for solving architecture optimization problems.
These optimization algorithms are implemented using \noun{pymoo} classes, enabling provisioning of any \noun{pymoo} algorithm for solving architecture optimization problems.
For example, \noun{SBArchOpt} provides the hierarchical sampling algorithm presented in Section~\ref{sec:HierSampling} as a \noun{pymoo} class. Similarly, a repair operator is implemented that either uses the problem-agnostic correction algorithm of Section~\ref{sec:HierSampling}, or uses the correction function supplied by the problem definition class.
A problem-specific correction function can correct both discrete and continuous design variables.
Additionally, \noun{SBArchOpt} implements result storage and restart functionalities, which are important for expensive and therefore long-running optimization problems.
Two pre-configured \noun{pymoo} algorithms are provided: \noun{ArchOptNSGA2} implementing NSGA-II~\cite{Deb2002} for architecture optimization, and \noun{DOEAlgorithm} that only executes the hierarchical sampling algorithm to perform a Design of Experiments (DoE).
The SBO algorithm developed in this work is implemented as \noun{ArchSBO}, and features automatic selection of infill and hidden constraints strategies. It can either use an RBF model or a GP model, both using the implementation in the Surrogate Modeling Toolbox (SMT)\footnotemark[\getrefnumber{note:smt}]~\cite{Saves2023SMT}.
To test and promote solving architecture optimization problems with SBO in general, \noun{SBArchOpt} also implements connections to other SBO libraries such as \noun{BoTorch}~\cite{Balandat2019}, \noun{Trieste}~\cite{Picheny2023}, \noun{HEBO}~\cite{CowenRivers2020}, \noun{SEGOMOE}~\cite{Bartoli2019} and \noun{SMARTy}~\cite{Bekemeyer2022}.

Finally, to support SAO algorithm development \noun{SBArchOpt} contains a library of test problems featuring various combinations of continuous or mixed-discrete, (non-)hierarchical, single-objective or multi-objective, (un)constrained problems with or without hidden constraints.
Many problems are implemented from literature or adopted from \noun{pymoo}'s test problem database, however also several new test problems are provided.
Realistic engineering problem behavior is exhibited by the jet engine architecture optimization problem used in Section~\ref{sec:JetEngineOpt}, a multi-stage launch vehicle architecture problem adopted from~\cite{GarciaSanchez2024} and several versions of a Guidance, Navigation \& Control (GNC) problem adopted from~\cite{Crawley2015,Apaza2021}, see also Section~\ref{sec:HierTestProblems}.

\section{Modeling Hierarchical Design Spaces in Bayesian Optimization}\label{sec:HierModel}

In order to apply SBO to architecture optimization problems, the surrogate model should be able to accurately model the behavior of the objective and constraint functions over the design space. For architecture optimization problems it means that the hierarchical structure between design variables should be integrated in the structure of the surrogate model.
Gaussian Process (GP) models as used in Bayesian Optimization (BO) have been extended to support hierarchical design variables by \noun{Saves \textit{et al.}}~\cite{Saves2023SMT} and implemented in version 2.0 of the Surrogate Modeling Toolbox (SMT)\footnotemark[\getrefnumber{note:smt}]. 
In previous works about hierarchical design spaces~\cite{HalleHannan2024,Lucas_b,Pelamatti2020a,Hutter2013,HuJoMe2009}, variables are structured according to their types and roles, which reflect the complex nature of real-world problems. Variables can be classified into three main categories: neutral, meta, and decreed. 
\begin{itemize}
    \item \textbf{Neutral Variables}: These variables not affected by other variables in the hierarchy. They can be continuous, integer, or categorical.
    \item \textbf{Meta Variables}: These variables determine whether other variables are activated or not, influencing the structure of the problem.
    They are also known as dimensional variables~\cite{Pelamatti2020a}, tree-structured variables~\cite{Bergstra2011} or conditional search space variables~\cite{Abdelkhalik2012}.
    \item \textbf{Decreed Variables}: These are the dependant conditionally active variables that are activated based on the values of the meta variables.
\end{itemize}
Modeling these hierarchical design spaces is critical in Bayesian optimization, particularly for complex problems such as system architecture optimization problems in this work.
Traditional GP and kernel functions often fall short in handling such complexity. Therefore, recent works have focused on developing kernels to accommodate hierarchical structures. One such innovation is the Alg-Kernel, which is designed to handle the intricacies of hierarchical and mixed-variable problems. The Alg-Kernel effectively captures the dependencies and interactions among different types of variables, providing a robust framework for surrogate modeling in hierarchical design spaces~\cite{Saves2023SMT}.
This work, however, only considered continuous and integer hierarchical variables. Here we provide an extension to previous work that also supports categorical conditionally active variables.
Based on the advances published in~\cite{Saves2023a}, we can define a new hierarchical kernel for categorical conditional variables based on one-hot encoding and on the algebraic distance defined in~\cite{Saves2023SMT}. 

For two given inputs, for example, the $r^{\text{th}}$ and  $s^{\text{th}}$ points in the DoE ($x_r$ and $x_s$), let $c^r_{i} $ and $c^s_{i} $ be the associated categorical variables taking respectively the $\ell^i_r$ and the $\ell^i_s$ level on the categorical variable $c_i$. Denote $e_{c^r_i}$ the one-hot encoding~\cite{one-hot} of $c^r_i$ that takes value $0$ everywhere but on the dimension $\ell_r^i$:  $e_{c^r_i} \in \mathbb{R}^{L_i} $ such that $\left( e_{c^r_i} \right)_{\ell^i_r}=1$ and $\left( e_{c^r_i} \right)_k=0$ for $k \neq \ell^i_r$.
The new hierarchical kernel between two points $r$ and $s$ and for categorical variables $x_i$ is formulated using similar principles as in~\cite{Saves2023a}.
Let $[R_i(\Theta_i)]$ be the  GP correlation kernel for the $i^{\text{th}}$ categorical variable $c_i$ that depends on $\Theta_i$, the hyperparameter matrix related to $c_i$. The variable $c_i$ features $L_i$ different levels and $[\Phi(\Theta_i)]$ is a chosen symmetric positive definite parameterization of the matrix $\Theta_i$.

\begin{itemize}
    \item If $\delta_i\left(x_r\right) = \delta_i\left(x_s\right) = 0$ (both variables are inactive): all the one-hot relaxed dimensions associated to this variable are also inactive meaning none of them is relevant.
    \item If $\delta_i\left(x_r\right) = \delta_i\left(x_s\right) = 1$ (both variables are active): all the one-hot relaxed dimensions associated to this variable are also active. 
    If $\ell^i_r = \ell^i_s, [R_i(\Theta_i)]_{{\ell^i_r},{\ell^i_r}} =1$. Otherwise, assuming that the chosen kernel is the exponential kernel,
    this leads to the  kernel
    \begin{equation}
        [R_i(\Theta_i)]_{{\ell^i_r},{\ell^i_s}} = \exp \left( -  \sqrt{2} [\Phi(\Theta_i)]_{{\ell^i_r},{\ell^i_r}}- \sqrt{2} [\Phi(\Theta_i)]_{{\ell^i_s},{\ell^i_s}}  \right)
    \end{equation}
    \item If $\delta_i\left(x_r\right) \neq \delta_i\left(x_s\right)$ (only one of the variables is active): for all the relaxed dimensions associated to this variable, because of the algebraic distance, there is a induced distance $1$ between both inputs.
    Assuming that the chosen kernel is the exponential kernel, this leads to the  kernel 
    \begin{equation}
        [R_i(\Theta_i)]_{{\ell^i_r},{\ell^i_s}} = \exp \left( -  \displaystyle{ \sum^{L_i}_{j=1} \    [\Phi(\Theta_i)]_{j,j}}  \right)
    \end{equation}
\end{itemize}
As shown in~\cite{Saves2022}, continuous relaxation (one-hot encoding) generalizes the Gower distance. Therefore, denoting $\theta_{\text{cov}}$ the unique correlation parameter for the considered variable, we can define the following adaptation for the particular case in which we are using the Gower distance kernel for a categorical variable $c_i$ with $L_i$ levels:
\begin{equation}
    \label{eq:gd_hier}    
\begin{aligned}
& \bullet d_i(x_r,x_s) = 0 \text{ if both } x_r \text{ and } x_s \text{ are inactive. } \\
& \bullet  d_i(x_r,x_s) = \sqrt{2} \ \theta_{\text{cov}} \text{ if both } x_r \text{ and } x_s \text{ are active.}\\
& \bullet  d_i(x_r,x_s) = \frac{1}{2} L_i \ \theta_{\text{cov}} \text{ if either } x_r \text{ or } x_s \text{ are active. } \\
\end{aligned}
\end{equation}
Finally, for the exponential homoscedastic hypersphere or for the homoscedastic hypersphere  kernel we apply the imputation method: whenever one of the variables is inactive, it is treated as if its value would be the first available level $L_0$.
These new categorical hierarchical kernels are implemented in SMT 2.3\footnotemark[\getrefnumber{note:smt}]. The rest of this work uses the Gower distance categorical hierarchical kernel for conditionally active variables defined in Equation~\eqref{eq:gd_hier}.

We now discuss other stages of an optimization run where information about design space hierarchy can be leveraged.

\section{Hierarchical Sampling and Correction Algorithms}\label{sec:HierSampling}

In this section we investigate how to most effectively generate design vectors covering the design space in such a way as to obtain as much information about the design space behavior as possible for a given function evaluation budget. This is relevant both for sampling the design space and when correcting invalid design vectors.
When sampling the design space (creating the DoE) for use as initial database of design points for SBO algorithms or initial population for evolutionary algorithms, the effects of rate diversity should be considered: if this effect is ignored, some regions in the design space may be over- or under-represented~\cite{Crawley2015}.
Correcting invalid design vectors might affect design space exploration by modifying results of evolutionary operators for evolutionary algorithms and infill optimization for SBO algorithms. Correction behavior is especially important for problems with high correction ratios due to the low chance of randomly finding a valid design vector.
First, we introduce the test problems used to investigate sampling and correction algorithms in Section~\ref{sec:HierTestProblems}. Then, sampling algorithms are investigated and selected in Section~\ref{sec:HierSamplingExp}, and correction algorithms in Section~\ref{sec:HierCorrExp}.

\subsection{Test Problems}\label{sec:HierTestProblems}

To compare correction and sampling strategies, the following test problems implemented in \noun{SBArchOpt} are used (see Table~\ref{tab:hier_probs} for detailed statistics):
\begin{itemize}
    \item Three instances of a multi-stage rocket design problem ("Rocket"), adopted from \noun{García Sánchez}~\cite{GarciaSanchez2024}. The problem attempts to minimize cost while maximizing payload mass for a given target orbit altitude, by selecting the number of stages, engine types and amounts per stage, and geometrical parameters of stages and rocket head. The problem has been modified to be less tightly constrained than the original problem.
    In addition to this multi-objective problem ("Rocket"), we also use two single-objective versions: one that minimizes cost ("RCost") and one that minimizes a weighted function of cost and payload mass ("RWt").
    The cost minimum lies in the group of single-stage rockets, which only makes up 0.3\% of all valid discrete design vectors, whereas the weighted minimum lies in one of the much larger multi-stage rocket groups.
    \item Four instances of a Guidance, Navigation \& Control (GNC) problem, adopted from~\cite{Crawley2015,Apaza2021}, which features mass and failure rate as minimization objectives by selecting the number and type of, and connections between sensors, flight computers and actuators.
    The original problem has been modified to use continuous variables for selecting object types, to turn the problem into a mixed-discrete problem (compared to fully discrete).
    We use the mixed-discrete GNC problem version without actuators as a multi-objective test problem ("MD GNC"), a single-objective version that minimizes weight only ("Wt GNC"), a single-objective version that minimizes failure rate ("FR GNC"), and a single-objective version that solves a scalarized objective composed of weight and failure rate ("SO GNC").
    The weight minimum lies in the group containing only one sensor and one computer, which is only represented by one $x$ of $x_{\mathrm{valid,discr}}$ (1 / 327 = 0.3\%).
    \item A version of the simple jet engine problem solved in Section~\ref{sec:JetEngineOpt} that uses surrogate models for evaluation ("Jet SM").
    This version represents the same optimization problem, however replaces the multidisciplinary thermodynamic cycle analysis evaluation by random forest regressors for each output, to reduce evaluation times (milliseconds, compared to minutes for the original problem).
    We use the random forest regressors implemented in \noun{scikit-learn}\footnote{\href{https://scikit-learn.org/}{https://scikit-learn.org/}}.
    Just as the original problem, this version features a hidden constraint (\textit{i.e.} failed evaluations).
\end{itemize}
Table~\ref{tab:hier_probs} provides some more statistics: both single-objective and multi-objective problems are included, and all except the GNC problems contain design constraints.
The GNC problems feature a relatively high $\mathrm{IR}$, showing that they indeed feature hierarchical design spaces.
The GNC and Jet SM problems need correction algorithms as they get about half their $\mathrm{IR}$ from the need for correction shown by a $\mathrm{CRF}$ a little over 50\%.
$\mathrm{MRD}$ is a common occurrence in SAO problems too, as can be seen.

\begin{table}
\caption{Test problems used for testing hierarchical sampling and optimization strategies.}\label{tab:hier_probs}
\begin{tabular}{lllllllllll}
\toprule
Problem &\noun{SBArchOpt} Class &  $n_f$ &$n_g$& $n_{x_c}$&$n_{x_d}$  &$n_{\mathrm{valid,discr}}$&$\mathrm{IR}$&$\mathrm{CR}$ & $\mathrm{CRF}$ &$\mathrm{MRD}$\\
\midrule
 RCost &\noun{SOLCRocketArch}& 1& 2& 6& 8 &18522& 3.7 &1.0 &0\%&94\%\\
 RWt &\noun{SOLCRocketArch}& 1& 2& 6& 8 &18522& 3.7 &1.0 &0\%&94\%\\
 Rocket &\noun{LCRocketArch}& 2& 2& 6& 8 &18522& 3.7 &1.0 &0\%&94\%\\
 SO GNC &\noun{SOMDGNCNoAct}& 1& & 6& 11 &327& 150 &16.1 &56\%&88\%\\
 FR GNC &\noun{SOMDGNCNoAct}& 1& & 6& 11 &327& 150 &16.1 &56\%&88\%\\
 Wt GNC &\noun{SOMDGNCNoAct}& 1& & 6& 11 &327& 150 &16.1 &56\%&88\%\\
 MD GNC &\noun{MDGNCNoAct}& 2& & 6& 11 &327& 150 &16.1 &56\%&88\%\\
 Jet SM &\noun{SimpleTurbofanArchModel}& 1& 5& 9& 6 &70& 3.9 &2.1 &55\%&60\%\\
\end{tabular}
\end{table}

In subsequent sections, optimization performance is compared using the procedure described in Appendix~\ref{sec:perf_rank}.
A rank of 1 indicates best performance, higher values indicate progressively worse performance. The best performing infill is then selected by looking at the highest proportion of rank $\leq$ 2 and rank 1.
Performance is compared based on $\Delta \mathrm{HV}$ ($\Delta$ hypervolume) regret. $\Delta \mathrm{HV}$ represents the distance to the known optimum (or Pareto front in case of multi-objective optimization) normalized to the range of objective values.
Regret represents the cumulative $\Delta \mathrm{HV}$ as seen over the optimization run. For more details, refer to Appendix~\ref{sec:perf_rank}.

\subsection{Sampling}\label{sec:HierSamplingExp}

Many sampling methods have been developed in the past, such as random sampling or Latin Hypercube Sampling (LHS)~\cite{Martins2022}.
Such methods can still be applied for hierarchical design spaces, for example by rejecting sampled design invalid design vectors (\textit{i.e.} vectors with violated value constraints) as done in \noun{ConfigSpace}~ \cite{Lindauer2019}. Instead of rejecting invalid design vectors, it is also possible to apply correction. Using conventional sampling methods with correction however might result in over- or under-representation of certain design space regions~\cite{Crawley2015}.
Another way is to sample design vectors directly from the set of valid design vectors $x_{\mathrm{valid,discr}}$: this ensures generated design vectors are valid and therefore correction is not needed after sampling.
To do this, $x_{\mathrm{valid,discr}}$ may be divided into interesting subdivisions (\textit{e.g.} based on subproblems as defined by dimensional variables~\cite{Pelamatti2020a}) and sampling might be separated in two steps: first decide how many samples to draw from each subdivision, then sample from the subdivisions. As it is mostly the discrete design variables that decide the activeness of other design variables, in this work the hierarchical sampling procedure is applied to the discrete design vectors first. After sampling discrete design vectors, continuous design variables are filled by Sobol' sampling~\cite{Renardy2021}, ignoring inactive continuous design variables.

One way of defining subdivisions is by using design variable activeness information. For example, groups can be defined based on the number of active design variables $n_{\mathrm{act}}$, similar to what \noun{Kaltenecker \textit{et al.}}~\cite{Kaltenecker2019} did for sampling software product lines.
An extension of this is to group by the active variables $x_{\mathrm{act}}$ directly. In that case, design vectors with the same number of active variables where these variables do not completely overlap are treated as different subdivisions. This approach is similar to how \noun{Frank \textit{et al.}}~\cite{Frank2016} define their architecture optimization problems: they divide by selections of values from a morphological matrix constrained by a compatibility matrix, then for each set of selections with the same active design variables they define a separate optimization problem.
Table~\ref{tab:hier_sampling} shows how the valid discrete design vectors $x_{\mathrm{valid,discr}}$ are separated into groups based on $x_{\mathrm{act}}$.

\begin{table}
\caption{Grouping and weighting process of the hierarchical sampling algorithm: discrete design vectors are grouped by $x_{\mathrm{act}}$ and weighted by $w = n_{\mathrm{act}}$ (number of active $x$); $w_{\mathrm{rel}}$ represents the relative weighting for sampling design vectors from groups (\textit{e.g.} $w_{\mathrm{rel}} = 36\%$ means 36\% of $x$ is sampled from that group). Red background indicates an inactive variable.}\label{tab:hier_sampling}
\begin{tabular}{lllllccc}
\toprule
  $x_0$ & $x_1$ & $x_2$ & $x_3$ & $x_4$  & $x_{\mathrm{act}}$&$w = n_{\mathrm{act}}$&$w_{\mathrm{rel}}$ \\
\midrule
  {\cellcolor[HTML]{E1F3DC}} \color[HTML]{000000} 0 & {\cellcolor[HTML]{E1F3DC}} \color[HTML]{000000} 0 & {\cellcolor[HTML]{E1F3DC}} \color[HTML]{000000} 0 & {\cellcolor[HTML]{FFC7CE}} \color[HTML]{9C0006}   & {\cellcolor[HTML]{E1F3DC}} \color[HTML]{000000} 0  & \multirow{3}{*}{\begin{tabular}{c}$x_0$, $x_1$ \\ $x_2$, $x_4$\end{tabular}}&\multirow{3}{*}{4} &\multirow{3}{*}{36\%} \\
  {\cellcolor[HTML]{E1F3DC}} \color[HTML]{000000} 0 & {\cellcolor[HTML]{E1F3DC}} \color[HTML]{000000} 0 & {\cellcolor[HTML]{E1F3DC}} \color[HTML]{000000} 0 & {\cellcolor[HTML]{FFC7CE}} \color[HTML]{9C0006}   & {\cellcolor[HTML]{8ED08B}} \color[HTML]{000000} 1  &  && \\
  {\cellcolor[HTML]{E1F3DC}} \color[HTML]{000000} 0 & {\cellcolor[HTML]{E1F3DC}} \color[HTML]{000000} 0 & {\cellcolor[HTML]{E1F3DC}} \color[HTML]{000000} 0 & {\cellcolor[HTML]{FFC7CE}} \color[HTML]{9C0006}   & {\cellcolor[HTML]{2C944C}} \color[HTML]{F1F1F1} 2  &  && \\
 \midrule
  {\cellcolor[HTML]{E1F3DC}} \color[HTML]{000000} 0 & {\cellcolor[HTML]{E1F3DC}} \color[HTML]{000000} 0 & {\cellcolor[HTML]{8ED08B}} \color[HTML]{000000} 1 & {\cellcolor[HTML]{E1F3DC}} \color[HTML]{000000} 0 & {\cellcolor[HTML]{FFC7CE}} \color[HTML]{9C0006}    & \multirow{4}{*}{\begin{tabular}{c}$x_0$, $x_1$ \\ $x_2$, $x_3$\end{tabular}}&\multirow{4}{*}{4} &\multirow{4}{*}{36\%}\\
  {\cellcolor[HTML]{E1F3DC}} \color[HTML]{000000} 0 & {\cellcolor[HTML]{E1F3DC}} \color[HTML]{000000} 0 & {\cellcolor[HTML]{8ED08B}} \color[HTML]{000000} 1 & {\cellcolor[HTML]{8ED08B}} \color[HTML]{000000} 1 & {\cellcolor[HTML]{FFC7CE}} \color[HTML]{9C0006}    &  && \\
  {\cellcolor[HTML]{E1F3DC}} \color[HTML]{000000} 0 & {\cellcolor[HTML]{E1F3DC}} \color[HTML]{000000} 0 & {\cellcolor[HTML]{2C944C}} \color[HTML]{F1F1F1} 2 & {\cellcolor[HTML]{E1F3DC}} \color[HTML]{000000} 0 & {\cellcolor[HTML]{FFC7CE}} \color[HTML]{9C0006}    &  && \\
  {\cellcolor[HTML]{E1F3DC}} \color[HTML]{000000} 0 & {\cellcolor[HTML]{E1F3DC}} \color[HTML]{000000} 0 & {\cellcolor[HTML]{2C944C}} \color[HTML]{F1F1F1} 2 & {\cellcolor[HTML]{8ED08B}} \color[HTML]{000000} 1 & {\cellcolor[HTML]{FFC7CE}} \color[HTML]{9C0006}    &  && \\
 \midrule
  {\cellcolor[HTML]{E1F3DC}} \color[HTML]{000000} 0 & {\cellcolor[HTML]{8ED08B}} \color[HTML]{000000} 1 & {\cellcolor[HTML]{FFC7CE}} \color[HTML]{9C0006}   & {\cellcolor[HTML]{FFC7CE}} \color[HTML]{9C0006}   & {\cellcolor[HTML]{FFC7CE}} \color[HTML]{9C0006}    & $x_0$, $x_1$&2 &18\% \\
 \midrule
  {\cellcolor[HTML]{8ED08B}} \color[HTML]{000000} 1 & {\cellcolor[HTML]{FFC7CE}} \color[HTML]{9C0006}   & {\cellcolor[HTML]{FFC7CE}} \color[HTML]{9C0006}   & {\cellcolor[HTML]{FFC7CE}} \color[HTML]{9C0006}   & {\cellcolor[HTML]{FFC7CE}} \color[HTML]{9C0006}    & $x_0$&1 &9\% \\
\bottomrule
\end{tabular}
\end{table}

Another way is to define subdivisions based on the value of certain discrete design variables.
\noun{Pelamatti \textit{et al.}}~\cite{Pelamatti2020a} defined subdivisions based on dedicated dimensional variables. This approach, however, requires the user to define which variables act as such grouping variables, thereby requiring the user to think in terms of the design variables rather than the architecture design space (\textit{i.e.} domain-specific design space model).
The goal of hierarchical sampling is to mitigate max rate diversity ($\mathrm{MRD}$) issues: subdivisions can therefore be defined based on values of high-$\mathrm{RD}$ variables.
Subdivisions are made by recursively selecting $\mathrm{argmax} \mathrm{RD}_j$, subject to $\mathrm{RD}_j \geq \mathrm{RD}_{\mathrm{min}}$ where $\mathrm{RD}$ is determined for the current (partial) subdivision and $\mathrm{RD}_{\mathrm{min}}$ is set such that only high-$\mathrm{RD}$ variables are used for grouping.
Table~\ref{tab:mrd_sampling} shows an example of grouping by $\mathrm{MRD}$ with $\mathrm{RD}_{min} = 50\%$: group 3 is defined based on $\mathrm{RD}_0$ ($\mathrm{RD}$ of $x_0$), group 1 and 2 are defined based on $\mathrm{RD}_1$. No further groups are defined because the next-highest $\mathrm{RD}$ is $\mathrm{RD}_2$, which is lower than $\mathrm{RD}_{min}$.
In this work we use $\mathrm{RD}_{min} = 80\%$ as it is a good compromise between reducing group-internal $\mathrm{MRD}$ and the number of groups created.

\begin{table}
\caption{Grouping by $\mathrm{MRD}$: discrete design vectors are grouped by values of high-$\mathrm{RD}$ variables.
The example uses $\mathrm{RD}_{min} = 50\%$.
Red background indicates an inactive variable.}\label{tab:mrd_sampling}
\begin{tabular}{cccccc}
\toprule
   &$x_0$ & $x_1$ & $x_2$ & $x_3$ & $x_4$ \\
\midrule

$x_j = 0$ & 89\% & 87\% & 42\% & 50\% & 33\% \\
$x_j = 1$ & 11\% & 13\% & 29\% & 50\% & 33\% \\
$x_j = 2$ &      &      & 29\% &      & 33\% \\
$\mathrm{RD}$    & 78\% & 74\% & 13\% & 0\%  & 0\% \\

\midrule

\multirow{7}{*}{Group 1}   &{\cellcolor[HTML]{E1F3DC}} \color[HTML]{000000} 0 & {\cellcolor[HTML]{E1F3DC}} \color[HTML]{000000} 0 & {\cellcolor[HTML]{E1F3DC}} \color[HTML]{000000} 0 & {\cellcolor[HTML]{FFC7CE}} \color[HTML]{9C0006}   & {\cellcolor[HTML]{E1F3DC}} \color[HTML]{000000} 0 \\
   &{\cellcolor[HTML]{E1F3DC}} \color[HTML]{000000} 0 & {\cellcolor[HTML]{E1F3DC}} \color[HTML]{000000} 0 & {\cellcolor[HTML]{E1F3DC}} \color[HTML]{000000} 0 & {\cellcolor[HTML]{FFC7CE}} \color[HTML]{9C0006}   & {\cellcolor[HTML]{8ED08B}} \color[HTML]{000000} 1  \\
   &{\cellcolor[HTML]{E1F3DC}} \color[HTML]{000000} 0 & {\cellcolor[HTML]{E1F3DC}} \color[HTML]{000000} 0 & {\cellcolor[HTML]{E1F3DC}} \color[HTML]{000000} 0 & {\cellcolor[HTML]{FFC7CE}} \color[HTML]{9C0006}   & {\cellcolor[HTML]{2C944C}} \color[HTML]{F1F1F1} 2  \\
   &{\cellcolor[HTML]{E1F3DC}} \color[HTML]{000000} 0 & {\cellcolor[HTML]{E1F3DC}} \color[HTML]{000000} 0 & {\cellcolor[HTML]{8ED08B}} \color[HTML]{000000} 1 & {\cellcolor[HTML]{E1F3DC}} \color[HTML]{000000} 0 & {\cellcolor[HTML]{FFC7CE}} \color[HTML]{9C0006} \\
   &{\cellcolor[HTML]{E1F3DC}} \color[HTML]{000000} 0 & {\cellcolor[HTML]{E1F3DC}} \color[HTML]{000000} 0 & {\cellcolor[HTML]{8ED08B}} \color[HTML]{000000} 1 & {\cellcolor[HTML]{8ED08B}} \color[HTML]{000000} 1 & {\cellcolor[HTML]{FFC7CE}} \color[HTML]{9C0006} \\
   &{\cellcolor[HTML]{E1F3DC}} \color[HTML]{000000} 0 & {\cellcolor[HTML]{E1F3DC}} \color[HTML]{000000} 0 & {\cellcolor[HTML]{2C944C}} \color[HTML]{F1F1F1} 2 & {\cellcolor[HTML]{E1F3DC}} \color[HTML]{000000} 0 & {\cellcolor[HTML]{FFC7CE}} \color[HTML]{9C0006} \\
   &{\cellcolor[HTML]{E1F3DC}} \color[HTML]{000000} 0 & {\cellcolor[HTML]{E1F3DC}} \color[HTML]{000000} 0 & {\cellcolor[HTML]{2C944C}} \color[HTML]{F1F1F1} 2 & {\cellcolor[HTML]{8ED08B}} \color[HTML]{000000} 1 & {\cellcolor[HTML]{FFC7CE}} \color[HTML]{9C0006} \\
 \midrule
 Group 2  &{\cellcolor[HTML]{E1F3DC}} \color[HTML]{000000} 0 & {\cellcolor[HTML]{8ED08B}} \color[HTML]{000000} 1 & {\cellcolor[HTML]{FFC7CE}} \color[HTML]{9C0006}   & {\cellcolor[HTML]{FFC7CE}} \color[HTML]{9C0006}   & {\cellcolor[HTML]{FFC7CE}} \color[HTML]{9C0006} \\
 \midrule
 Group 3  &{\cellcolor[HTML]{8ED08B}} \color[HTML]{000000} 1 & {\cellcolor[HTML]{FFC7CE}} \color[HTML]{9C0006}   & {\cellcolor[HTML]{FFC7CE}} \color[HTML]{9C0006}   & {\cellcolor[HTML]{FFC7CE}} \color[HTML]{9C0006}   & {\cellcolor[HTML]{FFC7CE}} \color[HTML]{9C0006} \\

\bottomrule
\end{tabular}
\end{table}

After subdivisions are defined, it should be determined how many samples to take from each subdivision by assigning weights $w$ to each group. This can be done using any distribution in principle, however we consider three to be most relevant for architecture optimization problems: uniform weighting ($w = 1$) as applied in~\cite{Kaltenecker2019}, weighting by number of active design variables ($w = n_{\mathrm{act}}$) as applied in~\cite{Pelamatti2020a}, or weighting by group size ($w = \sqrt{n_{x,\mathrm{grp}}}$, where $n_{x,\mathrm{grp}}$ is the number of $x$ in the group). The reason for the latter two is that although smaller subdivisions (\textit{i.e.} subdivisions with less active variables and therefore less possible discrete design vectors) should not be under-represented compared to larger subdivisions, larger subdivisions might need more samples to give a good overview of subdivision behavior to the optimization algorithm.
Weighting by $n_{x,\mathrm{grp}}$ uses $w = \sqrt{n_{x,\mathrm{grp}}}$, because if $n_{x,\mathrm{grp}}$ is used directly it is equivalent to hierarchical sampling without grouping. The square root is applied to prevent oversampling large groups.
Table~\ref{tab:hier_sampling} shows how the valid discrete design vectors $x_{\mathrm{valid,discr}}$ are separated into groups based on $x_{\mathrm{act}}$ and weights are assigned based on $n_{\mathrm{act}}$. The relative weighting $w_{\mathrm{rel}}$ then determines how many of the requested samples are sampled from each group. For example, if 100 samples are requested, 36 will come from the first group, 36 from the second group, and 18 and 9 from the last two groups, respectively.
If there are less discrete vectors in a group than requested, vectors can be selected multiple times if there are also continuous variables in the problem.

After valid discrete design vectors have been sampled, active continuous variables are assigned values using Sobol' sampling~\cite{Renardy2021}.
To summarize, the hierarchical sampling algorithm consists of the following steps (visualized in Table~\ref{tab:hier_sampling}):
\begin{enumerate}
    \item Generate all possible valid discrete design vectors $x_{\mathrm{valid,discr}}$ and obtain associated activeness information $\delta$;
    \item Group design vectors (\textit{e.g.} by active variables $x_{\mathrm{act}}$) and weight each group (\textit{e.g.} by the number of active variables $n_{\mathrm{act}}$);
    \item Sample discrete design vectors according to the weights of each group;
    \item Assign values to active continuous variables using Sobol' sampling.
\end{enumerate}
Table~\ref{tab:sampling_strats} provides an overview of discussed sampling strategies for hierarchical design spaces.
For problems with large design spaces, it might not be possible to generate $x_{\mathrm{valid,discr}}$ due to memory or time limits.
An overview of past applications of the architecture optimization method presented in~\cite{Bussemaker2024adsg} suggests that this limitation might occur above design space sizes of 100's of millions of valid design vectors, meaning that for many architecture optimization problems this would not pose a problem.
If, however, indeed $x_{\mathrm{valid,discr}}$ is not available for a given optimization problem, it excludes the use of the previously discussed hierarchical sampling algorithms and therefore requires the use of a conventional (non-hierarchical) sampling algorithm combined with a correction step.
\begin{table}
\caption{Overview of sampling strategies for hierarchical design spaces. $x_{\mathrm{valid,discr}}$ refers to all valid discrete design vectors.}\label{tab:sampling_strats}
\begin{tabular}{lllll}
\toprule
Strategy& Needs & Grouping& Weighting \\
\midrule
Non-hierarchical & Correction &  & \\
\midrule
 \multirow{1}{*}{Hierarchical}& \multirow{1}{*}{$x_{\mathrm{valid,discr}}$} & None &  & \\
 && $n_{\mathrm{act}}$&  \\
 && $n_{\mathrm{act}}$& $n_{\mathrm{act}}$ \\
 && $n_{\mathrm{act}}$& $n_{x,\mathrm{grp}}$ \\
 && $x_{\mathrm{act}}$& \\
 && $x_{\mathrm{act}}$& $n_{\mathrm{act}}$ \\
 && $x_{\mathrm{act}}$& $n_{x,\mathrm{grp}}$ \\
 && $\mathrm{MRD}$& \\
 && $\mathrm{MRD}$& $n_{\mathrm{act}}$ \\
 && $\mathrm{MRD}$& $n_{x,\mathrm{grp}}$ \\
 \bottomrule
\end{tabular}
\end{table}

First, the sampling strategies listed in Table~\ref{tab:sampling_strats} are tested on NSGA-II.
We use the implementation of NSGA-II available in \noun{SBArchOpt} as highlighted in Section~\ref{subsec:SBArchOpt}.
Problem-specific correction is used and $x_{\mathrm{valid,discr}}$ is assumed to be available.
NSGA-II is executed with a DoE size of $10 \cdot n_x$, 25 generations and 100 repetitions.
Optimization performance is compared based on $\Delta \mathrm{HV}$ regret using the procedure described in Appendix~\ref{sec:perf_rank}.
$\Delta \mathrm{HV}$ ($\Delta$ hypervolume) represents the distance to the known optimum (or Pareto front in case of multi-objective optimization) normalized to the range of objective values.
Table~\ref{tab:sampling_nsga2_res} presents the results of sampling strategies for NSGA-II.
The hierarchical non-grouping sampling performs worst, especially on the rocket problems and weight-minimizing GNC problem (rank 8 and 455\% penalty compared to the best strategy). This is because in those problems (part of) the optimum lies in architectures represented by only 0.3\% of $x_{\mathrm{valid,discr}}$ (1 sensor/computer for the GND problem; 1 stage for the Rocket problem). The non-grouping sampler uniformly samples over all $x_{\mathrm{valid,discr}}$ and therefore has a high chance of not sampling these architectures. Problems with high $\mathrm{MRD}$ potentially suffer from this effect, depending on where the optimum lies.
From the other samplers, the hierarchical samplers without weighting consistently perform better than samplers with weighting.
Among the non-weighted hierarchical samplers, the $n_{\mathrm{act}}$ sampler performs best, with the other hierarchical samplers incurring between 6\% and 7\% relative penalty.

For the BO algorithm, the non-hierarchical sampler and hierarchical non-weighted samplers are compared.
The BO algorithm is executed with $n_{\mathrm{doe}} = 3 \cdot n_x$, 40 infill points and 24 repetitions. For the Jet SM problem, $n_{\mathrm{doe}} = 10 \cdot n_x$ for BO will be used, to correct for the fact that this problem features approximately a 60\% failure rate and therefore needs a larger DoE.
The weight-minimizing GNC problem ("Wt GNC") is not included, as it proved too easy for the BO algorithm to solve which lead to arbitrary differences in $\Delta \mathrm{HV}$ regret.
Table~\ref{tab:sampling_bo_res} presents sampling results for the BO algorithm.
It shows the hierarchical $x_{\mathrm{act}}$ sampler performs best.
The hierarchical $n_{\mathrm{act}}$ sampler performs worst with an 18\% mean performance penalty.

Non-hierarchical sampling performs slightly worse than the best sampler, at a 12\% mean performance penalty both for NSGA-II and for the BO algorithm.
Therefore, although the availability of $x_{\mathrm{valid,discr}}$ improves algorithm performance, the non-availability does not prevent the problem from being solved.
This observation is in line with results published in~\cite{Bussemaker2024adsg}, where it is shown that using a design space encoder ("ADORE Fast") that does not provide $x_{\mathrm{valid,discr}}$ (regardless of design space size) does indeed reduce optimization performance for both NSGA-II and BO, however not so severely as to prevent the optimization problem from being solved.

When comparing between NSGA-II and the BO algorithm, therefore, $x_{\mathrm{act}}$ and $\mathrm{MRD}$ remain as good candidates for hierarchical sampling.
We select hierarchical sampling based on $x_{\mathrm{act}}$ for its slightly better performance on the BO algorithm.
If $x_{\mathrm{valid,discr}}$ is not available, the non-hierarchical sampler is used instead.

\begin{landscape}

\begin{table}
\centering
\small
\caption{Best performing sampling strategies on various test problems running NSGA-II, ranked by $\Delta \mathrm{HV}$ regret (lower rank is better). Penalty represents the mean $\Delta \mathrm{HV}$ regret increase compared to the best infill. Best performing strategy is underlined; darker colors represent better results.
Abbreviation: grp. = grouping, wt. = weighting, hier. = hierarchical.
}\label{tab:sampling_nsga2_res}

\input{sampling_rank_nsga2_mod}
\end{table}

\begin{table}
\centering
\small
\caption{Best performing sampling strategies on various test problems running the BO algorithm, ranked by $\Delta \mathrm{HV}$ regret (lower rank is better). Penalty represents the mean $\Delta \mathrm{HV}$ regret increase compared to the best infill. Best performing strategy is underlined; darker colors represent better results.
Abbreviation: grp. = grouping, wt. = weighting, hier. = hierarchical.
}\label{tab:sampling_bo_res}
\input{sampling_rank_sbo_mod}
\end{table}

\end{landscape}

\subsection{Correction}\label{sec:HierCorrExp}

Correction is the operation of creating a valid design vector from an invalid design vector. Here we only discuss the process of correction for discrete design vectors, as we assume that only discrete variables determine the hierarchical structure in terms of activeness and value constraints in this work. Continuous variables are therefore only subject to imputation if they are inactive.

The mechanism of correction depends on the hierarchical structure of the optimization problem and therefore can be implemented on a per-problem basis, as in~\cite{Bussemaker2021c}. There, design variables are corrected during parsing of the design vector into an architecture description: design variables not representing a valid option value are corrected one-by-one to the closest valid value.
For example, when a design variable selects the third compressor stage for bleed off-take when only two compressor stages are available, the design variable is modified to select the second compressor stage for bleed off-take instead.
This is called a greedy correction algorithm, as it takes the locally best choice for correcting invalid values for each invalid design variable. It is arguably the most convenient way to implement {problem-specific} correction, because of this step-by-step correction mechanism.

\begin{figure}
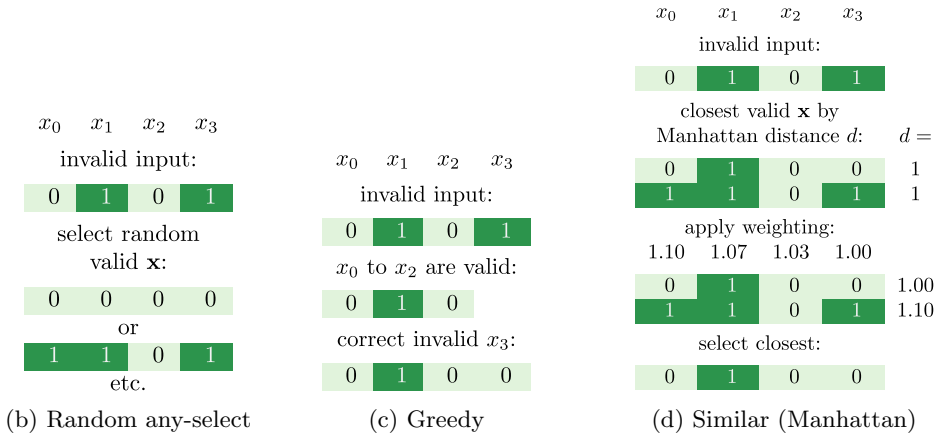

\centering
    \begin{subfigure}[b]{0.45\textwidth}
        \centering
        \includesvg[inkscapelatex=false,width=.5\textwidth]{fig3_corr_example.svg}
        \caption{Example $x_{\mathrm{valid,discr}}$}
    \end{subfigure}
    \par\vspace{.5cm}
    \begin{subfigure}[b]{0.28\textwidth}
        \centering
        \includesvg[inkscapelatex=false,width=.75\textwidth]{fig4_corr_random.svg}
        \caption{Random any-select}\label{fig:eager_corr_random}
    \end{subfigure}
    \hfill
    \begin{subfigure}[b]{0.28\textwidth}
        \centering
        \includesvg[inkscapelatex=false,width=.75\textwidth]{fig5_corr_greedy.svg}
        \caption{Greedy}\label{fig:eager_corr_greedy}
    \end{subfigure}
    \hfill
    \begin{subfigure}[b]{0.4\textwidth}
        \centering
        \includesvg[inkscapelatex=false,width=.75\textwidth]{fig6_corr_similar.svg}
        \caption{Similar (Manhattan)}\label{fig:eager_corr_similar}
    \end{subfigure}
\caption{A visualization of different eager correction strategies}\label{fig:eager_corr}
\end{figure}

{Problem-agnostic} correction algorithms instead could reduce implementation effort and potentially improve optimizer performance. Such correction algorithms can be defined both for when $x_{\mathrm{valid,discr}}$ is available and for when it is not: correction algorithms with $x_{\mathrm{valid,discr}}$ available are called \textit{eager} algorithms (see Figure~\ref{fig:eager_corr}), as they require the upfront exposure of valid design vectors; when $x_{\mathrm{valid,discr}}$ is not available the correction algorithm is known as a \textit{lazy} algorithm, as design vector validity is checked during the correction process rather than upfront.
As eager correction algorithms have access to $x_{\mathrm{valid,discr}}$, the simplest algorithm selects any of the available vectors as a replacement (called "Any-select"). We can define one that always returns the first (or any index for that matter) valid design vector, however that would not help with exploration at all: for problems with high imputation ratios there would still be a low chance to generate a new valid design vector. A better option therefore is to select a random valid design vector, as shown in Figure~\ref{fig:eager_corr_random}.

Another approach is to select a design vector that is close to the design vector to be corrected $x_{\mathrm{corr}}$. For eager algorithms this can be done in two ways: by applying a greedy algorithm or by selecting the most similar valid design vector.
The eager greedy algorithm (see Figure~\ref{fig:eager_corr_greedy}) repeatedly filters $x_{\mathrm{valid,discr}}$ based on the selected value of each design variable, starting from the left of $x_{\mathrm{corr}}$. If a given design variable value leaves no valid design vectors to be selected from, the closest value that does is selected. This is repeated until either only one valid design vector remains, or until all design variables have been processed. If in the latter case multiple design vectors remain, either the first (non-randomized configuration) or a random design vector is chosen (randomized configuration).
Eager selection based on similarity (see Figure~\ref{fig:eager_corr_similar}) is done by calculating the weighted distances from $x_{\mathrm{corr}}$ to all vectors in $x_{\mathrm{valid,discr}}$ and selecting the valid vector with the minimum distance to replace $x_{\mathrm{corr}}$.
Weighting factors linearly varying from 1.1 to 1.0 are applied in order to favor changes on the right side of the design vector over changes on the left, assuming that left-side design variables represent higher-impact decisions~\cite{Bussemaker2022c}.
As distance metrics either Euclidean or Manhattan distance can be used. If multiple valid vectors have the same minimum distance to $x_{\mathrm{corr}}$ either the first (non-randomized configuration) or a random design vector (randomized configuration) can be selected.

Lazy correction algorithms do not have access to $x_{\mathrm{valid,discr}}$ and instead provide some way to generate design vectors and check with a user-provided validation function $\mathrm{isCorrect}(\boldsymbol{x})$ whether that design vector is valid or not. This generation and checking process continues until a valid design vector is found.
Analogous to the any-selection eager algorithm, also an any-selection lazy algorithm can be defined. The difference in randomized and non-randomized selection however is that for non-randomized selection, the result may be remembered and returned in subsequent correction requests, whereas for randomized selection this is not the case.
For randomized selection, an average of $\mathrm{CR}$ design vectors (see also the discussion in Section~\ref{subsec:DSchallenges}) will have to be generated before finding a valid design vector, which represents the main limitation of lazy over eager algorithms.
Greedy selection is not possible for lazy algorithms, because $\mathrm{isCorrect}(\boldsymbol{x})$ only operates on complete design vectors and not fractions of design vectors as would be needed for one-by-one correction.
Lazy correction by similarity is done by modifying $x_{\mathrm{corr}}$ with some $\Delta_{\mathrm{corr}}$ vector, which can either be generated depth-first or distance-first.
Depth-first $\Delta_{\mathrm{corr}}$ generation directly applies the generated Cartesian product of $\Delta$ values for each design variable.
Distance-first generation first generates all possible $\Delta_{\mathrm{corr}}$ vectors and then sorts them by Euclidean or Manhattan distance before applying them to $x_{\mathrm{corr}}$. The same distance-weighting process as for eager similarity selection is used here.
\noun{ConfigSpace} uses a randomized version of lazy similarity correction: candidate design vectors are generated by replacing one design variable by a random valid value at a time. We do not consider this method, as it does not allow to correct multiple design variables at a time.
Table~\ref{tab:corr_strats} presents an overview of discussed correction strategies.

\begin{table}
\caption{Overview of correction algorithms for hierarchical design spaces. Euc and Manh refer to Euclidean and Manhattan distance metrics, respectively. See Figure~\ref{fig:eager_corr} for a visualization of eager correction algorithms.}\label{tab:corr_strats}
\begin{tabular}{llll}
\toprule
Type& Algorithm& Configuration&  Needs\\
\midrule
Problem-specific& &  &   Custom correction function\\
Eager& Any-select&  &   $x_{\mathrm{valid,discr}}$\\
 & Greedy&  &  $x_{\mathrm{valid,discr}}$\\
 & Similar&  Distance (Euc/Manh)&  $x_{\mathrm{valid,discr}}$\\
 Lazy& Any-select& & $\mathrm{isCorrect}(\boldsymbol{x})$\\
 & \multirow{2}{*}{Similar}& Depth- or distance-first& \multirow{2}{*}{$\mathrm{isCorrect}(\boldsymbol{x})$}\\
 & & Distance (Euc/Manh)& \\
\end{tabular}
\end{table}

The identified correction strategies are tested on NSGA-II first, with a DoE size of $10 \cdot n_x$, 25 generations and 100 repetitions.
Eager correction strategies are tested with the hierarchical $x_{\mathrm{act}}$ sampling algorithm; lazy strategies with the non-hierarchical sampler as here the assumption is that $x_{\mathrm{valid,discr}}$ is not available.
Correction strategies are tested for the various algorithm-specific configuration options (see Table~\ref{tab:corr_strats}).
Eager and lazy correction are additionally compared to problem-specific correction.
Table~\ref{tab:corr_nsga2_res} presents results of correction strategies for NSGA-II.
Eager correction performs better than lazy correction, and similarly to problem-specific correction.
Eager Any-select performs best, with Eager similar (Manhattan or Euclidean distance) performing similarly.
Lazy similar (depth-first) correction performs best among lazy correction strategies.
Lazy correction, however, takes between 1 and 2 orders of magnitude longer than problem-specific and eager correction: 2 to 50 ms, compared to 0.1 to 1 ms for eager and problem-specific correction. Additionally, lazy correction time increases linearly with $\mathrm{CR}$, because it is based on a trial-and-error approach.

\newgeometry{bottom=3cm} 
\begin{landscape} 

\begin{table}
\centering
\small
\caption{Best performing correction strategies on various test problems running NGSA-II, ranked by $\Delta \mathrm{HV}$ regret (lower rank is better). Penalty represents the mean $\Delta \mathrm{HV}$ regret increase compared to the best infill. Best performing strategy is underlined; darker colors represent better results.}\label{tab:corr_nsga2_res}
\input{corr_rank_nsga2_mod}
\end{table}

\begin{table}
\centering
\small
\caption{Best performing correction strategies on various test problems running the BO algorithm, ranked by $\Delta \mathrm{HV}$ regret (lower rank is better). Penalty represents the mean $\Delta \mathrm{HV}$ regret increase compared to the best infill. Best performing strategy is underlined; darker colors represent better results.}\label{tab:corr_bo_res}
\input{corr_rank_sbo_mod}
\end{table}

\end{landscape}
\restoregeometry

For the BO algorithm, best performing eager and lazy correction algorithms and problem-specific correction are compared.
The BO algorithm is executed with $n_{\mathrm{doe}} = 3 \cdot n_x$ ($n_{\mathrm{doe}} = 10 \cdot n_x$ for the Jet SM problem), 40 infill points and 24 repetitions.
Table~\ref{tab:corr_bo_res} presents sampling results for the BO algorithm.
Problem-specific correction with hierarchical sampling performs best.
Eager correction performs similar to problem-specific correction, however attains rank 1 less often.
Lazy correction and problem-specific correction with non-hierarchical sampling perform worst.
Based on these investigations, we conclude that problem-specific correction is sufficient for good optimizer performance, both for NSGA-II and BO.
If the user prefers not to implement problem-specific correction instead, and if $x_{\mathrm{valid,discr}}$ is available, then Eager Any-select correction should be used.
These results are used in the next section where the influence of hierarchy information integration strategies are investigated.

\section{Hierarchical Design Space Integration Strategies}\label{sec:HierIntegration}

According to \noun{Zaefferer \textit{et al.}}~\cite{Zaefferer2018a} there are three high-level strategies for considering design space hierarchy when implementing and solving an optimization problem:
\begin{enumerate}
    \item \textit{Naive}: no modification of the optimization algorithm at all, thereby effectively ignoring the effect of design variable hierarchy;
    \item \textit{Correction and imputation}: here correction and imputation are applied to ensure that all evaluated design vectors are valid (see also Section~\ref{subsec:DSchallenges});
    \item \textit{Explicit consideration}: the hierarchical structure is explicitly made available to and used by the optimization algorithm.
\end{enumerate}
The naive approach means that there might be a discrepancy between which design vectors the optimizer thinks are being evaluated and which design vectors actually are evaluated. For example, design vectors where only inactive design variables differ in value all represent the same system architecture instance. Not making this information available to the optimizer might lead to wasted computational resources, because exploration could be performed in sections of the design space that have no influence on performance. The optimization might also stall if the imputation ratio (see Section~\ref{subsec:DSchallenges}) is too high, because of the low chance of randomly finding valid design vectors.
As discussed in Section~\ref{subsec:DSchallenges}, \textit{correction and imputation} ensure that design vectors are valid by ensuring value constraints are satisfied and setting inactive design variables to some default value, respectively. Applying these operations avoids the aforementioned design vector mapping problems and only requires that the optimization algorithm accepts modified design vectors as output of the evaluation call.
Another way to support correction and imputation is through an \textit{ask-and-tell} interface~\cite{Collette2013}: here a process external to the optimizer has control over the optimization loop, "asking" the optimizer for one or more design vectors to evaluate and "telling" the optimizer the results after evaluation is finished. Results are "told" to the optimizer together with the corresponding design vectors, which means the ask-and-tell pattern allows integrating correction and imputation steps without any further modifications. Therefore, any optimization framework or algorithm that implements an ask-and-tell interface is compatible with correction and imputation, for example \noun{pymoo}~\cite{Blank2020}, \noun{BoTorch}~\cite{Balandat2019}, \noun{Trieste}~\cite{Picheny2023} and \noun{HEBO}~\cite{CowenRivers2020}.

Correction and imputation can also be implemented using a \textit{repair} operator~\cite{SalcedoSanz2009}: a problem-specific function that modifies design vectors, for example to satisfy an \textit{a-priori} constraint or otherwise improve the design points using heuristics.
The advantage of a repair operator over ask-and-tell or modifying design vectors in the evaluation function is that the correction and imputation operators are now available as a standalone function rather than always tied to evaluation. This allows correction and imputation to be applied during other steps in the optimization process than only evaluation, for example when generating initial design points or when searching the design space for the best infill point for surrogate-based optimization algorithms~\cite{GarridoMerchan2020}.
In this work, repair is only applied to discrete variables, however in principle it can also be applied to continuous variables, as was done for example in the work by \noun{Zaefferer \& Horn~\cite{Zaefferer2018a}.}

\begin{figure}
\centering
\includesvg[inkscapelatex=false,width=.9\textwidth]{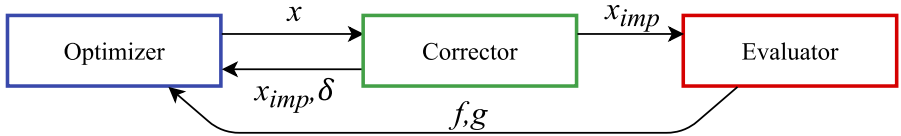}
\caption{High-level strategies for dealing with hierarchical optimization}\label{fig:hier_strat}
\end{figure}

The most invasive way of supporting hierarchical design spaces is through \textit{explicit consideration} by optimization algorithms. Here correction and imputation become an integral part of the optimizer and \textit{activeness} information as defined by the $\delta$-function (see Section~\ref{subsec:DSchallenges}) is made available to the optimizer. The availability of activeness information makes the following possible:
\begin{itemize}
    \item Generating all possible valid design vectors $x_{\mathrm{valid,discr}}$, and therefore using hierarchical sampling and problem-agnostic correction algorithms, as investigated in Section~\ref{sec:HierSampling}.
    \item Using hierarchical kernels in GP models used by BO~\cite{Hutter2013,Lu2018,Levesque2017,Saves2023SMT}.
\end{itemize}
Figure~\ref{fig:hier_strat} compares the three high-level integration strategies. It highlights the existence of a standalone corrector function that corrects and imputes design vectors and optionally returns activeness information.

One way to explicitly consider the hierarchical structure is by formally modeling the hierarchical structure and making this model available to the optimization algorithm.
It for example enables problem decomposition approaches~\cite{Frank2016,Pelamatti2020a} and the development of dedicated evolutionary search operators~\cite{Abdelkhalik2012,Nyew2015,Apaza2021}.
Additionally, because the model provides all information needed without needing to interrogate the problem definition (\textit{i.e.} as needed for a repair operator) this opens up the possibility for physically separating the optimizer from the function evaluation, for example to enable remote ask-and-tell execution.
Hierarchical design space models can be classified according to different levels of complexity:
\begin{enumerate}
    \item Single-level: one set of variables determining activeness of a disjoint set of conditional variables, \textit{e.g.}~\cite{Pelamatti2020a,Audet2022};
    \item Tree-structured: conditional variables can also determine activeness of other variables, \textit{e.g.} \noun{BoTorch}~\cite{Balandat2019};
    \item Directed acyclic graph: activeness can be determined by multiple variables, \textit{e.g.}~\cite{Hutter2013,HalleHannan2024} and \noun{ConfigSpace}~\cite{Lindauer2019};
    \item Directed graph: also supports cyclic dependencies, \textit{e.g.} the Design Space Graph (DSG)~\cite{Bussemaker2024adsg}.
\end{enumerate}
\noun{ConfigSpace}\footnote{\href{https://automl.github.io/ConfigSpace/}{https://automl.github.io/ConfigSpace/}} implements one of the most general hierarchical design space definitions known to the authors. \noun{ConfigSpace} is an open-source Python library that supports the definition of mixed-discrete design variables, conditional activation of design variables, and value constraints~\cite{Lindauer2019}.
Conditional activation is specified by if-clauses that can be composed of equality and inequality checks and Boolean conjunctions (AND, OR). Conditional activation depending on multiple variables is supported, as long as there are no cyclic dependencies.
Value constraints can be specified using forbidden clauses: generic if-clauses between two or more design variables that define a violated value constraint when the condition is true.
Once defined, the design space object can then be used to query which design variables are conditionally active, query whether a given design vector is valid, correct and impute a design vector, get activeness information for a given design vector, and generate random valid design vectors.
\noun{ConfigSpace} is used by several optimization frameworks to support hierarchical design spaces, for example SMAC3~\cite{Lindauer2022}, BOAH~\cite{Lindauer2019}, OpenBox~\cite{Li2021} and SMT~\cite{Saves2023SMT}.
It should be noted that also if the optimization problem does not expose an explicit hierarchical design space model, the problem still contains some model of the hierarchical structure, either explicitly defined or implicitly embedded in the evaluation code (\textit{e.g.} as for the turbofan optimization benchmark problem of~\cite{Bussemaker2021c}). It might be more convenient for the user to think in terms of a domain-specific model of the design space compared to directly thinking in terms of design variables, and thus it might help adoption of architecture optimization.
Examples of design space modeling techniques developed in the context of systems engineering and architecture optimization include the Architecture Decision Graph~\cite{Simmons2008}, the Adaptive Reconfigurable Matrix of Alternatives~\cite{Mavris2008}, the Architecture Decision Diagram~\cite{Apaza2021}, feature models~\cite{Benavides2010,Czarnecki2012}, function-means models~\cite{Gedell2012}, and the Architecture Design Space Graph (ADSG)~\cite{Bussemaker2022c}.
Retiarii~\cite{Zhang2020} and NASBench~\cite{Ying2019} provide domain-specific formats for modeling deep neural network architecture design spaces~\cite{Elsken2019}.
Table~\ref{tab:hier_strat} presents a detailed overview of discussed integration strategies, and capabilities gained at each level of integration.

\begin{table}
\caption{Different strategies for dealing with design variable hierarchy in optimization algorithms.}\label{tab:hier_strat}
\begin{tabular}{lllllllll}
\toprule
General strategy & Integration  &Usage& \footnotemark[1] & \footnotemark[2] & \footnotemark[3]  &\footnotemark[4]& \footnotemark[5]& \footnotemark[6]\\
\midrule
Naive & N/A  &N/A&  &  &   &&  & \\
Correction \& Imputation & $x$-output &Evaluation& $\checkmark$ &  &   &&  & \\
 & Ask-and-tell  &Evaluation& $\checkmark$ &  &   &&  & \\
 & Repair operator  &Evaluation, sampling& $\checkmark$ & $\checkmark$ &   &&  & \\
Explicit consideration & Activeness $\delta_i(\textbf{x})$ &Sampling& $\checkmark$ & $\checkmark$ & $\checkmark$  &&  & \\
 & Activeness $\delta_i(\textbf{x})$ &Sampling \& modeling& $\checkmark$ & $\checkmark$ & $\checkmark$  & $\checkmark$  & &\\
 & Formal model  &Sampling \& modeling& $\checkmark$ & $\checkmark$ & $\checkmark$  &$\checkmark$  & $\checkmark$ & $\checkmark$\\
 \bottomrule
\end{tabular}
\footnotetext[1]{All evaluated $x$ are valid.}
\footnotetext[2]{All $x$ in sampling and infill search are valid.}
\footnotetext[3]{Availability of all valid discrete design vectors $x_{\mathrm{valid,discr}}$.}
\footnotetext[4]{Hierarchical kernels for surrogate modeling.}
\footnotetext[5]{Dedicated search operators and possibility for problem decomposition.}
\footnotetext[6]{Physical separation between optimization and evaluation code.}
\end{table}

The correction and sampling algorithms developed in Section~\ref{sec:HierSampling} need the highest level of hierarchy integration: activeness information.
%
We now compare the different levels of integration to investigate if applying higher levels of integration (see Table~\ref{tab:hier_strat}) indeed results in better optimizer performance.
Ask-and-tell integration has the same results as $x$-output integration: they only differ in the way the algorithm is called. Similarly, having a formal model changes nothing in the information available to the algorithm compared to the activeness integration. Therefore, experiments are only run for the naive, $x$-output, repair and activeness integration strategies.
For BO, the activeness strategy is run in two configurations: one where activeness is only used for hierarchical sampling, and one where activeness is available additionally for the GP models.
An overview of tested integration strategies and available capabilities for each strategy is provided in Table~\ref{tab:hier_strat_tested}.
NSGA-II is executed with $n_{\mathrm{doe}} = 10 \cdot n_x$, 25 generations and 100 repetitions.
The BO algorithm is executed with $n_{\mathrm{doe}} = 3 \cdot n_x$ ($n_{\mathrm{doe}} = 10 \cdot n_x$ for the Jet SM problem), 40 infill points and 16 repetitions.
The hierarchical test problems presented in Section~\ref{sec:HierTestProblems} and shown in Table~\ref{tab:hier_probs} are used to compare the hierarchical optimization strategies.

\begin{table}
\caption{Tested hierarchy integration strategies and impact on available capabilities.  
General strategy refers to strategies presented in Table~\ref{tab:hier_strat}.
Abbreviations: hier. = hierarchical, corr. = correction, expl. cons. = explicit consideration.}\label{tab:hier_strat_tested}
\begin{tabular}{llllll}
\toprule
& Naive & X out & Repair & Hier. sampling & Activeness \\
\midrule
General strategy & Naive & Corr. & Corr. & Expl. cons. & Expl. cons. \\
Valid $x$-output & & $\checkmark$ & $\checkmark$ & $\checkmark$ & $\checkmark$ \\
Repair operator & & & $\checkmark$ & $\checkmark$ & $\checkmark$ \\
Hierarchical sampling & & & & $\checkmark$ & $\checkmark$ \\
Hierarchical GP (BO only) & & & & & $\checkmark$ \\
\bottomrule
\end{tabular}
\end{table}

Table~\ref{tab:hier_strat_nsga2} shows that for NSGA-II, the level of integration does not influence performance much.
Repair performs best, indicating that hierarchical sampling (as used in Activeness) is not necessarily beneficial.
Table~\ref{tab:hier_strat_sbo} shows that for BO a higher level of integration improves optimizer performance. 
Naive, X-out and Repair integration is penalized significantly (134\%, 149\% and 89\%, respectively), showing that in these cases the BO algorithm is less well able to suggest valid infill design points.
It should be noted that if it is not possible to generate $x_{\mathrm{valid,discr}}$ due to memory or time limits, the Activeness integration can still be used, however without hierarchical sampling.
The unavailability of $x_{\mathrm{valid,discr}}$ and therefore hierarchical sampling does reduce optimizer performance, as shown and discussed in Section~\ref{sec:HierSamplingExp}, however not as severe as switching to Repair integration.
The unavailability of $x_{\mathrm{valid,discr}}$ therefore does not represent a major obstacle to solving SAO problems.

If we consider only querying a Gaussian Process (GP) model, as shown by \noun{Saves \textit{et al.}} in~\cite{Saves2023SMT}, taking into account the activeness in the GP using hierarchical kernels leads to better results than using non-hierarchical kernels.
However, Table~\ref{tab:hier_strat_sbo} shows a reduction in performance when using hierarchical kernels (Activeness) compared to non-hierarchical kernels (Hierarchical sampling).
This result is due to two factors: the potentially high concentration of training points in some areas of the design space, and the structure of the used test problems.
High concentrations of training points in certain areas of the design space occur because the algorithm is attempting to close-in on an optimum in these areas, and is a common occurrence in SBO.
The Rocket and GNC problems have relatively smooth objective functions, and therefore in their hierarchical versions the non-hierarchical GP models are still able to accurately represent the objective functions.
It is expected that the first effect will be reduced for higher-dimensional problems, and that the second effect will be reduced for problems with more non-linear and non-smooth objective and constraint functions. Therefore we keep both the non-hierarchical and hierarchical GP models into consideration for the subsequent more complex investigation.

\begin{landscape}

\begin{table}
\caption{Comparison of hierarchical optimization strategies on various test problems running NSGA-II, ranked by $\Delta \mathrm{HV}$ regret (lower rank is better). Penalty represents the mean $\Delta \mathrm{HV}$ regret increase compared to the best infill. Best performing strategy is underlined; darker colors represent better results.}\label{tab:hier_strat_nsga2}
\input{hier_rank_nsga2}
\end{table}

\begin{table}
\caption{Comparison of hierarchical optimization strategies on various test problems running the BO algorithm, ranked by $\Delta \mathrm{HV}$ regret (lower rank is better). Penalty represents the mean $\Delta \mathrm{HV}$ regret increase compared to the best infill. Best performing strategy is underlined; darker colors represent better results.}\label{tab:hier_strat_sbo}
\input{hier_rank_sbo}
\end{table}

\end{landscape}

\section{Application: Jet Engine Architecture}\label{sec:JetEngineOpt}

\begin{figure}
\centering
\includesvg[inkscapelatex=false,width=\textwidth]{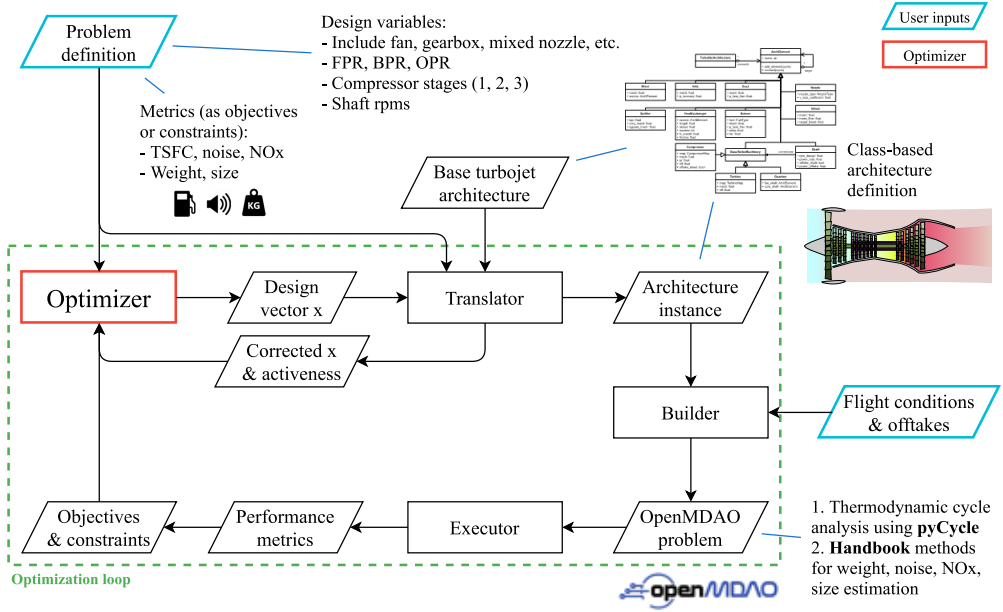}
\caption{Overview of the jet engine optimization testing framework. The user provides the problem definition (in terms of design variables and metrics selected from a database) and the flight conditions and power offtakes to size the engine for.
Engine schematic adapted from original by K. Aainsqatsi, available at \href{https://commons.wikimedia.org/wiki/File:Turbofan_operation.svg}{\footnotesize https://commons.wikimedia.org/wiki/File:Turbofan\_operation.svg}.}\label{fig:jet_engine_prob}
\end{figure}

To demonstrate the application of the architecture optimization strategies presented in this work, a jet engine optimization problem is solved using a Bayesian Optimization (BO) algorithm.
The jet engine optimization problem is a benchmark problem specifically developed to provide realistic architecture optimization challenges and behavior~\cite{Bussemaker2021c}.
It is defined using the jet engine optimization testing framework presented by \noun{Bussemaker \textit{et al.}}~\cite{Bussemaker2021c}, the purpose of which is the provide a flexible way to define jet engine optimization problems with the purpose of testing optimization algorithms for SAO, see Figure~\ref{fig:jet_engine_prob} for an overview.
The user defines the optimization problem by selecting from available architectural choices (that define the design variables) and metrics and inputting the flight conditions and power offtakes to size the engine for.
Available architectural choices include whether to add a fan and bypass flow, the number of compressor and turbine stages, the use of intercooling and inter-turbine burning, and where to apply bleed and power offtakes. The selected choices are used to define the continuous, integer, and categorical design variables for the optimizer.
A translator code is provided that translates a design vector generated by the optimizer into an architecture instance defined using objects. These objects contain all information required to build the analysis problem, including input parameters (from flight conditions, offtake requirements, or design vectors) and airflow connection sequences (\textit{e.g.} compressor to combustor, combustor to turbine, etc.).
A builder code then takes these objects and constructs an OpenMDAO~\cite{Gray2019} problem that performs thermodynamic cycle analysis and engine sizing using pyCycle~\cite{Hendricks2019}, with the main output being the Thrust-Specific Fuel Consumption (TSFC) of the engine. Handbook methods are added to calculate additional metrics such as noise level, NOx emissions, weight, and size.
Thermodynamic cycle analysis takes between 1 and 5 minutes to complete. However, it is not guaranteed to converge to a feasible solution, leading to the presence of a hidden constraint. If the hidden constraint is violated, metrics are set to NaN (not a number).
The testing framework enables specification of a wide variety of test problems, all based on realistic engineering behavior, however with varying number of design variables, objectives, and constraints.
The code is available open source\footnote{\href{https://github.com/jbussemaker/OpenTurbofanArchitecting/}{https://github.com/jbussemaker/OpenTurbofanArchitecting/}}, for more information the reader is referred to~\cite{Bussemaker2021c}.

In this investigation, we use the simple problem formulation defined in~\cite{Bussemaker2021c}:

\color{black}

\begin{equation}\begin{array}{cclll}   
\mathrm{minimize} & & \mathrm{TSFC} && \\
\mathrm{w.r.t.} && \mathrm{IncludeFan} \in \{ \mathrm{False}, \mathrm{True} \} \\
 
 && \mathrm{if\ IncludeFan} = \mathrm{True}: \\
 && \hspace{5mm}2.0 \leq \mathrm{BPR} \leq 12.5 \\
 && \hspace{5mm}1.1 \leq \mathrm{FPR} \leq 1.8 \\
 && \hspace{5mm}\mathrm{MixedNozzle} \in \{ \mathrm{False}, \mathrm{True} \} \\
 && \hspace{5mm}\mathrm{IncludeGearbox} \in \{ \mathrm{False}, \mathrm{True} \} \\
 && \hspace{5mm}\mathrm{if\ IncludeGearbox} = \mathrm{True}: \\
 && \hspace{5mm}\hspace{5mm}1.0 \leq \mathrm{GearRatio} \leq 5.0 \\
 
 && 1.1 \leq \mathrm{OPR} \leq 60.0 \\
 && n_{\mathrm{shafts}} \in \{ 1, 2, 3 \} \\
 && \mathrm{if\ } n_{\mathrm{shafts}} > 1: \\
 && \hspace{5mm} 0.1 \leq \mathrm\mathrm_{\mathrm{factor,i}} \leq 0.9 & & i=2,\dots,n_{\mathrm{shafts}} \\
 && 1000 \leq \mathrm{RPM}_i \leq 20000 & & i=1,\dots,n_{\mathrm{shafts}} \\
 && \mathrm{PowerOfftake} \in \{ 1,\dots,n_{\mathrm{shafts}} \} \\
 && \mathrm{BleedOfftake} \in \{ 1,\dots,n_{\mathrm{shafts}} \} \\
\mathrm{subject}\ \mathrm{to} && M_{\mathrm{jet}} \leq 1.0 \\
 && \mathrm{PR}_{\mathrm{factor,sum}} \leq 0.9 \\
 && \mathrm{PR}_{\mathrm{max,i}} \leq 15.0 & & i=1,2,3 \\
\end{array}\end{equation}
which is a single-objective (TSFC minimization) problem with several architectural choices: fan inclusion IncludeFan, number of compressor stages $n_{\mathrm{shafts}}$, gearbox inclusion IncludeGearbox, mixed nozzle selection MixedNozzle and power offtake locations PowerOfftake and BleedOfftake.
The problem includes several levels of activation hierarchy: bypass ratio BPR, fan pressure ratio FPR, gearbox inclusion IncludeGearbox and mixed nozzle selection MixedNozzle are only active if the fan is included ($\mathrm{IncludeFan} = \mathrm{True}$); the gear ratio GearRatio is only active if $\mathrm{IncludeGearbox} = \mathrm{True}$; and shaft-related pressure ratio factors $\mathrm{PR}_{\mathrm{factor,i}}$ and rotational speeds $\mathrm{RPM}_i$ are only active if the respective number of shafts are selected.
The power offtake PowerOfftake and bleed offtake BleedOfftake selections are always active, however are value-constrained by the available shafts.

In total, there are 70 valid discrete design vectors. However, the Cartesian product of discrete variables leads to 216 combinations: the discrete imputation ratio therefore is $\mathrm{IR}_d = 216 / 70 = 3.1$ (see Eq.~\eqref{eq:IRd}).
The continuous imputation ratio $\mathrm{IR}_c = 1.26$ (see Eq.~\eqref{eq:IRc}), meaning that there are on average $9 / 1.26 = 7.14$ continuous variables active (as seen over all valid discrete design vectors).
The overall imputation ratio is $\mathrm{IR} = 3.89$ (see Eq.~\eqref{eq:IR}).
The correction ratio $\mathrm{CR} = 2.10$ (see Eq.~\eqref{eq:CR}), which leads to correction ratio fraction $\mathrm{CRF} = 55\%$ (see Eq.~\eqref{eq:CRF}). Thus, a little over half of the design space hierarchy is due to value constraints (\textit{i.e.} the need for correction).
The problem additionally features 5 design constraints, constraining the output jet Mach number $M_{\mathrm{jet}}$ and pressure ratio distributions over the selected compressor stages ($\mathrm{PR}_{\mathrm{factor,sum}}$ and $\mathrm{PR}_{\mathrm{max,i}}$).
Additionally, the underlying thermodynamic cycle analysis and sizing code does not always converge, leading to a hidden constraint being violated in approximately 50\% of design points generated in a random DoE.
We use the problem implemented in \noun{SBArchOpt} as \noun{SimpleTurbofanArch}.

The BO algorithm is executed 24 times with $n_{\mathrm{doe}} = 113$, $n_{\mathrm{infill}} = 187$ (a budget of 300 evaluations), and $n_{\mathrm{batch}} = 4$.
The algorithm is executed for (see Table~\ref{tab:hier_strat_tested}) Repair integration (no hierarchy information exposed, however the repair operator is available), Hierarchical Sampling (therefore an MD GP used during the optimization), and Activeness (both hierarchical sampling and hierarchical GP is used).
The effectiveness of NSGA-II has already been demonstrated in~\cite{Bussemaker2021c}, however for completeness we compare BO results against NSGA-II results. NSGA-II is executed with repair operator available and using hierarchical sampling.

\begin{figure}
\centering
\includesvg[inkscapelatex=false,width=0.9\textwidth]{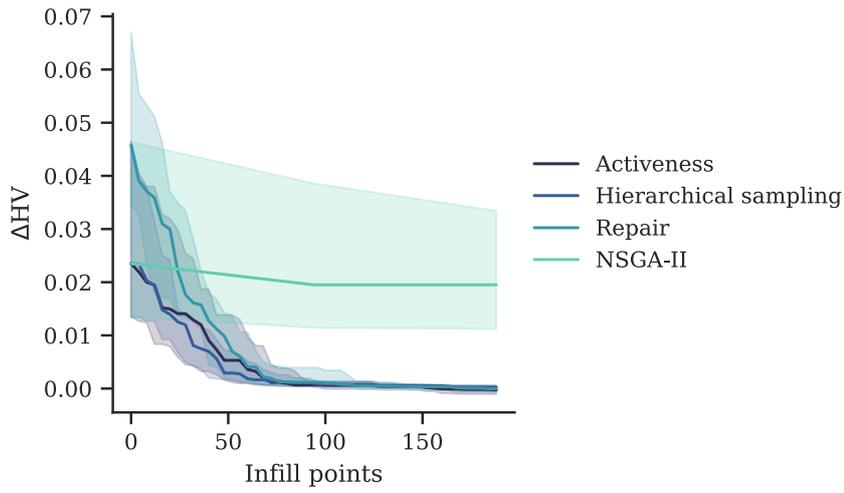}
\caption{Comparison of NSGA-II to the BO algorithm with different levels of hierarchy integration (see also Table~\ref{tab:hier_strat_tested}): Repair (no hierarchy information; repair operator available), Hierarchical sampling (no hierarchical GP) and Activeness (hierarchical sampling and hierarchical GP)}\label{fig:jet_eng_bo}
\end{figure}

\begin{figure}
\centering
\includesvg[inkscapelatex=false,width=0.9\textwidth]{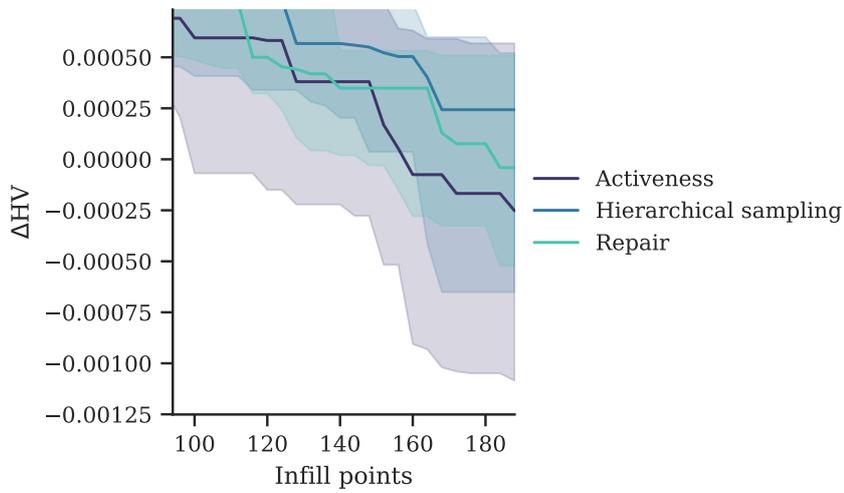}
\caption{Details of final optimization phase of comparison shown in Figure~\ref{fig:jet_eng_bo}, only showing the BO algorithm}\label{fig:jet_eng_bo_zoom}
\end{figure}

\begin{table}
\caption{Comparison of median optimal TSFC (minimization) values achieved for the jet engine problem.
Refer to Table~\ref{tab:hier_strat_tested} for comparison of hierarchy integration strategies.
$\Delta \mathrm{HV}$ regret was not available for NSGA-II with 3250 evaluations.
}\label{tab:jet_engine_opt}
\begin{tabular}{lllll}
\toprule
Algorithm & Hierarchy integration & $N_{\mathrm{fe}}$ & TSFC [g/kNs] & $\Delta \mathrm{HV}$ regret \\
\midrule
BO & Activeness            & 300  & 6.633           & 0.91 \\
BO & Hierarchical sampling & 300  & 6.653 (+0.30\%) & 0.86 (-4.8\%) \\
BO & Repair                & 300  & 6.640 (+0.11\%) & 1.59 (+75\%) \\
NSGA-II &                  & 300  & 7.455 (+12.4\%) & 3.67 (+305\%) \\
NSGA-II &                  & 3250 & 6.640 (+0.11\%) & $-$ \\
\bottomrule
\end{tabular}
\end{table}

Figures~\ref{fig:jet_eng_bo} and~\ref{fig:jet_eng_bo_zoom} present the results of the jet engine optimization for the BO algorithm with the three compared hierarchy integration strategies.
Table~\ref{tab:jet_engine_opt} presents achieved median optimal TSFC and $\Delta \mathrm{HV}$ regret values.
All BO algorithm configurations approximate the previously-found optimum within 0.2\%, or slightly improve it, within the 300 function evaluations. The previously-found optimum was found with NSGA-II and an evaluation budget of 3250; BO therefore can be considered to be able to find the same result in 92\% less function evaluations.
For the BO algorithm, Activeness and Hierarchical sampling perform with similar $\Delta \mathrm{HV}$ regret (see Table~\ref{tab:jet_engine_opt}). However, Activeness is able to find a slightly better TSFC.
Repair yields a significantly higher $\Delta \mathrm{HV}$ regret due to less effective initial sampling, however in the end finds a better TSFC than Hierarchical sampling.
NSGA-II is not able to improve much upon its initial population within 300 function evaluations.
These results demonstrate that BO can be used to solve SAO problems, enabling a significant reduction in number of function evaluations needed compared to evolutionary algorithms like NSGA-II.
Integration of design space hierarchy information, both by using the hierarchical sampling algorithm and using hierarchical GP, results in better optimizer performance.

\section{Conclusions and Outlook}\label{sec:Conclusions}

System Architecture Optimization (SAO) problems are challenging to solve due to their mixed-discrete and hierarchical design spaces combined with constrained, multi-objective, black-box, and expensive to evaluate solution spaces that potentially are subject to hidden constraints.
Four metrics are defined for characterizing SAO design space hierarchy: imputation ratio $\mathrm{IR}$, correction ratio $\mathrm{CR}$, correction fraction $\mathrm{CRF}$, and maximum rate diversity $\mathrm{MRD}$.
The mixed-discrete, black-box nature of SAO problems require the use of gradient-free, global optimization algorithms. Multi-objective Evolutionary Algorithms (MOEA), in particular NSGA-II, have been shown in the past to be effective, although at a high function evaluation cost.
The potential of Surrogate-Based Optimization (SBO) and in particular \textit{Bayesian Optimization (BO)} algorithms is identified, mainly due to their high potential in solving global optimization problems with a restricted function evaluation budget.

This work extends a non-hierarchical BO algorithm to support optimization in hierarchical SAO design spaces.
First, hierarchical Gaussian Process (GP) kernels developed in previous research and used in the SMT library are extended to also support \textit{categorical hierarchical variables}. This enables the usage of hierarchical GP models in any mixed-discrete hierarchical design space.

Then, hierarchical sampling and correction algorithms are developed to deal with rate diversity and correction ratio effects, respectively.
The best performing and thereafter applied \textit{hierarchical sampling} algorithm groups discrete design vectors by active design variables $x_{\mathrm{act}}$, then uniformly samples discrete design vectors from the defined groups, and applies Sobol' sampling to generate samples for continuous variables.
Hierarchical sampling performs better than non-hierarchical sampling, however cannot be used if valid discrete design vectors $x_{\mathrm{valid,discr}}$ are not available.
\textit{Problem-agnostic correction algorithms} are developed to deal with imputation ratio effects. Eager correction algorithms have access to valid discrete design vectors $x_{\mathrm{valid,discr}}$, whereas lazy correction algorithms do not.
Using problem-agnostic correction does not necessarily lead to better optimizer performance when compared for both NSGA-II and BO, and therefore it is not applied in the rest of the work.

Information about the hierarchical design space can be integrated into an optimization algorithm at various levels: the information can be ignored (naive approach), design vectors can be modified either by the evaluation function or by a standalone repair operator (correction and imputation approach), and activeness information can additionally be made available and used by the optimizer (explicit consideration approach).
It is shown that for NSGA-II the level of integration does not influence optimizer performance much.
For BO, {higher levels of integration} increase optimizer performance, except that using a hierarchical GP for infill search reduces performance.

The developed BO algorithm is demonstrated on a jet engine optimization problem that features design space hierarchy, design constraints, and hidden constraints.
It is shown that BO is able to find the optimum in 92\% less function evaluations than NSGA-II. 
The higher levels of hierarchy integration perform better.
The developed BO algorithm is implemented as as \noun{ArchSBO} in \noun{SBArchOpt}\footnote{\href{https://sbarchopt.readthedocs.io/}{https://sbarchopt.readthedocs.io/}}~\cite{Bussemaker2023a}.
Used test problems are also available in \noun{SBArchOpt}.
The implemented algorithm is provided with a default configuration that should be appropriate for a wide range of optimization problems.

Although BO has been shown to be effective at solving moderate-size architecture optimization problems, use of BO should also be enabled for larger design space sizes, \textit{i.e.} tens to hundreds of design variables~\cite{Donelli2023}, for example using dimension reduction techniques such as EGORSE~\cite{Priem2023}.
Multi-fidelity BO may also be considered to make BO more effective and/or further reduce required computational resources.
To aid adoption of SAO, it should be easy to formulate and run optimizations.
Formulating architecture optimization problems can be made easier by developing and applying modeling techniques that do not require the explicit definition of design variables and design variable hierarchy, \textit{e.g.}~\cite{Bussemaker2022c,Zhang2020}.
Running optimizations can be made easier by enabling configuration and execution of optimization algorithms from the architecture modeling environment directly, rather than requiring the user to switch environments and/or write additional code.
It should be investigated whether formulating kernels for GP models directly based on similarity between architecture instances is possible. Specifically, for architecture instances that are represented as (directed) graphs, \textit{e.g.~\cite{Bussemaker2024adsg}, graph kernels~\cite{Nikolentzos2021,Sirico2023} should be investigated.}

\section*{Funding}
The research presented in this paper has been performed in the framework of the COLOSSUS project (Collaborative System of Systems Exploration of Aviation Products, Services and Business Models) and has received funding from the European Union Horizon Europe Programme under grant agreement n${^\circ}$ 101097120.

\section*{Acknowledgments}
The authors would like to recognize Luca Boggero and Björn Nagel for supporting the secondment of Jasper Bussemaker at ONERA DTIS in Toulouse, during which the research presented in this paper was mostly performed.

\section*{Data Availability and Replication of Results}

The dataset containing experimental results on hierarchical sampling and correction algorithms, hierarchical integration strategies, and the jet engine optimization application case is available at \href{https://doi.org/10.5281/zenodo.10683733}{https://doi.org/10.5281/zenodo.10683733}.
All test problems used in this work are available in the open-source \noun{SBArchOpt} library.
Code used to run the experiments is available open-source too, at \href{https://github.com/jbussemaker/ArchitectureOptimizationExperiments/tree/hidden-constraints}{{\small github.com/jbussemaker/ArchitectureOptimizationExperiments} (\noun{hidden-constraints} branch)}.
The repository maps generated results (Tables and Figures) to executable Python functions to enable reproducing results.
Developed results have been implemented in \noun{SBArchOpt} version 1.4.

\begin{appendices}

\section{Optimizer Performance Comparison Method}\label{sec:perf_rank}

This appendix describes the method for comparing optimization algorithm performance by ranking various algorithm configurations based on $\Delta \mathrm{HV}$ regret.
$\Delta \mathrm{HV}$ ($\Delta$ hypervolume) represents the distance to the known optimum (or Pareto front in case of multi-objective optimization) normalized to the range of objective values. For single-objective optimization, the range is calculated from the difference between the known optimum and the maximum objective value encountered during the optimization; for multi-objective optimization the hypervolume of the known Pareto front is taken as normalization value.
To compare multiple optimization runs with different initial Design of Experiments (DoEs), $\Delta \mathrm{HV}_{\mathrm{ratio}} = \Delta \mathrm{HV} / \Delta \mathrm{HV}_0$ is used, where $\Delta \mathrm{HV}_0$ is the value at the first iteration (\textit{i.e.} after the DoE has been evaluated).
Regret represents the cumulative error, in our case $\Delta \mathrm{HV}_{\mathrm{ratio}}$, over the course of an optimization~\cite{Garnett2023}: lower values are better, and the value gives an indication of both how closely the optimum is approach by the end of the optimization and by how quickly this was achieved. $\Delta \mathrm{HV}$ regret at iteration $i$ is calculated from:
\begin{equation}
    \mathrm{Regret}(\Delta \mathrm{HV})_i = \mathrm{Regret}(\Delta \mathrm{HV})_{i-1} + \frac{1}{2} n_{\mathrm{step}} \left( \Delta \mathrm{HV}_{\mathrm{ratio},i-1} + \Delta \mathrm{HV}_{\mathrm{ratio},i} \right)
\end{equation}
where $n_{\mathrm{step}} = n_{\mathrm{batch}}$ when comparing performance per function evaluation, and $n_{\mathrm{step}} = 1$ when comparing by infill iteration.
Performance is sampled $n_{\mathrm{samples}}$ times for the same configuration to correct for randomness.

For a given test problem, the best performing algorithm configuration has rank 1; higher ranks indicate lower performance. Similarly-performing configurations have the same rank, as tested by an independent two-sample t-test as implemented by Scipy's \noun{ttest\_ind\_from\_stats}\footnote{\href{https://docs.scipy.org/doc/scipy/reference/generated/scipy.stats.ttest_ind_from_stats.html}{https://docs.scipy.org/doc/scipy/reference/generated/scipy.stats.ttest\_ind\_from\_stats.html}} function. The algorithm for determining rank is listed in Algorithm~\ref{algo:perf_rank}.
The best performing algorithm configuration is then selected by counting ranks over multiple test problems: the best performing configuration is the one with the highest proportion of rank 1 within the set of configurations with highest proportion of rank $\leq$ 2.

\begin{algorithm}
\caption{Determine performance rank $R_i$ for each algorithm configuration $C_i$ given performance measure $p_i$ with standard deviation $\sigma_i$.}\label{algo:perf_rank}
\begin{algorithmic}[1]
\Require $C$, $p$, $\sigma$, $n_{\mathrm{samples}}$, $\mathrm{perfMin}$
\Ensure $R$
\State $R \leftarrow zeros$ \Comment{Initialize ranks to 0 (unevaluated)}
\If{$\mathrm{perfMin}$} \Comment{Get initial best performing configuration}
    \State $i_{\mathrm{comp}} \leftarrow \arg \min p$
\Else
    \State $i_{\mathrm{comp}} \leftarrow \arg \max p$
\EndIf
\State $R_{i_{\mathrm{comp}}} \leftarrow 1$ \Comment{Set initial best performing configuration to rank 1}
\While{$any(R = 0)$} \Comment{Loop while there are unevaluated ranks}
    \State $i_{\mathrm{uneval}} \leftarrow \arg R = 0$ \Comment{Get unevaluated ranks}
    \If{$\mathrm{perfMin}$} \Comment{Get best performing unevaluated configuration}
        \State $j_{\mathrm{comp}} \leftarrow \arg \min p(i_{\mathrm{uneval}})$
    \Else
        \State $j_{\mathrm{comp}} \leftarrow \arg \max p(i_{\mathrm{uneval}})$
    \EndIf
    \State $p_{\mathrm{same}} \leftarrow \mathrm{tTestInd}(p_{i_{\mathrm{comp}}}, \sigma_{i_{\mathrm{comp}}}, n_{\mathrm{samples}}, p_{j_{\mathrm{comp}}}, \sigma_{j_{\mathrm{comp}}}, n_{\mathrm{samples}})$
    \If{$p_{\mathrm{same}} \leq 10\%$} \Comment{If probability of being the same is low}
        \State $R_{j_{\mathrm{comp}}} \leftarrow R_{i_{\mathrm{comp}}}+1$ \Comment{Assign higher rank to compared configuration}
        \State $i_{\mathrm{comp}} \leftarrow j_{\mathrm{comp}}$ \Comment{Update reference configuration}
    \Else
        \State $R_{j_{\mathrm{comp}}} \leftarrow R_{i_{\mathrm{comp}}}$ \Comment{Assign same rank}
    \EndIf
\EndWhile
\end{algorithmic}
\end{algorithm}

\end{appendices}

\bibliography{article}

\end{document}

%% file: sampling_rank_nsga2_mod.tex
\begin{tabular}{lccccccccccc}
\toprule
Strategy (grp.; wt.) & RCost & RWt & Rocket & SO GNC & FR GNC & Wt GNC & MD GNC & Jet SM & Rank 1 & Rank $\leq$ 2 & Penalty \\
\midrule
Non-hier. & {\cellcolor[HTML]{05712F}} \color[HTML]{F1F1F1} 2 & {\cellcolor[HTML]{00441B}} \color[HTML]{F1F1F1} 1 & {\cellcolor[HTML]{05712F}} \color[HTML]{F1F1F1} 2 & {\cellcolor[HTML]{8ED08B}} \color[HTML]{000000} 5 & {\cellcolor[HTML]{56B567}} \color[HTML]{F1F1F1} 4 & {\cellcolor[HTML]{2C944C}} \color[HTML]{F1F1F1} 3 & {\cellcolor[HTML]{BCE4B5}} \color[HTML]{000000} 6 & {\cellcolor[HTML]{56B567}} \color[HTML]{F1F1F1} 4 & {\cellcolor[HTML]{DEEBF7}} \color[HTML]{000000} 12\% & {\cellcolor[HTML]{9DCAE1}} \color[HTML]{000000} 38\% & {\cellcolor[HTML]{006C2C}} \color[HTML]{F1F1F1} 12\% \\
Hier. & {\cellcolor[HTML]{BCE4B5}} \color[HTML]{000000} 6 & {\cellcolor[HTML]{8ED08B}} \color[HTML]{000000} 5 & {\cellcolor[HTML]{BCE4B5}} \color[HTML]{000000} 6 & {\cellcolor[HTML]{05712F}} \color[HTML]{F1F1F1} 2 & {\cellcolor[HTML]{00441B}} \color[HTML]{F1F1F1} 1 & {\cellcolor[HTML]{F7FCF5}} \color[HTML]{000000} 8 & {\cellcolor[HTML]{2C944C}} \color[HTML]{F1F1F1} 3 & {\cellcolor[HTML]{00441B}} \color[HTML]{F1F1F1} 1 & {\cellcolor[HTML]{C6DBEF}} \color[HTML]{000000} 25\% & {\cellcolor[HTML]{9DCAE1}} \color[HTML]{000000} 38\% & {\cellcolor[HTML]{F7FCF5}} \color[HTML]{000000} 455\% \\
\underline{Hier. ($n_{\mathrm{act}}$)} & {\cellcolor[HTML]{00441B}} \color[HTML]{F1F1F1} 1 & {\cellcolor[HTML]{05712F}} \color[HTML]{F1F1F1} 2 & {\cellcolor[HTML]{00441B}} \color[HTML]{F1F1F1} 1 & {\cellcolor[HTML]{8ED08B}} \color[HTML]{000000} 5 & {\cellcolor[HTML]{8ED08B}} \color[HTML]{000000} 5 & {\cellcolor[HTML]{00441B}} \color[HTML]{F1F1F1} 1 & {\cellcolor[HTML]{BCE4B5}} \color[HTML]{000000} 6 & {\cellcolor[HTML]{05712F}} \color[HTML]{F1F1F1} 2 & {\cellcolor[HTML]{9DCAE1}} \color[HTML]{000000} \underline{38\%} & {\cellcolor[HTML]{4191C6}} \color[HTML]{F1F1F1} \underline{62\%} & {\cellcolor[HTML]{00441B}} \color[HTML]{F1F1F1} \underline{0\%} \\
Hier. ($n_{\mathrm{act}}$; $n_{\mathrm{act}}$) & {\cellcolor[HTML]{05712F}} \color[HTML]{F1F1F1} 2 & {\cellcolor[HTML]{05712F}} \color[HTML]{F1F1F1} 2 & {\cellcolor[HTML]{2C944C}} \color[HTML]{F1F1F1} 3 & {\cellcolor[HTML]{2C944C}} \color[HTML]{F1F1F1} 3 & {\cellcolor[HTML]{2C944C}} \color[HTML]{F1F1F1} 3 & {\cellcolor[HTML]{2C944C}} \color[HTML]{F1F1F1} 3 & {\cellcolor[HTML]{56B567}} \color[HTML]{F1F1F1} 4 & {\cellcolor[HTML]{05712F}} \color[HTML]{F1F1F1} 2 & {\cellcolor[HTML]{F7FBFF}} \color[HTML]{000000} 0\% & {\cellcolor[HTML]{9DCAE1}} \color[HTML]{000000} 38\% & {\cellcolor[HTML]{006328}} \color[HTML]{F1F1F1} 9\% \\
Hier. ($n_{\mathrm{act}}$; $n_{x,\mathrm{grp}}$) & {\cellcolor[HTML]{56B567}} \color[HTML]{F1F1F1} 4 & {\cellcolor[HTML]{56B567}} \color[HTML]{F1F1F1} 4 & {\cellcolor[HTML]{8ED08B}} \color[HTML]{000000} 5 & {\cellcolor[HTML]{00441B}} \color[HTML]{F1F1F1} 1 & {\cellcolor[HTML]{05712F}} \color[HTML]{F1F1F1} 2 & {\cellcolor[HTML]{BCE4B5}} \color[HTML]{000000} 6 & {\cellcolor[HTML]{00441B}} \color[HTML]{F1F1F1} 1 & {\cellcolor[HTML]{00441B}} \color[HTML]{F1F1F1} 1 & {\cellcolor[HTML]{9DCAE1}} \color[HTML]{000000} 38\% & {\cellcolor[HTML]{6AAED6}} \color[HTML]{F1F1F1} 50\% & {\cellcolor[HTML]{D7EFD1}} \color[HTML]{000000} 81\% \\
Hier. ($x_{\mathrm{act}}$) & {\cellcolor[HTML]{00441B}} \color[HTML]{F1F1F1} 1 & {\cellcolor[HTML]{05712F}} \color[HTML]{F1F1F1} 2 & {\cellcolor[HTML]{05712F}} \color[HTML]{F1F1F1} 2 & {\cellcolor[HTML]{BCE4B5}} \color[HTML]{000000} 6 & {\cellcolor[HTML]{BCE4B5}} \color[HTML]{000000} 6 & {\cellcolor[HTML]{05712F}} \color[HTML]{F1F1F1} 2 & {\cellcolor[HTML]{BCE4B5}} \color[HTML]{000000} 6 & {\cellcolor[HTML]{2C944C}} \color[HTML]{F1F1F1} 3 & {\cellcolor[HTML]{DEEBF7}} \color[HTML]{000000} 12\% & {\cellcolor[HTML]{6AAED6}} \color[HTML]{F1F1F1} 50\% & {\cellcolor[HTML]{005B25}} \color[HTML]{F1F1F1} 7\% \\
Hier. ($x_{\mathrm{act}}$; $n_{\mathrm{act}}$) & {\cellcolor[HTML]{05712F}} \color[HTML]{F1F1F1} 2 & {\cellcolor[HTML]{00441B}} \color[HTML]{F1F1F1} 1 & {\cellcolor[HTML]{2C944C}} \color[HTML]{F1F1F1} 3 & {\cellcolor[HTML]{56B567}} \color[HTML]{F1F1F1} 4 & {\cellcolor[HTML]{56B567}} \color[HTML]{F1F1F1} 4 & {\cellcolor[HTML]{56B567}} \color[HTML]{F1F1F1} 4 & {\cellcolor[HTML]{56B567}} \color[HTML]{F1F1F1} 4 & {\cellcolor[HTML]{00441B}} \color[HTML]{F1F1F1} 1 & {\cellcolor[HTML]{C6DBEF}} \color[HTML]{000000} 25\% & {\cellcolor[HTML]{9DCAE1}} \color[HTML]{000000} 38\% & {\cellcolor[HTML]{067230}} \color[HTML]{F1F1F1} 15\% \\
Hier. ($x_{\mathrm{act}}$; $n_{x,\mathrm{grp}}$) & {\cellcolor[HTML]{56B567}} \color[HTML]{F1F1F1} 4 & {\cellcolor[HTML]{56B567}} \color[HTML]{F1F1F1} 4 & {\cellcolor[HTML]{8ED08B}} \color[HTML]{000000} 5 & {\cellcolor[HTML]{00441B}} \color[HTML]{F1F1F1} 1 & {\cellcolor[HTML]{05712F}} \color[HTML]{F1F1F1} 2 & {\cellcolor[HTML]{BCE4B5}} \color[HTML]{000000} 6 & {\cellcolor[HTML]{00441B}} \color[HTML]{F1F1F1} 1 & {\cellcolor[HTML]{00441B}} \color[HTML]{F1F1F1} 1 & {\cellcolor[HTML]{9DCAE1}} \color[HTML]{000000} 38\% & {\cellcolor[HTML]{6AAED6}} \color[HTML]{F1F1F1} 50\% & {\cellcolor[HTML]{DDF2D8}} \color[HTML]{000000} 84\% \\
Hier. (MRD) & {\cellcolor[HTML]{00441B}} \color[HTML]{F1F1F1} 1 & {\cellcolor[HTML]{05712F}} \color[HTML]{F1F1F1} 2 & {\cellcolor[HTML]{00441B}} \color[HTML]{F1F1F1} 1 & {\cellcolor[HTML]{56B567}} \color[HTML]{F1F1F1} 4 & {\cellcolor[HTML]{56B567}} \color[HTML]{F1F1F1} 4 & {\cellcolor[HTML]{2C944C}} \color[HTML]{F1F1F1} 3 & {\cellcolor[HTML]{8ED08B}} \color[HTML]{000000} 5 & {\cellcolor[HTML]{00441B}} \color[HTML]{F1F1F1} 1 & {\cellcolor[HTML]{9DCAE1}} \color[HTML]{000000} 38\% & {\cellcolor[HTML]{6AAED6}} \color[HTML]{F1F1F1} 50\% & {\cellcolor[HTML]{005622}} \color[HTML]{F1F1F1} 6\% \\
Hier. (MRD; $n_{\mathrm{act}}$) & {\cellcolor[HTML]{2C944C}} \color[HTML]{F1F1F1} 3 & {\cellcolor[HTML]{2C944C}} \color[HTML]{F1F1F1} 3 & {\cellcolor[HTML]{56B567}} \color[HTML]{F1F1F1} 4 & {\cellcolor[HTML]{2C944C}} \color[HTML]{F1F1F1} 3 & {\cellcolor[HTML]{2C944C}} \color[HTML]{F1F1F1} 3 & {\cellcolor[HTML]{8ED08B}} \color[HTML]{000000} 5 & {\cellcolor[HTML]{2C944C}} \color[HTML]{F1F1F1} 3 & {\cellcolor[HTML]{00441B}} \color[HTML]{F1F1F1} 1 & {\cellcolor[HTML]{DEEBF7}} \color[HTML]{000000} 12\% & {\cellcolor[HTML]{DEEBF7}} \color[HTML]{000000} 12\% & {\cellcolor[HTML]{238B45}} \color[HTML]{F1F1F1} 25\% \\
Hier. (MRD; $n_{x,\mathrm{grp}}$) & {\cellcolor[HTML]{8ED08B}} \color[HTML]{000000} 5 & {\cellcolor[HTML]{56B567}} \color[HTML]{F1F1F1} 4 & {\cellcolor[HTML]{8ED08B}} \color[HTML]{000000} 5 & {\cellcolor[HTML]{00441B}} \color[HTML]{F1F1F1} 1 & {\cellcolor[HTML]{05712F}} \color[HTML]{F1F1F1} 2 & {\cellcolor[HTML]{E1F3DC}} \color[HTML]{000000} 7 & {\cellcolor[HTML]{05712F}} \color[HTML]{F1F1F1} 2 & {\cellcolor[HTML]{00441B}} \color[HTML]{F1F1F1} 1 & {\cellcolor[HTML]{C6DBEF}} \color[HTML]{000000} 25\% & {\cellcolor[HTML]{6AAED6}} \color[HTML]{F1F1F1} 50\% & {\cellcolor[HTML]{F7FCF5}} \color[HTML]{000000} 118\% \\
\bottomrule
\end{tabular}

%% file: sampling_rank_sbo_mod.tex
\begin{tabular}{lcccccccccc}
\toprule
Strategy (grp.; wt.) & RCost & RWt & Rocket & SO GNC & FR GNC & MD GNC & Jet SM & Rank 1 & Rank $\leq$ 2 & Penalty \\
\midrule
Non-hier. & {\cellcolor[HTML]{00441B}} \color[HTML]{F1F1F1} 1 & {\cellcolor[HTML]{F7FCF5}} \color[HTML]{000000} 2 & {\cellcolor[HTML]{00441B}} \color[HTML]{F1F1F1} 1 & {\cellcolor[HTML]{F7FCF5}} \color[HTML]{000000} 2 & {\cellcolor[HTML]{00441B}} \color[HTML]{F1F1F1} 1 & {\cellcolor[HTML]{00441B}} \color[HTML]{F1F1F1} 1 & {\cellcolor[HTML]{00441B}} \color[HTML]{F1F1F1} 1 & {\cellcolor[HTML]{2B7BBA}} \color[HTML]{F1F1F1} 71\% & {\cellcolor[HTML]{08306B}} \color[HTML]{F1F1F1} 100\% & {\cellcolor[HTML]{94D390}} \color[HTML]{000000} 12\% \\
Hier. & {\cellcolor[HTML]{00441B}} \color[HTML]{F1F1F1} 1 & {\cellcolor[HTML]{00441B}} \color[HTML]{F1F1F1} 1 & {\cellcolor[HTML]{F7FCF5}} \color[HTML]{000000} 2 & {\cellcolor[HTML]{F7FCF5}} \color[HTML]{000000} 2 & {\cellcolor[HTML]{00441B}} \color[HTML]{F1F1F1} 1 & {\cellcolor[HTML]{00441B}} \color[HTML]{F1F1F1} 1 & {\cellcolor[HTML]{F7FCF5}} \color[HTML]{000000} 2 & {\cellcolor[HTML]{539ECD}} \color[HTML]{F1F1F1} 57\% & {\cellcolor[HTML]{08306B}} \color[HTML]{F1F1F1} 100\% & {\cellcolor[HTML]{86CC85}} \color[HTML]{000000} 11\% \\
Hier. ($n_{\mathrm{act}}$) & {\cellcolor[HTML]{F7FCF5}} \color[HTML]{000000} 2 & {\cellcolor[HTML]{00441B}} \color[HTML]{F1F1F1} 1 & {\cellcolor[HTML]{00441B}} \color[HTML]{F1F1F1} 1 & {\cellcolor[HTML]{00441B}} \color[HTML]{F1F1F1} 1 & {\cellcolor[HTML]{00441B}} \color[HTML]{F1F1F1} 1 & {\cellcolor[HTML]{F7FCF5}} \color[HTML]{000000} 2 & {\cellcolor[HTML]{00441B}} \color[HTML]{F1F1F1} 1 & {\cellcolor[HTML]{2B7BBA}} \color[HTML]{F1F1F1} 71\% & {\cellcolor[HTML]{08306B}} \color[HTML]{F1F1F1} 100\% & {\cellcolor[HTML]{E9F7E5}} \color[HTML]{000000} 18\% \\
\underline{Hier. ($x_{\mathrm{act}}$)} & {\cellcolor[HTML]{00441B}} \color[HTML]{F1F1F1} 1 & {\cellcolor[HTML]{00441B}} \color[HTML]{F1F1F1} 1 & {\cellcolor[HTML]{00441B}} \color[HTML]{F1F1F1} 1 & {\cellcolor[HTML]{00441B}} \color[HTML]{F1F1F1} 1 & {\cellcolor[HTML]{00441B}} \color[HTML]{F1F1F1} 1 & {\cellcolor[HTML]{00441B}} \color[HTML]{F1F1F1} 1 & {\cellcolor[HTML]{00441B}} \color[HTML]{F1F1F1} 1 & {\cellcolor[HTML]{08306B}} \color[HTML]{F1F1F1} \underline{100\%} & {\cellcolor[HTML]{08306B}} \color[HTML]{F1F1F1} \underline{100\%} & {\cellcolor[HTML]{00441B}} \color[HTML]{F1F1F1} \underline{0\%} \\
Hier. (MRD) & {\cellcolor[HTML]{00441B}} \color[HTML]{F1F1F1} 1 & {\cellcolor[HTML]{00441B}} \color[HTML]{F1F1F1} 1 & {\cellcolor[HTML]{F7FCF5}} \color[HTML]{000000} 2 & {\cellcolor[HTML]{00441B}} \color[HTML]{F1F1F1} 1 & {\cellcolor[HTML]{00441B}} \color[HTML]{F1F1F1} 1 & {\cellcolor[HTML]{F7FCF5}} \color[HTML]{000000} 2 & {\cellcolor[HTML]{00441B}} \color[HTML]{F1F1F1} 1 & {\cellcolor[HTML]{2B7BBA}} \color[HTML]{F1F1F1} 71\% & {\cellcolor[HTML]{08306B}} \color[HTML]{F1F1F1} 100\% & {\cellcolor[HTML]{006328}} \color[HTML]{F1F1F1} 2\% \\
\bottomrule
\end{tabular}

%% file: corr_rank_nsga2_mod.tex
\begin{tabular}{lllccccccccccc}
\toprule
Correction & Config & Sampling & RCost & RWt & Rocket & SO GNC & FR GNC & Wt GNC & MD GNC & Jet SM & Rank 1 & Rank $\leq$ 2 & Penalty \\
\midrule
\underline{Eager Any-select} & \underline{} & \underline{Hier. $x_{\mathrm{act}}$ } & {\cellcolor[HTML]{00441B}} \color[HTML]{F1F1F1} 1 & {\cellcolor[HTML]{00441B}} \color[HTML]{F1F1F1} 1 & {\cellcolor[HTML]{00441B}} \color[HTML]{F1F1F1} 1 & {\cellcolor[HTML]{00441B}} \color[HTML]{F1F1F1} 1 & {\cellcolor[HTML]{00441B}} \color[HTML]{F1F1F1} 1 & {\cellcolor[HTML]{C8E9C1}} \color[HTML]{000000} 4 & {\cellcolor[HTML]{238B45}} \color[HTML]{F1F1F1} 2 & {\cellcolor[HTML]{00441B}} \color[HTML]{F1F1F1} 1 & {\cellcolor[HTML]{2070B4}} \color[HTML]{F1F1F1} \underline{75\%} & {\cellcolor[HTML]{08509B}} \color[HTML]{F1F1F1} \underline{88\%} & {\cellcolor[HTML]{00441B}} \color[HTML]{F1F1F1} \underline{0\%} \\
Eager Greedy &  & Hier. $x_{\mathrm{act}}$  & {\cellcolor[HTML]{00441B}} \color[HTML]{F1F1F1} 1 & {\cellcolor[HTML]{238B45}} \color[HTML]{F1F1F1} 2 & {\cellcolor[HTML]{00441B}} \color[HTML]{F1F1F1} 1 & {\cellcolor[HTML]{75C477}} \color[HTML]{000000} 3 & {\cellcolor[HTML]{F7FCF5}} \color[HTML]{000000} 5 & {\cellcolor[HTML]{00441B}} \color[HTML]{F1F1F1} 1 & {\cellcolor[HTML]{C8E9C1}} \color[HTML]{000000} 4 & {\cellcolor[HTML]{00441B}} \color[HTML]{F1F1F1} 1 & {\cellcolor[HTML]{6AAED6}} \color[HTML]{F1F1F1} 50\% & {\cellcolor[HTML]{4191C6}} \color[HTML]{F1F1F1} 62\% & {\cellcolor[HTML]{208843}} \color[HTML]{F1F1F1} 5\% \\
Eager Similar & Manh & Hier. $x_{\mathrm{act}}$  & {\cellcolor[HTML]{00441B}} \color[HTML]{F1F1F1} 1 & {\cellcolor[HTML]{00441B}} \color[HTML]{F1F1F1} 1 & {\cellcolor[HTML]{00441B}} \color[HTML]{F1F1F1} 1 & {\cellcolor[HTML]{00441B}} \color[HTML]{F1F1F1} 1 & {\cellcolor[HTML]{75C477}} \color[HTML]{000000} 3 & {\cellcolor[HTML]{238B45}} \color[HTML]{F1F1F1} 2 & {\cellcolor[HTML]{00441B}} \color[HTML]{F1F1F1} 1 & {\cellcolor[HTML]{238B45}} \color[HTML]{F1F1F1} 2 & {\cellcolor[HTML]{4191C6}} \color[HTML]{F1F1F1} 62\% & {\cellcolor[HTML]{08509B}} \color[HTML]{F1F1F1} 88\% & {\cellcolor[HTML]{026F2E}} \color[HTML]{F1F1F1} 3\% \\
Eager Similar & Euc & Hier. $x_{\mathrm{act}}$  & {\cellcolor[HTML]{00441B}} \color[HTML]{F1F1F1} 1 & {\cellcolor[HTML]{00441B}} \color[HTML]{F1F1F1} 1 & {\cellcolor[HTML]{238B45}} \color[HTML]{F1F1F1} 2 & {\cellcolor[HTML]{238B45}} \color[HTML]{F1F1F1} 2 & {\cellcolor[HTML]{C8E9C1}} \color[HTML]{000000} 4 & {\cellcolor[HTML]{238B45}} \color[HTML]{F1F1F1} 2 & {\cellcolor[HTML]{00441B}} \color[HTML]{F1F1F1} 1 & {\cellcolor[HTML]{238B45}} \color[HTML]{F1F1F1} 2 & {\cellcolor[HTML]{9DCAE1}} \color[HTML]{000000} 38\% & {\cellcolor[HTML]{08509B}} \color[HTML]{F1F1F1} 88\% & {\cellcolor[HTML]{087432}} \color[HTML]{F1F1F1} 3\% \\
Lazy Any-select &  & Non-hier.   & {\cellcolor[HTML]{00441B}} \color[HTML]{F1F1F1} 1 & {\cellcolor[HTML]{00441B}} \color[HTML]{F1F1F1} 1 & {\cellcolor[HTML]{238B45}} \color[HTML]{F1F1F1} 2 & {\cellcolor[HTML]{C8E9C1}} \color[HTML]{000000} 4 & {\cellcolor[HTML]{75C477}} \color[HTML]{000000} 3 & {\cellcolor[HTML]{00441B}} \color[HTML]{F1F1F1} 1 & {\cellcolor[HTML]{75C477}} \color[HTML]{000000} 3 & {\cellcolor[HTML]{75C477}} \color[HTML]{000000} 3 & {\cellcolor[HTML]{9DCAE1}} \color[HTML]{000000} 38\% & {\cellcolor[HTML]{6AAED6}} \color[HTML]{F1F1F1} 50\% & {\cellcolor[HTML]{4EB264}} \color[HTML]{F1F1F1} 8\% \\
Lazy Similar & Depth-f. & Non-hier.   & {\cellcolor[HTML]{238B45}} \color[HTML]{F1F1F1} 2 & {\cellcolor[HTML]{00441B}} \color[HTML]{F1F1F1} 1 & {\cellcolor[HTML]{238B45}} \color[HTML]{F1F1F1} 2 & {\cellcolor[HTML]{75C477}} \color[HTML]{000000} 3 & {\cellcolor[HTML]{75C477}} \color[HTML]{000000} 3 & {\cellcolor[HTML]{238B45}} \color[HTML]{F1F1F1} 2 & {\cellcolor[HTML]{238B45}} \color[HTML]{F1F1F1} 2 & {\cellcolor[HTML]{238B45}} \color[HTML]{F1F1F1} 2 & {\cellcolor[HTML]{DEEBF7}} \color[HTML]{000000} 12\% & {\cellcolor[HTML]{2070B4}} \color[HTML]{F1F1F1} 75\% & {\cellcolor[HTML]{46AE60}} \color[HTML]{F1F1F1} 8\% \\
Lazy Similar & Manh; Dist-f. & Non-hier.   & {\cellcolor[HTML]{75C477}} \color[HTML]{000000} 3 & {\cellcolor[HTML]{00441B}} \color[HTML]{F1F1F1} 1 & {\cellcolor[HTML]{75C477}} \color[HTML]{000000} 3 & {\cellcolor[HTML]{75C477}} \color[HTML]{000000} 3 & {\cellcolor[HTML]{C8E9C1}} \color[HTML]{000000} 4 & {\cellcolor[HTML]{75C477}} \color[HTML]{000000} 3 & {\cellcolor[HTML]{00441B}} \color[HTML]{F1F1F1} 1 & {\cellcolor[HTML]{75C477}} \color[HTML]{000000} 3 & {\cellcolor[HTML]{C6DBEF}} \color[HTML]{000000} 25\% & {\cellcolor[HTML]{C6DBEF}} \color[HTML]{000000} 25\% & {\cellcolor[HTML]{C4E8BD}} \color[HTML]{000000} 15\% \\
Lazy Similar & Euc; Dist-f. & Non-hier.   & {\cellcolor[HTML]{75C477}} \color[HTML]{000000} 3 & {\cellcolor[HTML]{238B45}} \color[HTML]{F1F1F1} 2 & {\cellcolor[HTML]{75C477}} \color[HTML]{000000} 3 & {\cellcolor[HTML]{75C477}} \color[HTML]{000000} 3 & {\cellcolor[HTML]{75C477}} \color[HTML]{000000} 3 & {\cellcolor[HTML]{00441B}} \color[HTML]{F1F1F1} 1 & {\cellcolor[HTML]{00441B}} \color[HTML]{F1F1F1} 1 & {\cellcolor[HTML]{75C477}} \color[HTML]{000000} 3 & {\cellcolor[HTML]{C6DBEF}} \color[HTML]{000000} 25\% & {\cellcolor[HTML]{9DCAE1}} \color[HTML]{000000} 38\% & {\cellcolor[HTML]{95D391}} \color[HTML]{000000} 12\% \\
Problem-specific  &  & Hier. $x_{\mathrm{act}}$  & {\cellcolor[HTML]{238B45}} \color[HTML]{F1F1F1} 2 & {\cellcolor[HTML]{238B45}} \color[HTML]{F1F1F1} 2 & {\cellcolor[HTML]{238B45}} \color[HTML]{F1F1F1} 2 & {\cellcolor[HTML]{75C477}} \color[HTML]{000000} 3 & {\cellcolor[HTML]{C8E9C1}} \color[HTML]{000000} 4 & {\cellcolor[HTML]{00441B}} \color[HTML]{F1F1F1} 1 & {\cellcolor[HTML]{75C477}} \color[HTML]{000000} 3 & {\cellcolor[HTML]{238B45}} \color[HTML]{F1F1F1} 2 & {\cellcolor[HTML]{DEEBF7}} \color[HTML]{000000} 12\% & {\cellcolor[HTML]{4191C6}} \color[HTML]{F1F1F1} 62\% & {\cellcolor[HTML]{2D954D}} \color[HTML]{F1F1F1} 6\% \\
Problem-specific  &  & Non-hier.   & {\cellcolor[HTML]{75C477}} \color[HTML]{000000} 3 & {\cellcolor[HTML]{00441B}} \color[HTML]{F1F1F1} 1 & {\cellcolor[HTML]{238B45}} \color[HTML]{F1F1F1} 2 & {\cellcolor[HTML]{238B45}} \color[HTML]{F1F1F1} 2 & {\cellcolor[HTML]{238B45}} \color[HTML]{F1F1F1} 2 & {\cellcolor[HTML]{75C477}} \color[HTML]{000000} 3 & {\cellcolor[HTML]{238B45}} \color[HTML]{F1F1F1} 2 & {\cellcolor[HTML]{75C477}} \color[HTML]{000000} 3 & {\cellcolor[HTML]{DEEBF7}} \color[HTML]{000000} 12\% & {\cellcolor[HTML]{4191C6}} \color[HTML]{F1F1F1} 62\% & {\cellcolor[HTML]{4AAF61}} \color[HTML]{F1F1F1} 8\% \\
\bottomrule
\end{tabular}

%% file: corr_rank_sbo_mod.tex
\begin{tabular}{lllcccccccccc}
\toprule
Correction & Config & Sampling & RCost & RWt & Rocket & SO GNC & FR GNC & MD GNC & Jet SM & Rank 1 & Rank $\leq$ 2 & Penalty \\
\midrule
Eager Any-select &  & Hier. $x_{\mathrm{act}}$  & {\cellcolor[HTML]{00441B}} \color[HTML]{F1F1F1} 1 & {\cellcolor[HTML]{00441B}} \color[HTML]{F1F1F1} 1 & {\cellcolor[HTML]{00441B}} \color[HTML]{F1F1F1} 1 & {\cellcolor[HTML]{F7FCF5}} \color[HTML]{000000} 2 & {\cellcolor[HTML]{00441B}} \color[HTML]{F1F1F1} 1 & {\cellcolor[HTML]{F7FCF5}} \color[HTML]{000000} 2 & {\cellcolor[HTML]{00441B}} \color[HTML]{F1F1F1} 1 & {\cellcolor[HTML]{2B7BBA}} \color[HTML]{F1F1F1} 71\% & {\cellcolor[HTML]{08306B}} \color[HTML]{F1F1F1} 100\% & {\cellcolor[HTML]{006227}} \color[HTML]{F1F1F1} 2\% \\
Eager Similar & Manh & Hier. $x_{\mathrm{act}}$  & {\cellcolor[HTML]{00441B}} \color[HTML]{F1F1F1} 1 & {\cellcolor[HTML]{00441B}} \color[HTML]{F1F1F1} 1 & {\cellcolor[HTML]{F7FCF5}} \color[HTML]{000000} 2 & {\cellcolor[HTML]{F7FCF5}} \color[HTML]{000000} 2 & {\cellcolor[HTML]{00441B}} \color[HTML]{F1F1F1} 1 & {\cellcolor[HTML]{00441B}} \color[HTML]{F1F1F1} 1 & {\cellcolor[HTML]{00441B}} \color[HTML]{F1F1F1} 1 & {\cellcolor[HTML]{2B7BBA}} \color[HTML]{F1F1F1} 71\% & {\cellcolor[HTML]{08306B}} \color[HTML]{F1F1F1} 100\% & {\cellcolor[HTML]{0A7633}} \color[HTML]{F1F1F1} 3\% \\
Lazy Similar & Depth-f. & Non-hier.   & {\cellcolor[HTML]{00441B}} \color[HTML]{F1F1F1} 1 & {\cellcolor[HTML]{00441B}} \color[HTML]{F1F1F1} 1 & {\cellcolor[HTML]{00441B}} \color[HTML]{F1F1F1} 1 & {\cellcolor[HTML]{F7FCF5}} \color[HTML]{000000} 2 & {\cellcolor[HTML]{00441B}} \color[HTML]{F1F1F1} 1 & {\cellcolor[HTML]{F7FCF5}} \color[HTML]{000000} 2 & {\cellcolor[HTML]{F7FCF5}} \color[HTML]{000000} 2 & {\cellcolor[HTML]{539ECD}} \color[HTML]{F1F1F1} 57\% & {\cellcolor[HTML]{08306B}} \color[HTML]{F1F1F1} 100\% & {\cellcolor[HTML]{3EA75A}} \color[HTML]{F1F1F1} 7\% \\
\underline{Problem-specific } & \underline{} & \underline{Hier. $x_{\mathrm{act}}$ } & {\cellcolor[HTML]{00441B}} \color[HTML]{F1F1F1} 1 & {\cellcolor[HTML]{F7FCF5}} \color[HTML]{000000} 2 & {\cellcolor[HTML]{00441B}} \color[HTML]{F1F1F1} 1 & {\cellcolor[HTML]{00441B}} \color[HTML]{F1F1F1} 1 & {\cellcolor[HTML]{00441B}} \color[HTML]{F1F1F1} 1 & {\cellcolor[HTML]{00441B}} \color[HTML]{F1F1F1} 1 & {\cellcolor[HTML]{00441B}} \color[HTML]{F1F1F1} 1 & {\cellcolor[HTML]{0B559F}} \color[HTML]{F1F1F1} \underline{86\%} & {\cellcolor[HTML]{08306B}} \color[HTML]{F1F1F1} \underline{100\%} & {\cellcolor[HTML]{00441B}} \color[HTML]{F1F1F1} \underline{0\%} \\
Problem-specific  &  & Non-hier.   & {\cellcolor[HTML]{00441B}} \color[HTML]{F1F1F1} 1 & {\cellcolor[HTML]{F7FCF5}} \color[HTML]{000000} 2 & {\cellcolor[HTML]{00441B}} \color[HTML]{F1F1F1} 1 & {\cellcolor[HTML]{F7FCF5}} \color[HTML]{000000} 2 & {\cellcolor[HTML]{F7FCF5}} \color[HTML]{000000} 2 & {\cellcolor[HTML]{00441B}} \color[HTML]{F1F1F1} 1 & {\cellcolor[HTML]{00441B}} \color[HTML]{F1F1F1} 1 & {\cellcolor[HTML]{539ECD}} \color[HTML]{F1F1F1} 57\% & {\cellcolor[HTML]{08306B}} \color[HTML]{F1F1F1} 100\% & {\cellcolor[HTML]{94D390}} \color[HTML]{000000} 12\% \\
\bottomrule
\end{tabular}

%% file: hier_rank_nsga2.tex
\begin{tabular}{lccccccccccc}
\toprule
 & RCost & RWt & Rocket & SO GNC & FR GNC & Wt GNC & MD GNC & Jet SM & Rank 1 & Rank $\leq$ 2 & Penalty \\
\midrule
Naive & {\cellcolor[HTML]{37A055}} \color[HTML]{F1F1F1} 2 & {\cellcolor[HTML]{37A055}} \color[HTML]{F1F1F1} 2 & {\cellcolor[HTML]{AEDEA7}} \color[HTML]{000000} 3 & {\cellcolor[HTML]{00441B}} \color[HTML]{F1F1F1} 1 & {\cellcolor[HTML]{AEDEA7}} \color[HTML]{000000} 3 & {\cellcolor[HTML]{37A055}} \color[HTML]{F1F1F1} 2 & {\cellcolor[HTML]{00441B}} \color[HTML]{F1F1F1} 1 & {\cellcolor[HTML]{00441B}} \color[HTML]{F1F1F1} 1 & {\cellcolor[HTML]{9DCAE1}} \color[HTML]{000000} 38\% & {\cellcolor[HTML]{2070B4}} \color[HTML]{F1F1F1} 75\% & {\cellcolor[HTML]{00441B}} \color[HTML]{F1F1F1} -3\% \\
X out & {\cellcolor[HTML]{AEDEA7}} \color[HTML]{000000} 3 & {\cellcolor[HTML]{00441B}} \color[HTML]{F1F1F1} 1 & {\cellcolor[HTML]{37A055}} \color[HTML]{F1F1F1} 2 & {\cellcolor[HTML]{37A055}} \color[HTML]{F1F1F1} 2 & {\cellcolor[HTML]{00441B}} \color[HTML]{F1F1F1} 1 & {\cellcolor[HTML]{37A055}} \color[HTML]{F1F1F1} 2 & {\cellcolor[HTML]{AEDEA7}} \color[HTML]{000000} 3 & {\cellcolor[HTML]{37A055}} \color[HTML]{F1F1F1} 2 & {\cellcolor[HTML]{C6DBEF}} \color[HTML]{000000} 25\% & {\cellcolor[HTML]{2070B4}} \color[HTML]{F1F1F1} 75\% & {\cellcolor[HTML]{005522}} \color[HTML]{F1F1F1} 1\% \\
\underline{Repair} & {\cellcolor[HTML]{00441B}} \color[HTML]{F1F1F1} 1 & {\cellcolor[HTML]{00441B}} \color[HTML]{F1F1F1} 1 & {\cellcolor[HTML]{37A055}} \color[HTML]{F1F1F1} 2 & {\cellcolor[HTML]{37A055}} \color[HTML]{F1F1F1} 2 & {\cellcolor[HTML]{37A055}} \color[HTML]{F1F1F1} 2 & {\cellcolor[HTML]{37A055}} \color[HTML]{F1F1F1} 2 & {\cellcolor[HTML]{37A055}} \color[HTML]{F1F1F1} 2 & {\cellcolor[HTML]{37A055}} \color[HTML]{F1F1F1} 2 & {\cellcolor[HTML]{C6DBEF}} \color[HTML]{000000} \underline{25\%} & {\cellcolor[HTML]{08306B}} \color[HTML]{F1F1F1} \underline{100\%} & {\cellcolor[HTML]{00441B}} \color[HTML]{F1F1F1} \underline{0\%} \\
Activeness & {\cellcolor[HTML]{00441B}} \color[HTML]{F1F1F1} 1 & {\cellcolor[HTML]{00441B}} \color[HTML]{F1F1F1} 1 & {\cellcolor[HTML]{00441B}} \color[HTML]{F1F1F1} 1 & {\cellcolor[HTML]{AEDEA7}} \color[HTML]{000000} 3 & {\cellcolor[HTML]{F7FCF5}} \color[HTML]{000000} 4 & {\cellcolor[HTML]{00441B}} \color[HTML]{F1F1F1} 1 & {\cellcolor[HTML]{F7FCF5}} \color[HTML]{000000} 4 & {\cellcolor[HTML]{00441B}} \color[HTML]{F1F1F1} 1 & {\cellcolor[HTML]{4191C6}} \color[HTML]{F1F1F1} 62\% & {\cellcolor[HTML]{4191C6}} \color[HTML]{F1F1F1} 62\% & {\cellcolor[HTML]{05712F}} \color[HTML]{F1F1F1} 3\% \\
\bottomrule
\end{tabular}

%% file: hier_rank_sbo.tex
\begin{tabular}{lcccccccccc}
\toprule
 & RCost & RWt & Rocket & SO GNC & FR GNC & MD GNC & Jet SM & Rank 1 & Rank $\leq$ 2 & Penalty \\
\midrule
Naive & {\cellcolor[HTML]{37A055}} \color[HTML]{F1F1F1} 2 & {\cellcolor[HTML]{37A055}} \color[HTML]{F1F1F1} 2 & {\cellcolor[HTML]{37A055}} \color[HTML]{F1F1F1} 2 & {\cellcolor[HTML]{00441B}} \color[HTML]{F1F1F1} 1 & {\cellcolor[HTML]{00441B}} \color[HTML]{F1F1F1} 1 & {\cellcolor[HTML]{AEDEA7}} \color[HTML]{000000} 3 & {\cellcolor[HTML]{37A055}} \color[HTML]{F1F1F1} 2 & {\cellcolor[HTML]{BAD6EB}} \color[HTML]{000000} 29\% & {\cellcolor[HTML]{0B559F}} \color[HTML]{F1F1F1} 86\% & {\cellcolor[HTML]{F7FCF5}} \color[HTML]{000000} 134\% \\
X out & {\cellcolor[HTML]{37A055}} \color[HTML]{F1F1F1} 2 & {\cellcolor[HTML]{37A055}} \color[HTML]{F1F1F1} 2 & {\cellcolor[HTML]{37A055}} \color[HTML]{F1F1F1} 2 & {\cellcolor[HTML]{AEDEA7}} \color[HTML]{000000} 3 & {\cellcolor[HTML]{37A055}} \color[HTML]{F1F1F1} 2 & {\cellcolor[HTML]{F7FCF5}} \color[HTML]{000000} 4 & {\cellcolor[HTML]{37A055}} \color[HTML]{F1F1F1} 2 & {\cellcolor[HTML]{F7FBFF}} \color[HTML]{000000} 0\% & {\cellcolor[HTML]{2B7BBA}} \color[HTML]{F1F1F1} 71\% & {\cellcolor[HTML]{F7FCF5}} \color[HTML]{000000} 149\% \\
Repair & {\cellcolor[HTML]{37A055}} \color[HTML]{F1F1F1} 2 & {\cellcolor[HTML]{37A055}} \color[HTML]{F1F1F1} 2 & {\cellcolor[HTML]{00441B}} \color[HTML]{F1F1F1} 1 & {\cellcolor[HTML]{37A055}} \color[HTML]{F1F1F1} 2 & {\cellcolor[HTML]{00441B}} \color[HTML]{F1F1F1} 1 & {\cellcolor[HTML]{37A055}} \color[HTML]{F1F1F1} 2 & {\cellcolor[HTML]{37A055}} \color[HTML]{F1F1F1} 2 & {\cellcolor[HTML]{BAD6EB}} \color[HTML]{000000} 29\% & {\cellcolor[HTML]{08306B}} \color[HTML]{F1F1F1} 100\% & {\cellcolor[HTML]{E8F6E3}} \color[HTML]{000000} 89\% \\
\underline{Hier sampl.} & {\cellcolor[HTML]{00441B}} \color[HTML]{F1F1F1} 1 & {\cellcolor[HTML]{00441B}} \color[HTML]{F1F1F1} 1 & {\cellcolor[HTML]{00441B}} \color[HTML]{F1F1F1} 1 & {\cellcolor[HTML]{00441B}} \color[HTML]{F1F1F1} 1 & {\cellcolor[HTML]{00441B}} \color[HTML]{F1F1F1} 1 & {\cellcolor[HTML]{00441B}} \color[HTML]{F1F1F1} 1 & {\cellcolor[HTML]{00441B}} \color[HTML]{F1F1F1} 1 & {\cellcolor[HTML]{08306B}} \color[HTML]{F1F1F1} \underline{100\%} & {\cellcolor[HTML]{08306B}} \color[HTML]{F1F1F1} \underline{100\%} & {\cellcolor[HTML]{00441B}} \color[HTML]{F1F1F1} \underline{0\%} \\
Activeness & {\cellcolor[HTML]{00441B}} \color[HTML]{F1F1F1} 1 & {\cellcolor[HTML]{00441B}} \color[HTML]{F1F1F1} 1 & {\cellcolor[HTML]{00441B}} \color[HTML]{F1F1F1} 1 & {\cellcolor[HTML]{00441B}} \color[HTML]{F1F1F1} 1 & {\cellcolor[HTML]{00441B}} \color[HTML]{F1F1F1} 1 & {\cellcolor[HTML]{AEDEA7}} \color[HTML]{000000} 3 & {\cellcolor[HTML]{00441B}} \color[HTML]{F1F1F1} 1 & {\cellcolor[HTML]{0B559F}} \color[HTML]{F1F1F1} 86\% & {\cellcolor[HTML]{0B559F}} \color[HTML]{F1F1F1} 86\% & {\cellcolor[HTML]{157F3B}} \color[HTML]{F1F1F1} 20\% \\
\bottomrule
\end{tabular}

%% file: article.bbl

\begin{thebibliography}{128}
\ifx \bisbn   \undefined \def \bisbn  #1{ISBN #1}\fi
\ifx \binits  \undefined \def \binits#1{#1}\fi
\ifx \bauthor  \undefined \def \bauthor#1{#1}\fi
\ifx \batitle  \undefined \def \batitle#1{#1}\fi
\ifx \bjtitle  \undefined \def \bjtitle#1{#1}\fi
\ifx \bvolume  \undefined \def \bvolume#1{\textbf{#1}}\fi
\ifx \byear  \undefined \def \byear#1{#1}\fi
\ifx \bissue  \undefined \def \bissue#1{#1}\fi
\ifx \bfpage  \undefined \def \bfpage#1{#1}\fi
\ifx \blpage  \undefined \def \blpage #1{#1}\fi
\ifx \burl  \undefined \def \burl#1{\textsf{#1}}\fi
\ifx \doiurl  \undefined \def \doiurl#1{\url{https://doi.org/#1}}\fi
\ifx \betal  \undefined \def \betal{\textit{et al.}}\fi
\ifx \binstitute  \undefined \def \binstitute#1{#1}\fi
\ifx \binstitutionaled  \undefined \def \binstitutionaled#1{#1}\fi
\ifx \bctitle  \undefined \def \bctitle#1{#1}\fi
\ifx \beditor  \undefined \def \beditor#1{#1}\fi
\ifx \bpublisher  \undefined \def \bpublisher#1{#1}\fi
\ifx \bbtitle  \undefined \def \bbtitle#1{#1}\fi
\ifx \bedition  \undefined \def \bedition#1{#1}\fi
\ifx \bseriesno  \undefined \def \bseriesno#1{#1}\fi
\ifx \blocation  \undefined \def \blocation#1{#1}\fi
\ifx \bsertitle  \undefined \def \bsertitle#1{#1}\fi
\ifx \bsnm \undefined \def \bsnm#1{#1}\fi
\ifx \bsuffix \undefined \def \bsuffix#1{#1}\fi
\ifx \bparticle \undefined \def \bparticle#1{#1}\fi
\ifx \barticle \undefined \def \barticle#1{#1}\fi
\bibcommenthead
\ifx \bconfdate \undefined \def \bconfdate #1{#1}\fi
\ifx \botherref \undefined \def \botherref #1{#1}\fi
\ifx \url \undefined \def \url#1{\textsf{#1}}\fi
\ifx \bchapter \undefined \def \bchapter#1{#1}\fi
\ifx \bbook \undefined \def \bbook#1{#1}\fi
\ifx \bcomment \undefined \def \bcomment#1{#1}\fi
\ifx \oauthor \undefined \def \oauthor#1{#1}\fi
\ifx \citeauthoryear \undefined \def \citeauthoryear#1{#1}\fi
\ifx \endbibitem  \undefined \def \endbibitem {}\fi
\ifx \bconflocation  \undefined \def \bconflocation#1{#1}\fi
\ifx \arxivurl  \undefined \def \arxivurl#1{\textsf{#1}}\fi
\csname PreBibitemsHook\endcsname

\bibitem[\protect\citeauthoryear{Crawley et~al.}{2015}]{Crawley2015}
\begin{bbook}
\bauthor{\bsnm{Crawley}, \binits{E.}},
\bauthor{\bsnm{Cameron}, \binits{B.}},
\bauthor{\bsnm{Selva}, \binits{D.}}:
\bbtitle{System Architecture: Strategy and Product Development for Complex Systems}.
\bpublisher{Pearson Education},
\blocation{England}
(\byear{2015}).
\doiurl{10.1007/978-1-4020-4399-4}
\end{bbook}
\endbibitem

\bibitem[\protect\citeauthoryear{Chan et~al.}{2022}]{chan2022aircraft}
\begin{bchapter}
\bauthor{\bsnm{Chan}, \binits{A.}},
\bauthor{\bsnm{Pires}, \binits{A.F.}},
\bauthor{\bsnm{Polacsek}, \binits{T.}},
\bauthor{\bsnm{Roussel}, \binits{S.}}:
\bctitle{The aircraft and its manufacturing system: From early requirements to global design}.
In: \bbtitle{International Conference on Advanced Information Systems Engineering}
(\byear{2022}).
\doiurl{10.1007/978-3-031-07472-1_10}
\end{bchapter}
\endbibitem

\bibitem[\protect\citeauthoryear{Iacobucci}{2012}]{Iacobucci2012}
\begin{botherref}
\oauthor{\bsnm{Iacobucci}, \binits{J.V.}}:
Rapid architecture alternative modeling (raam): a framework for capability-based analysis of system of systems architectures.
PhD thesis,
Georgia Institute of Technology
(2012)
\end{botherref}
\endbibitem

\bibitem[\protect\citeauthoryear{McDermott et~al.}{2020}]{McDermott2020}
\begin{bchapter}
\bauthor{\bsnm{McDermott}, \binits{T.A.}},
\bauthor{\bsnm{Folds}, \binits{D.J.}},
\bauthor{\bsnm{Hallo}, \binits{L.}}:
\bctitle{Addressing cognitive bias in systems engineering teams}.
In: \bbtitle{30th Annual INCOSE International Symposium},
\bconflocation{Virtual Event}
(\byear{2020}).
\doiurl{10.1002/j.2334-5837.2020.00721.x}
\end{bchapter}
\endbibitem

\bibitem[\protect\citeauthoryear{Judt and Lawson}{2016}]{Judt2016}
\begin{barticle}
\bauthor{\bsnm{Judt}, \binits{D.M.}},
\bauthor{\bsnm{Lawson}, \binits{C.P.}}:
\batitle{Development of an automated aircraft subsystem architecture generation and analysis tool}.
\bjtitle{Engineering Computations}
\bvolume{33}(\bissue{5}),
\bfpage{1327}--\blpage{1352}
(\byear{2016})
\doiurl{10.1108/EC-02-2014-0033}
\end{barticle}
\endbibitem

\bibitem[\protect\citeauthoryear{Bussemaker et~al.}{2021}]{Bussemaker2021}
\begin{bchapter}
\bauthor{\bsnm{Bussemaker}, \binits{J.H.}},
\bauthor{\bsnm{Bartoli}, \binits{N.}},
\bauthor{\bsnm{Lefebvre}, \binits{T.}},
\bauthor{\bsnm{Ciampa}, \binits{P.D.}},
\bauthor{\bsnm{Nagel}, \binits{B.}}:
\bctitle{Effectiveness of surrogate-based optimization algorithms for system architecture optimization}.
In: \bbtitle{{AIAA} {AVIATION} 2021 {FORUM}},
\bconflocation{Virtual Event}
(\byear{2021}).
\doiurl{10.2514/6.2021-3095}
\end{bchapter}
\endbibitem

\bibitem[\protect\citeauthoryear{Czarnecki et~al.}{2012}]{Czarnecki2012}
\begin{bchapter}
\bauthor{\bsnm{Czarnecki}, \binits{K.}},
\bauthor{\bsnm{Grünbacher}, \binits{P.}},
\bauthor{\bsnm{Rabiser}, \binits{R.}},
\bauthor{\bsnm{Schmid}, \binits{K.}},
\bauthor{\bsnm{Wasowski}, \binits{A.}}:
\bctitle{Cool features and tough decisions}.
In: \bbtitle{Proceedings of the Sixth International Workshop on Variability Modeling of Software-Intensive Systems - {VaMoS} {\textquotesingle}12}.
\bpublisher{{ACM} Press},
\blocation{Leipzig, Germany}
(\byear{2012}).
\doiurl{10.1145/2110147.2110167}
\end{bchapter}
\endbibitem

\bibitem[\protect\citeauthoryear{Gedell and Johannesson}{2012}]{Gedell2012}
\begin{barticle}
\bauthor{\bsnm{Gedell}, \binits{S.}},
\bauthor{\bsnm{Johannesson}, \binits{H.}}:
\batitle{Design rationale and system description aspects in product platform design: Focusing reuse in the design lifecycle phase}.
\bjtitle{Concurrent Engineering}
\bvolume{21}(\bissue{1}),
\bfpage{39}--\blpage{53}
(\byear{2012})
\doiurl{10.1177/1063293x12469216}
\end{barticle}
\endbibitem

\bibitem[\protect\citeauthoryear{Mavris et~al.}{2008}]{Mavris2008}
\begin{bchapter}
\bauthor{\bsnm{Mavris}, \binits{D.}},
\bauthor{\bsnm{{de Tenorio}}, \binits{C.}},
\bauthor{\bsnm{Armstrong}, \binits{M.}}:
\bctitle{Methodology for aircraft system architecture definition}.
In: \bbtitle{46th AIAA Aerospace Sciences Meeting and Exhibit},
pp. \bfpage{1}--\blpage{14}.
\bpublisher{American Institute of Aeronautics and Astronautics},
\blocation{Reston, Virigina}
(\byear{2008}).
\doiurl{10.2514/6.2008-149}
\end{bchapter}
\endbibitem

\bibitem[\protect\citeauthoryear{Chakraborty and Mavris}{2016}]{Chakraborty2016}
\begin{bchapter}
\bauthor{\bsnm{Chakraborty}, \binits{I.}},
\bauthor{\bsnm{Mavris}, \binits{D.N.}}:
\bctitle{{Integrated Assessment of Aircraft and Novel Subsystem Architectures in Early Design}}.
In: \bbtitle{54th AIAA Aerospace Sciences Meeting},
vol. \bseriesno{54},
pp. \bfpage{1268}--\blpage{1282}.
\bpublisher{American Institute of Aeronautics and Astronautics},
\blocation{Reston, Virginia}
(\byear{2016}).
\doiurl{10.2514/6.2016-0215}
\end{bchapter}
\endbibitem

\bibitem[\protect\citeauthoryear{Simmons}{2008}]{Simmons2008}
\begin{botherref}
\oauthor{\bsnm{Simmons}, \binits{W.L.}}:
A framework for decision support in systems architecting.
PhD thesis,
Massachusetts Institute of Technology
(2008)
\end{botherref}
\endbibitem

\bibitem[\protect\citeauthoryear{Herber}{2020}]{Herber2020a}
\begin{bchapter}
\bauthor{\bsnm{Herber}, \binits{D.R.}}:
\bctitle{Enhancements to the perfect matching approach for graph enumeration-based engineering challenges}.
In: \bbtitle{Volume 11A: 46th Design Automation Conference ({DAC})}.
\bpublisher{American Society of Mechanical Engineers}, \blocation{???}
(\byear{2020}).
\doiurl{10.1115/detc2020-22774}
\end{bchapter}
\endbibitem

\bibitem[\protect\citeauthoryear{Bussemaker et~al.}{2024}]{Bussemaker2024adsg}
\begin{bchapter}
\bauthor{\bsnm{Bussemaker}, \binits{J.H.}},
\bauthor{\bsnm{Boggero}, \binits{L.}},
\bauthor{\bsnm{Nagel}, \binits{B.}}:
\bctitle{System architecture design space exploration: Integration with computational environments and efficient optimization}.
In: \bbtitle{{AIAA} {AVIATION} 2024 {FORUM}},
\bconflocation{Las Vegas, NV, USA}
(\byear{2024}).
\doiurl{10.2514/6.2024-4647}
\end{bchapter}
\endbibitem

\bibitem[\protect\citeauthoryear{Sobieszczanski-Sobieski et~al.}{2015}]{SobieszczanskiSobieski2015}
\begin{bbook}
\bauthor{\bsnm{Sobieszczanski-Sobieski}, \binits{J.}},
\bauthor{\bsnm{Morris}, \binits{A.}},
\bauthor{\bsnm{{van Tooren}}, \binits{M.J.L.}}:
\bbtitle{Multidisciplinary Design Optimization Supported by Knowledge Based Engineering},
pp. \bfpage{1}--\blpage{378}.
\bpublisher{John Wiley {\&} Sons, Ltd},
\blocation{West Sussex, UK}
(\byear{2015}).
\doiurl{10.1002/9781118897072}
\end{bbook}
\endbibitem

\bibitem[\protect\citeauthoryear{Chaudemar and de~Saqui-Sannes}{2021}]{Chaudemar2021}
\begin{bchapter}
\bauthor{\bsnm{Chaudemar}, \binits{J.-C.}},
\bauthor{\bsnm{Saqui-Sannes}, \binits{P.}}:
\bctitle{{MBSE} and {MDAO} for early validation of design decisions: a bibliography survey}.
\bpublisher{{IEEE}}, \blocation{???}
(\byear{2021}).
\doiurl{10.1109/syscon48628.2021.9447140}
\end{bchapter}
\endbibitem

\bibitem[\protect\citeauthoryear{Helle et~al.}{2022}]{Helle2022a}
\begin{bchapter}
\bauthor{\bsnm{Helle}, \binits{P.}},
\bauthor{\bsnm{Schramm}, \binits{G.}},
\bauthor{\bsnm{Klostermann}, \binits{S.}},
\bauthor{\bsnm{Feo-Arenis}, \binits{S.}}:
\bctitle{Enabling multidisciplinary-analysis of {SysML} models in a heterogeneous tool landscape using parametric analysis models}.
In: \bbtitle{The Complex Systems Deisgn {\&} Management Conference (CSD{\&}M 2022)}
(\byear{2022})
\end{bchapter}
\endbibitem

\bibitem[\protect\citeauthoryear{Bussemaker et~al.}{2022}]{Bussemaker2022}
\begin{bchapter}
\bauthor{\bsnm{Bussemaker}, \binits{J.H.}},
\bauthor{\bsnm{Boggero}, \binits{L.}},
\bauthor{\bsnm{Ciampa}, \binits{P.D.}}:
\bctitle{From system architecting to system design and optimization: A link between {MBSE} and {MDAO}}.
In: \bbtitle{32nd Annual INCOSE International Symposium},
\bconflocation{Detroit, MI, USA}
(\byear{2022}).
\doiurl{10.1002/iis2.12935}
\end{bchapter}
\endbibitem

\bibitem[\protect\citeauthoryear{Sonneveld et~al.}{2023}]{Sonneveld2023}
\begin{botherref}
\oauthor{\bsnm{Sonneveld}, \binits{J.S.}},
\oauthor{\bsnm{Berg}, \binits{T.}},
\oauthor{\bsnm{La~Rocca}, \binits{G.}},
\oauthor{\bsnm{Valencia-Ibáñez}, \binits{S.}},
\oauthor{\bsnm{Manen}, \binits{B.}},
\oauthor{\bsnm{Bruggeman}, \binits{A.M.R.M.}}:
Dynamic workflow generation applied to aircraft moveable architecture optimization
(2023)
\doiurl{10.13009/EUCASS2023-544}
\end{botherref}
\endbibitem

\bibitem[\protect\citeauthoryear{Bruggeman et~al.}{2024}]{Bruggeman2024}
\begin{bchapter}
\bauthor{\bsnm{Bruggeman}, \binits{A.}},
\bauthor{\bsnm{Nikitin}, \binits{M.}},
\bauthor{\bsnm{La~Rocca}, \binits{G.}},
\bauthor{\bsnm{Bergsma}, \binits{O.}}:
\bctitle{Model-based approach for the simultaneous design of airframe components and their production process using dynamic mdao workflows}.
In: \bbtitle{AIAA SCITECH 2024 Forum}.
\bpublisher{American Institute of Aeronautics and Astronautics},
\blocation{Orlando, FL, USA}
(\byear{2024}).
\doiurl{10.2514/6.2024-1530}
\end{bchapter}
\endbibitem

\bibitem[\protect\citeauthoryear{Garg et~al.}{2024}]{Garg2024mdo}
\begin{bchapter}
\bauthor{\bsnm{Garg}, \binits{S.}},
\bauthor{\bsnm{García~Sánchez}, \binits{R.}},
\bauthor{\bsnm{Bussemaker}, \binits{J.H.}},
\bauthor{\bsnm{Boggero}, \binits{L.}},
\bauthor{\bsnm{Nagel}, \binits{B.}}:
\bctitle{Dynamic formulation and excecution of {MDAO} workflows for architecture optimization}.
In: \bbtitle{{AIAA} {AVIATION} 2024 {FORUM}},
\bconflocation{Las Vegas, NV, USA}
(\byear{2024}).
\doiurl{10.2514/6.2024-4402}
\end{bchapter}
\endbibitem

\bibitem[\protect\citeauthoryear{Frank et~al.}{2016}]{Frank2016}
\begin{bchapter}
\bauthor{\bsnm{Frank}, \binits{C.P.}},
\bauthor{\bsnm{Marlier}, \binits{R.}},
\bauthor{\bsnm{Pinon-Fischer}, \binits{O.J.}},
\bauthor{\bsnm{Mavris}, \binits{D.N.}}:
\bctitle{An evolutionary multi-architecture multi-objective optimization algorithm for design space exploration}.
In: \bbtitle{57th AIAA/ASCE/AHS/ASC Structures, Structural Dynamics, and Materials Conference},
\bconflocation{Reston, Virginia},
pp. \bfpage{1}--\blpage{19}
(\byear{2016}).
\doiurl{10.2514/6.2016-0414}
\end{bchapter}
\endbibitem

\bibitem[\protect\citeauthoryear{Apaza and Selva}{2021}]{Apaza2021}
\begin{bchapter}
\bauthor{\bsnm{Apaza}, \binits{G.}},
\bauthor{\bsnm{Selva}, \binits{D.}}:
\bctitle{Automatic composition of encoding scheme and search operators in system architecture optimization}.
In: \bbtitle{41st Computers and Information in Engineering Conference ({CIE})}.
\bpublisher{American Society of Mechanical Engineers},
\blocation{Virtual}
(\byear{2021}).
\doiurl{10.1115/detc2021-71399}
\end{bchapter}
\endbibitem

\bibitem[\protect\citeauthoryear{Bussemaker et~al.}{2024}]{Bussemaker2024hc}
\begin{bchapter}
\bauthor{\bsnm{Bussemaker}, \binits{J.H.}},
\bauthor{\bsnm{Saves}, \binits{P.}},
\bauthor{\bsnm{Bartoli}, \binits{N.}},
\bauthor{\bsnm{Lefebvre}, \binits{T.}},
\bauthor{\bsnm{Nagel}, \binits{B.}}:
\bctitle{Surrogate-based optimization of system architectures subject to hidden constraints}.
In: \bbtitle{{AIAA} {AVIATION} 2024 {FORUM}},
\bconflocation{Las Vegas, NV, USA}
(\byear{2024}).
\doiurl{10.2514/6.2024-4401}
\end{bchapter}
\endbibitem

\bibitem[\protect\citeauthoryear{Feurer and Hutter}{2019}]{Feurer2019}
\begin{bchapter}
\bauthor{\bsnm{Feurer}, \binits{M.}},
\bauthor{\bsnm{Hutter}, \binits{F.}}:
\bctitle{Hyperparameter optimization}.
In: \bbtitle{Automated Machine Learning},
pp. \bfpage{3}--\blpage{33}.
\bpublisher{Springer},
\blocation{Switzerland}
(\byear{2019}).
\doiurl{10.1007/978-3-030-05318-5_1}
\end{bchapter}
\endbibitem

\bibitem[\protect\citeauthoryear{Bischl et~al.}{2023}]{Bischl2023}
\begin{botherref}
\oauthor{\bsnm{Bischl}, \binits{B.}},
\oauthor{\bsnm{Binder}, \binits{M.}},
\oauthor{\bsnm{Lang}, \binits{M.}},
\oauthor{\bsnm{Pielok}, \binits{T.}},
\oauthor{\bsnm{Richter}, \binits{J.}},
\oauthor{\bsnm{Coors}, \binits{S.}},
\oauthor{\bsnm{Thomas}, \binits{J.}},
\oauthor{\bsnm{Ullmann}, \binits{T.}},
\oauthor{\bsnm{Becker}, \binits{M.}},
\oauthor{\bsnm{Boulesteix}, \binits{A.}},
\oauthor{\bsnm{Deng}, \binits{D.}},
\oauthor{\bsnm{Lindauer}, \binits{M.}}:
Hyperparameter optimization: Foundations, algorithms, best practices, and open challenges.
WIREs Data Mining and Knowledge Discovery
\textbf{13}(2)
(2023)
\doiurl{10.1002/widm.1484}
\end{botherref}
\endbibitem

\bibitem[\protect\citeauthoryear{Martins and Ning}{2022}]{Martins2022}
\begin{bbook}
\bauthor{\bsnm{Martins}, \binits{J.R.R.A.}},
\bauthor{\bsnm{Ning}, \binits{A.}}:
\bbtitle{Engineering Design Optimization}.
\bpublisher{Cambridge University Press},
\blocation{Cambridge}
(\byear{2022}).
\burl{https://mdobook.github.io/}
\end{bbook}
\endbibitem

\bibitem[\protect\citeauthoryear{Saves et~al.}{2021}]{saves2021constrained}
\begin{bchapter}
\bauthor{\bsnm{Saves}, \binits{P.}},
\bauthor{\bsnm{Bartoli}, \binits{N.}},
\bauthor{\bsnm{Diouane}, \binits{Y.}},
\bauthor{\bsnm{Lefebvre}, \binits{T.}},
\bauthor{\bsnm{Morlier}, \binits{J.}},
\bauthor{\bsnm{David}, \binits{C.}},
\bauthor{},
\bauthor{\bsnm{{Nguyen Van}}, \binits{E.}},
\bauthor{\bsnm{Defoort}, \binits{S.}}:
\bctitle{Constrained bayesian optimization over mixed categorical variables, with application to aircraft design}.
In: \bbtitle{AeroBest 2021}
(\byear{2021}).
\burl{https://hal.science/hal-03346341v1/file/DTIS21090postprint.pdf}
\end{bchapter}
\endbibitem

\bibitem[\protect\citeauthoryear{Saves et~al.}{2023}]{Saves2023a}
\begin{barticle}
\bauthor{\bsnm{Saves}, \binits{P.}},
\bauthor{\bsnm{Diouane}, \binits{Y.}},
\bauthor{\bsnm{Bartoli}, \binits{N.}},
\bauthor{\bsnm{Lefebvre}, \binits{T.}},
\bauthor{\bsnm{Morlier}, \binits{J.}}:
\batitle{A mixed-categorical correlation kernel for gaussian process}.
\bjtitle{Neurocomputing}
\bvolume{550},
\bfpage{126472}
(\byear{2023})
\doiurl{10.1016/j.neucom.2023.126472}
\end{barticle}
\endbibitem

\bibitem[\protect\citeauthoryear{Bussemaker and Ciampa}{2022}]{Bussemaker2022c}
\begin{bchapter}
\bauthor{\bsnm{Bussemaker}, \binits{J.H.}},
\bauthor{\bsnm{Ciampa}, \binits{P.D.}}:
\bctitle{{MBSE} in architecture design space exploration}.
In: \beditor{\bsnm{Madni}, \binits{A.M.}},
\beditor{\bsnm{Augustine}, \binits{N.}},
\beditor{\bsnm{Sievers}, \binits{M.}} (eds.)
\bbtitle{Handbook of Model-Based Systems Engineering}.
\bpublisher{Springer},
\blocation{Switzerland}
(\byear{2022}).
\doiurl{10.1007/978-3-030-27486-3_36-1}
\end{bchapter}
\endbibitem

\bibitem[\protect\citeauthoryear{Pelamatti et~al.}{2020}]{Pelamatti2020a}
\begin{barticle}
\bauthor{\bsnm{Pelamatti}, \binits{J.}},
\bauthor{\bsnm{Brevault}, \binits{L.}},
\bauthor{\bsnm{Balesdent}, \binits{M.}},
\bauthor{\bsnm{Talbi}, \binits{E.}},
\bauthor{\bsnm{Guerin}, \binits{Y.}}:
\batitle{Bayesian optimization of variable-size design space problems}.
\bjtitle{Optimization and Engineering}
(\byear{2020})
\doiurl{10.1007/s11081-020-09520-z}
\end{barticle}
\endbibitem

\bibitem[\protect\citeauthoryear{Armstrong et~al.}{2008}]{Armstrong2008}
\begin{bchapter}
\bauthor{\bsnm{Armstrong}, \binits{M.}},
\bauthor{\bsnm{Tenorio}, \binits{C.}},
\bauthor{\bsnm{Garcia}, \binits{E.}},
\bauthor{\bsnm{Mavris}, \binits{D.}}:
\bctitle{Function based architecture design space definition and exploration}.
In: \bbtitle{26th Congress of International Council of the Aeronautical Sciences},
\bconflocation{Anchorage, Alaska, USA}
(\byear{2008}).
\doiurl{10.2514/6.2008-8928}
\end{bchapter}
\endbibitem

\bibitem[\protect\citeauthoryear{Zaefferer and Horn}{2018}]{Zaefferer2018a}
\begin{bchapter}
\bauthor{\bsnm{Zaefferer}, \binits{M.}},
\bauthor{\bsnm{Horn}, \binits{D.}}:
\bctitle{A first analysis of kernels for kriging-based optimization in hierarchical search spaces}.
In: \bbtitle{Parallel Problem Solving from Nature, PPSN XI}
vol. \bseriesno{1},
pp. \bfpage{399}--\blpage{410}.
\bpublisher{Springer},
\blocation{Berlin, Heidelberg}
(\byear{2018}).
\doiurl{10.1007/978-3-319-99259-4_32}
\end{bchapter}
\endbibitem

\bibitem[\protect\citeauthoryear{Hutter and Osborne}{2013}]{Hutter2013}
\begin{botherref}
\oauthor{\bsnm{Hutter}, \binits{F.}},
\oauthor{\bsnm{Osborne}, \binits{M.A.}}:
A kernel for hierarchical parameter spaces
(2013)
\doiurl{10.48550/arXiv.1310.5738}
\end{botherref}
\endbibitem

\bibitem[\protect\citeauthoryear{Bergstra et~al.}{2011}]{Bergstra2011}
\begin{bchapter}
\bauthor{\bsnm{Bergstra}, \binits{J.}},
\bauthor{\bsnm{Bardenet}, \binits{R.}},
\bauthor{\bsnm{Bengio}, \binits{Y.}},
\bauthor{\bsnm{K\'{e}gl}, \binits{B.}}:
\bctitle{Algorithms for hyper-parameter optimization}.
In: \bbtitle{Advances in Neural Information Processing Systems 24},
\bconflocation{Granada, Spain}
(\byear{2011})
\end{bchapter}
\endbibitem

\bibitem[\protect\citeauthoryear{Jenatton et~al.}{2017}]{Jenatton2017}
\begin{bchapter}
\bauthor{\bsnm{Jenatton}, \binits{R.}},
\bauthor{\bsnm{Archambeau}, \binits{C.}},
\bauthor{\bsnm{Gonz{\'a}lez}, \binits{J.}},
\bauthor{\bsnm{Seeger}, \binits{M.}}:
\bctitle{{B}ayesian optimization with tree-structured dependencies}.
In: \bbtitle{Proceedings of the 34th International Conference on Machine Learning},
\bconflocation{Sydney, Australia}
(\byear{2017})
\end{bchapter}
\endbibitem

\bibitem[\protect\citeauthoryear{Abdelkhalik}{2012}]{Abdelkhalik2012}
\begin{barticle}
\bauthor{\bsnm{Abdelkhalik}, \binits{O.}}:
\batitle{Hidden genes genetic optimization for variable-size design space problems}.
\bjtitle{Journal of Optimization Theory and Applications}
\bvolume{156}(\bissue{2}),
\bfpage{450}--\blpage{468}
(\byear{2012})
\doiurl{10.1007/s10957-012-0122-6}
\end{barticle}
\endbibitem

\bibitem[\protect\citeauthoryear{Talbi}{2024}]{Talbi2024}
\begin{botherref}
\oauthor{\bsnm{Talbi}, \binits{P.E.-G.}}:
Metaheuristics for (variable-size) mixed optimization problems: A unified taxonomy and survey
(2024)
\doiurl{10.48550/ARXIV.2401.03880}
\end{botherref}
\endbibitem

\bibitem[\protect\citeauthoryear{Levesque et~al.}{2017}]{Levesque2017}
\begin{bchapter}
\bauthor{\bsnm{Levesque}, \binits{J.-C.}},
\bauthor{\bsnm{Durand}, \binits{A.}},
\bauthor{\bsnm{Gagne}, \binits{C.}},
\bauthor{\bsnm{Sabourin}, \binits{R.}}:
\bctitle{Bayesian optimization for conditional hyperparameter spaces}.
In: \bbtitle{2017 International Joint Conference on Neural Networks ({IJCNN})}.
\bpublisher{{IEEE}},
\blocation{Anchorage, Alaska, USA}
(\byear{2017}).
\doiurl{10.1109/ijcnn.2017.7965867}
\end{bchapter}
\endbibitem

\bibitem[\protect\citeauthoryear{Selva}{2012}]{Selva2012}
\begin{botherref}
\oauthor{\bsnm{Selva}, \binits{D.}}:
Rule-based system architecting of earth observation satellite systems.
PhD thesis,
Massachusetts Institute of Technology, Dept. of Aeronautics and Astronautics
(2012).
\url{http://hdl.handle.net/1721.1/76089}
\end{botherref}
\endbibitem

\bibitem[\protect\citeauthoryear{Weilkiens et~al.}{2015}]{Weilkiens2015}
\begin{bbook}
\bauthor{\bsnm{Weilkiens}, \binits{T.}},
\bauthor{\bsnm{Lamm}, \binits{J.G.}},
\bauthor{\bsnm{Roth}, \binits{S.}}:
\bbtitle{Model-Based System Architecture}.
\bpublisher{Wiley John and Sons},
\blocation{Hoboken, NJ, USA}
(\byear{2015}).
\doiurl{10.1002/9781119051930}
\end{bbook}
\endbibitem

\bibitem[\protect\citeauthoryear{{Le Digabel} and Wild}{2023}]{LeDigabel2023}
\begin{barticle}
\bauthor{\bsnm{{Le Digabel}}, \binits{S.}},
\bauthor{\bsnm{Wild}, \binits{S.M.}}:
\batitle{A taxonomy of constraints in black-box simulation-based optimization}.
\bjtitle{Optimization and Engineering}
(\byear{2023})
\doiurl{10.1007/s11081-023-09839-3}
\end{barticle}
\endbibitem

\bibitem[\protect\citeauthoryear{Ying et~al.}{2019}]{Ying2019}
\begin{bchapter}
\bauthor{\bsnm{Ying}, \binits{C.}},
\bauthor{\bsnm{Klein}, \binits{A.}},
\bauthor{\bsnm{Christiansen}, \binits{E.}},
\bauthor{\bsnm{Real}, \binits{E.}},
\bauthor{\bsnm{Murphy}, \binits{K.}},
\bauthor{\bsnm{Hutter}, \binits{F.}}:
\bctitle{{NAS}-bench-101: Towards reproducible neural architecture search}.
In: \bbtitle{Proceedings of the 36th International Conference on Machine Learning}.
\bpublisher{PMLR},
\blocation{Long Beach, CA, USA}
(\byear{2019})
\end{bchapter}
\endbibitem

\bibitem[\protect\citeauthoryear{Salcedo-Sanz}{2009}]{SalcedoSanz2009}
\begin{barticle}
\bauthor{\bsnm{Salcedo-Sanz}, \binits{S.}}:
\batitle{A survey of repair methods used as constraint handling techniques in evolutionary algorithms}.
\bjtitle{Computer Science Review}
\bvolume{3}(\bissue{3}),
\bfpage{175}--\blpage{192}
(\byear{2009})
\doiurl{10.1016/j.cosrev.2009.07.001}
\end{barticle}
\endbibitem

\bibitem[\protect\citeauthoryear{Koch et~al.}{2015}]{Koch2015}
\begin{bchapter}
\bauthor{\bsnm{Koch}, \binits{P.}},
\bauthor{\bsnm{Bagheri}, \binits{S.}},
\bauthor{\bsnm{Konen}, \binits{W.}},
\bauthor{\bsnm{Foussette}, \binits{C.}},
\bauthor{\bsnm{Krause}, \binits{P.}},
\bauthor{\bsnm{Bäck}, \binits{T.}}:
\bctitle{A new repair method for constrained optimization}.
In: \bbtitle{Proceedings of the 2015 Annual Conference on Genetic and Evolutionary Computation}.
\bpublisher{ACM},
\blocation{Madrid, Spain}
(\byear{2015}).
\doiurl{10.1145/2739480.2754658}
\end{bchapter}
\endbibitem

\bibitem[\protect\citeauthoryear{Benavides et~al.}{2010}]{Benavides2010}
\begin{barticle}
\bauthor{\bsnm{Benavides}, \binits{D.}},
\bauthor{\bsnm{Segura}, \binits{S.}},
\bauthor{\bsnm{Ruiz-Cort{\'{e}}s}, \binits{A.}}:
\batitle{Automated analysis of feature models 20 years later: A literature review}.
\bjtitle{Information Systems}
\bvolume{35}(\bissue{6}),
\bfpage{615}--\blpage{636}
(\byear{2010})
\doiurl{10.1016/j.is.2010.01.001}
\end{barticle}
\endbibitem

\bibitem[\protect\citeauthoryear{Bussemaker et~al.}{2021}]{Bussemaker2021c}
\begin{bchapter}
\bauthor{\bsnm{Bussemaker}, \binits{J.H.}},
\bauthor{\bsnm{{De Smedt}}, \binits{T.}},
\bauthor{\bsnm{{La Rocca}}, \binits{G.}},
\bauthor{\bsnm{Ciampa}, \binits{P.D.}},
\bauthor{\bsnm{Nagel}, \binits{B.}}:
\bctitle{System architecture optimization: An open source multidisciplinary aircraft jet engine architecting problem}.
In: \bbtitle{{AIAA} {AVIATION} 2021 {FORUM}},
\bconflocation{Virtual Event}
(\byear{2021}).
\doiurl{10.2514/6.2021-3078}
\end{bchapter}
\endbibitem

\bibitem[\protect\citeauthoryear{Jones et~al.}{1998}]{Jones1998}
\begin{barticle}
\bauthor{\bsnm{Jones}, \binits{D.R.}},
\bauthor{\bsnm{Schonlau}, \binits{M.}},
\bauthor{\bsnm{Welch}, \binits{W.J.}}:
\batitle{Efficient global optimization of expensive black-box functions}.
\bjtitle{Journal of Global Optimization}
\bvolume{13},
\bfpage{455}--\blpage{492}
(\byear{1998})
\doiurl{10.1023/A:1008306431147}
\end{barticle}
\endbibitem

\bibitem[\protect\citeauthoryear{Miettinen}{1998}]{Miettinen1998}
\begin{bbook}
\bauthor{\bsnm{Miettinen}, \binits{K.}}:
\bbtitle{Nonlinear Multiobjective Optimization}.
\bpublisher{Springer},
\blocation{USA}
(\byear{1998}).
\doiurl{10.1007/978-1-4615-5563-6}
\end{bbook}
\endbibitem

\bibitem[\protect\citeauthoryear{Rudenko and Schoenauer}{2004}]{Rudenko2004}
\begin{bchapter}
\bauthor{\bsnm{Rudenko}, \binits{O.}},
\bauthor{\bsnm{Schoenauer}, \binits{M.}}:
\bctitle{A steady performance stopping criterion for {P}areto-based evolutionary algorithms}.
In: \bbtitle{6th International Multi-Objective Programming and Goal Programming Conference},
\bconflocation{Hammamet, Tunisia}
(\byear{2004})
\end{bchapter}
\endbibitem

\bibitem[\protect\citeauthoryear{Priem et~al.}{2020}]{SEGO-UTB}
\begin{barticle}
\bauthor{\bsnm{Priem}, \binits{R.}},
\bauthor{\bsnm{Bartoli}, \binits{N.}},
\bauthor{\bsnm{Diouane}, \binits{Y.}},
\bauthor{\bsnm{Sgueglia}, \binits{A.}}:
\batitle{Upper trust bound feasibility criterion for mixed constrained {Bayesian} optimization with application to aircraft design}.
\bjtitle{Aerospace Science and Technology}
\bvolume{105},
\bfpage{105980}
(\byear{2020})
\end{barticle}
\endbibitem

\bibitem[\protect\citeauthoryear{Forrester et~al.}{2006}]{Forrester2006}
\begin{barticle}
\bauthor{\bsnm{Forrester}, \binits{A.I.J.}},
\bauthor{\bsnm{S{\'{o}}bester}, \binits{A.}},
\bauthor{\bsnm{Keane}, \binits{A.J.}}:
\batitle{Optimization with missing data}.
\bjtitle{Proceedings of the Royal Society A: Mathematical, Physical and Engineering Sciences}
\bvolume{462}(\bissue{2067}),
\bfpage{935}--\blpage{945}
(\byear{2006})
\doiurl{10.1098/rspa.2005.1608}
\end{barticle}
\endbibitem

\bibitem[\protect\citeauthoryear{Müller and Day}{2019}]{Mueller2019}
\begin{barticle}
\bauthor{\bsnm{Müller}, \binits{J.}},
\bauthor{\bsnm{Day}, \binits{M.}}:
\batitle{Surrogate optimization of computationally expensive black-box problems with hidden constraints}.
\bjtitle{{INFORMS} Journal on Computing}
\bvolume{31}(\bissue{4}),
\bfpage{689}--\blpage{702}
(\byear{2019})
\doiurl{10.1287/ijoc.2018.0864}
\end{barticle}
\endbibitem

\bibitem[\protect\citeauthoryear{Krengel and Hepperle}{2022}]{Krengel2022}
\begin{bchapter}
\bauthor{\bsnm{Krengel}, \binits{M.D.}},
\bauthor{\bsnm{Hepperle}, \binits{M.}}:
\bctitle{Effects of wing elasticity and basic load alleviation on conceptual aircraft designs}.
In: \bbtitle{{AIAA} {SCITECH} 2022 Forum}
(\byear{2022}).
\doiurl{10.2514/6.2022-0126}
\end{bchapter}
\endbibitem

\bibitem[\protect\citeauthoryear{Yang and Shami}{2020}]{Yang2020}
\begin{barticle}
\bauthor{\bsnm{Yang}, \binits{L.}},
\bauthor{\bsnm{Shami}, \binits{A.}}:
\batitle{On hyperparameter optimization of machine learning algorithms: Theory and practice}.
\bjtitle{Neurocomputing}
\bvolume{415},
\bfpage{295}--\blpage{316}
(\byear{2020})
\doiurl{10.1016/j.neucom.2020.07.061}
\end{barticle}
\endbibitem

\bibitem[\protect\citeauthoryear{Jones and Martins}{2020}]{Jones2020}
\begin{barticle}
\bauthor{\bsnm{Jones}, \binits{D.R.}},
\bauthor{\bsnm{Martins}, \binits{J.R.R.A.}}:
\batitle{The {DIRECT} algorithm: 25 years later}.
\bjtitle{Journal of Global Optimization}
\bvolume{79}(\bissue{3}),
\bfpage{521}--\blpage{566}
(\byear{2020})
\doiurl{10.1007/s10898-020-00952-6}
\end{barticle}
\endbibitem

\bibitem[\protect\citeauthoryear{Locatelli and Schoen}{2021}]{Locatelli2021}
\begin{barticle}
\bauthor{\bsnm{Locatelli}, \binits{M.}},
\bauthor{\bsnm{Schoen}, \binits{F.}}:
\batitle{({G}lobal) optimization: Historical notes and recent developments}.
\bjtitle{{EURO} Journal on Computational Optimization}
\bvolume{9},
\bfpage{100012}
(\byear{2021})
\doiurl{10.1016/j.ejco.2021.100012}
\end{barticle}
\endbibitem

\bibitem[\protect\citeauthoryear{Glover and Kochenberger}{2003}]{Glover2003}
\begin{bbook}
\bauthor{\bsnm{Glover}, \binits{F.}},
\bauthor{\bsnm{Kochenberger}, \binits{G.}}:
\bbtitle{Handbook of Metaheuristics}
vol. \bseriesno{57},
pp. \bfpage{457}--\blpage{474474}
(\byear{2003}).
\doiurl{10.1007/b101874}
\end{bbook}
\endbibitem

\bibitem[\protect\citeauthoryear{Petrowski and Ben-Hamida}{2017}]{Petrowski2017}
\begin{bbook}
\bauthor{\bsnm{Petrowski}, \binits{A.}},
\bauthor{\bsnm{Ben-Hamida}, \binits{S.}}:
\bbtitle{Evolutionary Algorithms},
p. \bfpage{256}.
\bpublisher{John Wiley \& Sons},
\blocation{London, UK}
(\byear{2017})
\end{bbook}
\endbibitem

\bibitem[\protect\citeauthoryear{Hamano et~al.}{2022}]{Hamano2022}
\begin{bchapter}
\bauthor{\bsnm{Hamano}, \binits{R.}},
\bauthor{\bsnm{Saito}, \binits{S.}},
\bauthor{\bsnm{Nomura}, \binits{M.}},
\bauthor{\bsnm{Shirakawa}, \binits{S.}}:
\bctitle{{CMA}-{ES} with margin}.
In: \bbtitle{Proceedings of the Genetic and Evolutionary Computation Conference}.
\bpublisher{{ACM}},
\blocation{Boston, US}
(\byear{2022}).
\doiurl{10.1145/3512290.3528827}
\end{bchapter}
\endbibitem

\bibitem[\protect\citeauthoryear{Deb et~al.}{2002}]{Deb2002}
\begin{barticle}
\bauthor{\bsnm{Deb}, \binits{K.}},
\bauthor{\bsnm{Pratap}, \binits{A.}},
\bauthor{\bsnm{Agarwal}, \binits{S.}},
\bauthor{\bsnm{Meyarivan}, \binits{T.}}:
\batitle{{A fast and elitist multiobjective genetic algorithm: NSGA-II}}.
\bjtitle{IEEE Transactions on Evolutionary Computation}
\bvolume{6}(\bissue{2}),
\bfpage{182}--\blpage{197}
(\byear{2002})
\doiurl{10.1109/4235.996017}
\end{barticle}
\endbibitem

\bibitem[\protect\citeauthoryear{Nyew et~al.}{2015}]{Nyew2015}
\begin{barticle}
\bauthor{\bsnm{Nyew}, \binits{H.M.}},
\bauthor{\bsnm{Abdelkhalik}, \binits{O.}},
\bauthor{\bsnm{Onder}, \binits{N.}}:
\batitle{Structured-chromosome evolutionary algorithms for variable-size autonomous interplanetary trajectory planning optimization}.
\bjtitle{Journal of Aerospace Information Systems}
\bvolume{12}(\bissue{3}),
\bfpage{314}--\blpage{328}
(\byear{2015})
\doiurl{10.2514/1.i010272}
\end{barticle}
\endbibitem

\bibitem[\protect\citeauthoryear{Gratton and Vicente}{2014}]{Gratton2014}
\begin{barticle}
\bauthor{\bsnm{Gratton}, \binits{S.}},
\bauthor{\bsnm{Vicente}, \binits{L.N.}}:
\batitle{A merit function approach for direct search}.
\bjtitle{{SIAM} Journal on Optimization}
\bvolume{24}(\bissue{4}),
\bfpage{1980}--\blpage{1998}
(\byear{2014})
\doiurl{10.1137/130917661}
\end{barticle}
\endbibitem

\bibitem[\protect\citeauthoryear{Lopez-Herrejon et~al.}{2015}]{LopezHerrejon2015}
\begin{barticle}
\bauthor{\bsnm{Lopez-Herrejon}, \binits{R.E.}},
\bauthor{\bsnm{Linsbauer}, \binits{L.}},
\bauthor{\bsnm{Egyed}, \binits{A.}}:
\batitle{A systematic mapping study of search-based software engineering for software product lines}.
\bjtitle{Information and Software Technology}
\bvolume{61},
\bfpage{33}--\blpage{51}
(\byear{2015})
\doiurl{10.1016/j.infsof.2015.01.008}
\end{barticle}
\endbibitem

\bibitem[\protect\citeauthoryear{Buonanno}{2005}]{Buonanno2005}
\begin{botherref}
\oauthor{\bsnm{Buonanno}, \binits{M.A.}}:
A method for aircraft concept exploration using multicriteria interactive genetic algorithms.
PhD thesis,
Georgia Institute of Technology
(August 2005)
\end{botherref}
\endbibitem

\bibitem[\protect\citeauthoryear{Pate et~al.}{2012}]{Pate2012}
\begin{barticle}
\bauthor{\bsnm{Pate}, \binits{D.J.}},
\bauthor{\bsnm{Patterson}, \binits{M.D.}},
\bauthor{\bsnm{German}, \binits{B.J.}}:
\batitle{{Optimizing Families of Reconfigurable Aircraft for Multiple Missions}}.
\bjtitle{Journal of Aircraft}
\bvolume{49}(\bissue{6}),
\bfpage{1988}--\blpage{2000}
(\byear{2012})
\doiurl{10.2514/1.C031667}
\end{barticle}
\endbibitem

\bibitem[\protect\citeauthoryear{Frank et~al.}{2018}]{Frank2018}
\begin{barticle}
\bauthor{\bsnm{Frank}, \binits{C.P.}},
\bauthor{\bsnm{Marlier}, \binits{R.A.}},
\bauthor{\bsnm{Pinon-Fischer}, \binits{O.J.}},
\bauthor{\bsnm{Mavris}, \binits{D.N.}}:
\batitle{Evolutionary multi-objective multi-architecture design space exploration methodology}.
\bjtitle{Optimization and Engineering}
\bvolume{19}(\bissue{2}),
\bfpage{359}--\blpage{381}
(\byear{2018})
\doiurl{10.1007/s11081-018-9373-x}
\end{barticle}
\endbibitem

\bibitem[\protect\citeauthoryear{Chugh et~al.}{2019}]{Chugh2019}
\begin{barticle}
\bauthor{\bsnm{Chugh}, \binits{T.}},
\bauthor{\bsnm{Sindhya}, \binits{K.}},
\bauthor{\bsnm{Hakanen}, \binits{J.}},
\bauthor{\bsnm{Miettinen}, \binits{K.}}:
\batitle{{A survey on handling computationally expensive multiobjective optimization problems with evolutionary algorithms}}.
\bjtitle{Soft Computing}
\bvolume{23}(\bissue{9}),
\bfpage{3137}--\blpage{3166}
(\byear{2019})
\doiurl{10.1007/s00500-017-2965-0}
\end{barticle}
\endbibitem

\bibitem[\protect\citeauthoryear{Regis and Shoemaker}{2013}]{Regis2013}
\begin{barticle}
\bauthor{\bsnm{Regis}, \binits{R.G.}},
\bauthor{\bsnm{Shoemaker}, \binits{C.A.}}:
\batitle{Combining radial basis function surrogates and dynamic coordinate search in high-dimensional expensive black-box optimization}.
\bjtitle{Engineering Optimization}
\bvolume{45}(\bissue{5}),
\bfpage{529}--\blpage{555}
(\byear{2013})
\doiurl{10.1080/0305215x.2012.687731}
\end{barticle}
\endbibitem

\bibitem[\protect\citeauthoryear{Bagheri et~al.}{2017}]{Bagheri2017}
\begin{barticle}
\bauthor{\bsnm{Bagheri}, \binits{S.}},
\bauthor{\bsnm{Konen}, \binits{W.}},
\bauthor{\bsnm{Emmerich}, \binits{M.}},
\bauthor{\bsnm{Bäck}, \binits{T.}}:
\batitle{Self-adjusting parameter control for surrogate-assisted constrained optimization under limited budgets}.
\bjtitle{Applied Soft Computing}
\bvolume{61},
\bfpage{377}--\blpage{393}
(\byear{2017})
\doiurl{10.1016/j.asoc.2017.07.060}
\end{barticle}
\endbibitem

\bibitem[\protect\citeauthoryear{Garnett}{2023}]{Garnett2023}
\begin{bbook}
\bauthor{\bsnm{Garnett}, \binits{R.}}:
\bbtitle{Bayesian Optimization}.
\bpublisher{Cambridge University Press},
\blocation{Cambridge, UK}
(\byear{2023}).
\doiurl{10.1017/9781108348973}
\end{bbook}
\endbibitem

\bibitem[\protect\citeauthoryear{Schonlau et~al.}{1998}]{Schonlau1998}
\begin{bchapter}
\bauthor{\bsnm{Schonlau}, \binits{M.}},
\bauthor{\bsnm{Welch}, \binits{W.J.}},
\bauthor{\bsnm{Jones}, \binits{D.R.}}:
\bctitle{Global versus local search in constrained optimization of computer models}.
In: \bbtitle{Lecture Notes - Monograph Series},
pp. \bfpage{11}--\blpage{25}.
\bpublisher{Institute of Mathematical Statistics},
\blocation{online}
(\byear{1998}).
\doiurl{10.1214/lnms/1215456182}
\end{bchapter}
\endbibitem

\bibitem[\protect\citeauthoryear{Sasena}{2002}]{Sasena2002}
\begin{botherref}
\oauthor{\bsnm{Sasena}, \binits{M.J.}}:
{Flexibility and Efficiency Enhancements for Constrained Global Design Optimization with Kriging Approximations}.
PhD thesis,
University of Michigan
(2002)
\end{botherref}
\endbibitem

\bibitem[\protect\citeauthoryear{Knowles}{2006}]{Knowles2006}
\begin{barticle}
\bauthor{\bsnm{Knowles}, \binits{J.}}:
\batitle{{ParEGO: a hybrid algorithm with on-line landscape approximation for expensive multiobjective optimization problems}}.
\bjtitle{IEEE Transactions on Evolutionary Computation}
\bvolume{10}(\bissue{1}),
\bfpage{50}--\blpage{66}
(\byear{2006})
\doiurl{10.1109/TEVC.2005.851274}
\end{barticle}
\endbibitem

\bibitem[\protect\citeauthoryear{Rojas-Gonzalez and {Van Nieuwenhuyse}}{2019}]{RojasGonzalez2019}
\begin{barticle}
\bauthor{\bsnm{Rojas-Gonzalez}, \binits{S.}},
\bauthor{\bsnm{{Van Nieuwenhuyse}}, \binits{I.}}:
\batitle{A survey on {K}riging-based infill algorithms for multiobjective simulation optimization}.
\bjtitle{Computers {\&} Operations Research}
\bvolume{116},
\bfpage{104869}
(\byear{2019})
\doiurl{10.1016/j.cor.2019.104869}
\end{barticle}
\endbibitem

\bibitem[\protect\citeauthoryear{Garrido-Merch{\'{a}}n and Hern{\'{a}}ndez-Lobato}{2020}]{GarridoMerchan2020}
\begin{barticle}
\bauthor{\bsnm{Garrido-Merch{\'{a}}n}, \binits{E.C.}},
\bauthor{\bsnm{Hern{\'{a}}ndez-Lobato}, \binits{D.}}:
\batitle{Dealing with categorical and integer-valued variables in bayesian optimization with gaussian processes}.
\bjtitle{Neurocomputing}
\bvolume{380},
\bfpage{20}--\blpage{35}
(\byear{2020})
\doiurl{10.1016/j.neucom.2019.11.004}
\end{barticle}
\endbibitem

\bibitem[\protect\citeauthoryear{Daulton et~al.}{2022}]{Daulton2022}
\begin{botherref}
\oauthor{\bsnm{Daulton}, \binits{S.}},
\oauthor{\bsnm{Wan}, \binits{X.}},
\oauthor{\bsnm{Eriksson}, \binits{D.}},
\oauthor{\bsnm{Balandat}, \binits{M.}},
\oauthor{\bsnm{Osborne}, \binits{M.A.}},
\oauthor{\bsnm{Bakshy}, \binits{E.}}:
Bayesian optimization over discrete and mixed spaces via probabilistic reparameterization
(2022)
\doiurl{10.48550/ARXIV.2210.10199}
\end{botherref}
\endbibitem

\bibitem[\protect\citeauthoryear{Pelamatti et~al.}{2019}]{Pelamatti2019}
\begin{barticle}
\bauthor{\bsnm{Pelamatti}, \binits{J.}},
\bauthor{\bsnm{Brevault}, \binits{L.}},
\bauthor{\bsnm{Balesdent}, \binits{M.}},
\bauthor{\bsnm{Talbi}, \binits{E.}},
\bauthor{\bsnm{Guerin}, \binits{Y.}}:
\batitle{{Efficient global optimization of constrained mixed variable problems}}.
\bjtitle{Journal of Global Optimization}
\bvolume{73}(\bissue{3}),
\bfpage{583}--\blpage{613}
(\byear{2019})
\doiurl{10.1007/s10898-018-0715-1}
\end{barticle}
\endbibitem

\bibitem[\protect\citeauthoryear{Zuniga and Sinoquet}{2020}]{Munoz2020}
\begin{barticle}
\bauthor{\bsnm{Zuniga}, \binits{M.M.}},
\bauthor{\bsnm{Sinoquet}, \binits{D.}}:
\batitle{Global optimization for mixed categorical-continuous variables based on gaussian process models with a randomized categorical space exploration step}.
\bjtitle{{INFOR}: Information Systems and Operational Research}
\bvolume{58}(\bissue{2}),
\bfpage{310}--\blpage{341}
(\byear{2020})
\doiurl{10.1080/03155986.2020.1730677}
\end{barticle}
\endbibitem

\bibitem[\protect\citeauthoryear{Dreczkowski et~al.}{2023}]{Dreczkowski2023}
\begin{botherref}
\oauthor{\bsnm{Dreczkowski}, \binits{K.}},
\oauthor{\bsnm{Grosnit}, \binits{A.}},
\oauthor{\bsnm{Ammar}, \binits{H.B.}}:
Framework and benchmarks for combinatorial and mixed-variable bayesian optimization
(2023)
\doiurl{10.48550/ARXIV.2306.09803}
\end{botherref}
\endbibitem

\bibitem[\protect\citeauthoryear{Audet et~al.}{2023}]{Audet2022}
\begin{barticle}
\bauthor{\bsnm{Audet}, \binits{C.}},
\bauthor{\bsnm{Hall{\'e}-Hannan}, \binits{E.}},
\bauthor{\bsnm{{Le Digabel}}, \binits{S.}}:
\batitle{A general mathematical framework for constrained mixed-variable blackbox optimization problems with meta and categorical variables}.
\bjtitle{Operations Research Forum}
\bvolume{4},
\bfpage{1}--\blpage{37}
(\byear{2023})
\end{barticle}
\endbibitem

\bibitem[\protect\citeauthoryear{Horn et~al.}{2019}]{Horn2019}
\begin{bchapter}
\bauthor{\bsnm{Horn}, \binits{D.}},
\bauthor{\bsnm{Stork}, \binits{J.}},
\bauthor{\bsnm{Schü{\ss}ler}, \binits{N.}},
\bauthor{\bsnm{Zaefferer}, \binits{M.}}:
\bctitle{Surrogates for hierarchical search spaces}.
In: \bbtitle{Proceedings of the Genetic and Evolutionary Computation Conference}.
\bpublisher{{ACM}},
\blocation{Prague, CZ}
(\byear{2019}).
\doiurl{10.1145/3321707.3321765}
\end{bchapter}
\endbibitem

\bibitem[\protect\citeauthoryear{Lu et~al.}{2018}]{Lu2018}
\begin{bchapter}
\bauthor{\bsnm{Lu}, \binits{X.}},
\bauthor{\bsnm{Gonzalez}, \binits{J.}},
\bauthor{\bsnm{Dai}, \binits{Z.}},
\bauthor{\bsnm{Lawrence}, \binits{N.}}:
\bctitle{Structured variationally auto-encoded optimization}.
In: \bbtitle{Proceedings of the 35th International Conference on Machine Learning}.
\bpublisher{PMLR},
\blocation{Stockholm, SE}
(\byear{2018})
\end{bchapter}
\endbibitem

\bibitem[\protect\citeauthoryear{Saves et~al.}{2024}]{Saves2023SMT}
\begin{barticle}
\bauthor{\bsnm{Saves}, \binits{P.}},
\bauthor{\bsnm{Lafage}, \binits{R.}},
\bauthor{\bsnm{Bartoli}, \binits{N.}},
\bauthor{\bsnm{Diouane}, \binits{Y.}},
\bauthor{\bsnm{Bussemaker}, \binits{J.H.}},
\bauthor{\bsnm{Lefebvre}, \binits{T.}},
\bauthor{\bsnm{Hwang}, \binits{J.T.}},
\bauthor{\bsnm{Morlier}, \binits{J.}},
\bauthor{\bsnm{Martins}, \binits{J.R.R.A.}}:
\batitle{{SMT 2.0}: A surrogate modeling toolbox with a focus on hierarchical and mixed variables gaussian processes}.
\bjtitle{Advances in Engineering Software}
\bvolume{188},
\bfpage{103571}
(\byear{2024})
\doiurl{10.1016/j.advengsoft.2023.103571}
\end{barticle}
\endbibitem

\bibitem[\protect\citeauthoryear{Bouhlel et~al.}{2018}]{Bouhlel2018}
\begin{barticle}
\bauthor{\bsnm{Bouhlel}, \binits{M.A.}},
\bauthor{\bsnm{Bartoli}, \binits{N.}},
\bauthor{\bsnm{Regis}, \binits{R.G.}},
\bauthor{\bsnm{Otsmane}, \binits{A.}},
\bauthor{\bsnm{Morlier}, \binits{J.}}:
\batitle{Efficient global optimization for high-dimensional constrained problems by using the kriging models combined with the partial least squares method}.
\bjtitle{Engineering Optimization}
\bvolume{50}(\bissue{12}),
\bfpage{2038}--\blpage{2053}
(\byear{2018})
\doiurl{10.1080/0305215x.2017.1419344}
\end{barticle}
\endbibitem

\bibitem[\protect\citeauthoryear{Priem et~al.}{2023}]{Priem2023}
\begin{bchapter}
\bauthor{\bsnm{Priem}, \binits{R.}},
\bauthor{\bsnm{Bartoli}, \binits{N.}},
\bauthor{\bsnm{Diouane}, \binits{Y.}},
\bauthor{\bsnm{Dubreuil}, \binits{S.}},
\bauthor{\bsnm{Saves}, \binits{P.}}:
\bctitle{High-dimensional efficient global optimization using both random and supervised embeddings}.
In: \bbtitle{{AIAA} {AVIATION} 2023 Forum}.
\bpublisher{American Institute of Aeronautics and Astronautics},
\blocation{San Diego, CA, USA}
(\byear{2023}).
\doiurl{10.2514/6.2023-4448}
\end{bchapter}
\endbibitem

\bibitem[\protect\citeauthoryear{Saves et~al.}{2024}]{Saves2024}
\begin{botherref}
\oauthor{\bsnm{Saves}, \binits{P.}},
\oauthor{\bsnm{Diouane}, \binits{Y.}},
\oauthor{\bsnm{Bartoli}, \binits{N.}},
\oauthor{\bsnm{Lefebvre}, \binits{T.}},
\oauthor{\bsnm{Morlier}, \binits{J.}}:
High-dimensional mixed-categorical gaussian processes with application to multidisciplinary design optimization for a green aircraft.
Structural and Multidisciplinary Optimization
\textbf{67}(5)
(2024)
\doiurl{10.1007/s00158-024-03785-z}
\end{botherref}
\endbibitem

\bibitem[\protect\citeauthoryear{Hutter et~al.}{2011}]{Hutter2011}
\begin{bchapter}
\bauthor{\bsnm{Hutter}, \binits{F.}},
\bauthor{\bsnm{Hoos}, \binits{H.H.}},
\bauthor{\bsnm{Leyton-Brown}, \binits{K.}}:
\bctitle{Sequential model-based optimization for general algorithm configuration}.
In: \bbtitle{Lecture Notes in Computer Science},
pp. \bfpage{507}--\blpage{523}.
\bpublisher{Springer},
\blocation{Berlin Heidelberg}
(\byear{2011}).
\doiurl{10.1007/978-3-642-25566-3_40}
\end{bchapter}
\endbibitem

\bibitem[\protect\citeauthoryear{Lindauer et~al.}{2022}]{Lindauer2022}
\begin{barticle}
\bauthor{\bsnm{Lindauer}, \binits{M.}},
\bauthor{\bsnm{Eggensperger}, \binits{K.}},
\bauthor{\bsnm{Feurer}, \binits{M.}},
\bauthor{\bsnm{Biedenkapp}, \binits{A.}},
\bauthor{\bsnm{Deng}, \binits{D.}},
\bauthor{\bsnm{Benjamins}, \binits{C.}},
\bauthor{\bsnm{Ruhkopf}, \binits{T.}},
\bauthor{\bsnm{Sass}, \binits{R.}},
\bauthor{\bsnm{Hutter}, \binits{F.}}:
\batitle{{SMAC3}: A versatile {B}ayesian optimization package for hyperparameter optimization}.
\bjtitle{Journal of Machine Learning Research}
\bvolume{23}(\bissue{54}),
\bfpage{1}--\blpage{9}
(\byear{2022})
\end{barticle}
\endbibitem

\bibitem[\protect\citeauthoryear{Ozaki et~al.}{2022}]{Ozaki2022}
\begin{barticle}
\bauthor{\bsnm{Ozaki}, \binits{Y.}},
\bauthor{\bsnm{Tanigaki}, \binits{Y.}},
\bauthor{\bsnm{Watanabe}, \binits{S.}},
\bauthor{\bsnm{Nomura}, \binits{M.}},
\bauthor{\bsnm{Onishi}, \binits{M.}}:
\batitle{Multiobjective tree-structured parzen estimator}.
\bjtitle{Journal of Artificial Intelligence Research}
\bvolume{73},
\bfpage{1209}--\blpage{1250}
(\byear{2022})
\doiurl{10.1613/jair.1.13188}
\end{barticle}
\endbibitem

\bibitem[\protect\citeauthoryear{Eggensperger et~al.}{2015}]{Eggensperger2015}
\begin{botherref}
\oauthor{\bsnm{Eggensperger}, \binits{K.}},
\oauthor{\bsnm{Hutter}, \binits{F.}},
\oauthor{\bsnm{Hoos}, \binits{H.}},
\oauthor{\bsnm{Leyton-Brown}, \binits{K.}}:
Efficient benchmarking of hyperparameter optimizers via surrogates.
Proceedings of the {AAAI} Conference on Artificial Intelligence
\textbf{29}(1)
(2015)
\doiurl{10.1609/aaai.v29i1.9375}
\end{botherref}
\endbibitem

\bibitem[\protect\citeauthoryear{Gamot et~al.}{2023}]{gamot2023hidden}
\begin{barticle}
\bauthor{\bsnm{Gamot}, \binits{J.}},
\bauthor{\bsnm{Balesdent}, \binits{M.}},
\bauthor{\bsnm{Tremolet}, \binits{A.}},
\bauthor{\bsnm{Wuilbercq}, \binits{R.}},
\bauthor{\bsnm{Melab}, \binits{N.}},
\bauthor{\bsnm{Talbi}, \binits{E.-G.}}:
\batitle{Hidden-variables genetic algorithm for variable-size design space optimal layout problems with application to aerospace vehicles}.
\bjtitle{Engineering Applications of Artificial Intelligence}
\bvolume{121},
\bfpage{105941}
(\byear{2023})
\end{barticle}
\endbibitem

\bibitem[\protect\citeauthoryear{Greenhill et~al.}{2020}]{Greenhill2020}
\begin{barticle}
\bauthor{\bsnm{Greenhill}, \binits{S.}},
\bauthor{\bsnm{Rana}, \binits{S.}},
\bauthor{\bsnm{Gupta}, \binits{S.}},
\bauthor{\bsnm{Vellanki}, \binits{P.}},
\bauthor{\bsnm{Venkatesh}, \binits{S.}}:
\batitle{Bayesian optimization for adaptive experimental design: A review}.
\bjtitle{{IEEE} Access}
\bvolume{8},
\bfpage{13937}--\blpage{13948}
(\byear{2020})
\doiurl{10.1109/access.2020.2966228}
\end{barticle}
\endbibitem

\bibitem[\protect\citeauthoryear{Calandra et~al.}{2016}]{Calandra2016}
\begin{bchapter}
\bauthor{\bsnm{Calandra}, \binits{R.}},
\bauthor{\bsnm{Peters}, \binits{J.}},
\bauthor{\bsnm{Rasmussen}, \binits{C.E.}},
\bauthor{\bsnm{Deisenroth}, \binits{M.P.}}:
\bctitle{Manifold gaussian processes for regression}.
In: \bbtitle{2016 International Joint Conference on Neural Networks ({IJCNN})}.
\bpublisher{{IEEE}},
\blocation{Vancouver, CA}
(\byear{2016}).
\doiurl{10.1109/ijcnn.2016.7727626}
\end{bchapter}
\endbibitem

\bibitem[\protect\citeauthoryear{Bouhlel et~al.}{2016}]{Bouhlel2016}
\begin{barticle}
\bauthor{\bsnm{Bouhlel}, \binits{M.A.}},
\bauthor{\bsnm{Bartoli}, \binits{N.}},
\bauthor{\bsnm{Otsmane}, \binits{A.}},
\bauthor{\bsnm{Morlier}, \binits{J.}}:
\batitle{Improving {K}riging surrogates of high-dimensional design models by {Partial Least Squares} dimension reduction}.
\bjtitle{Structural and Multidisciplinary Optimization}
\bvolume{53}(\bissue{5}),
\bfpage{935}--\blpage{952}
(\byear{2016})
\doiurl{10.1007/s00158-015-1395-9}
\end{barticle}
\endbibitem

\bibitem[\protect\citeauthoryear{Saves et~al.}{2022}]{Saves2022b}
\begin{bchapter}
\bauthor{\bsnm{Saves}, \binits{P.}},
\bauthor{\bsnm{Bartoli}, \binits{N.}},
\bauthor{\bsnm{Diouane}, \binits{Y.}},
\bauthor{\bsnm{Lefebvre}, \binits{T.}},
\bauthor{\bsnm{Morlier}, \binits{J.}},
\bauthor{\bsnm{David}, \binits{C.}},
\bauthor{\bsnm{Van}, \binits{E.N.}},
\bauthor{\bsnm{Defoort}, \binits{S.}}:
\bctitle{Multidisciplinary design optimization with mixed categorical variables for aircraft design}.
In: \bbtitle{{AIAA} {SCITECH} 2022 Forum}.
\bpublisher{American Institute of Aeronautics and Astronautics},
\blocation{San Diego, CA, USA}
(\byear{2022}).
\doiurl{10.2514/6.2022-0082}
\end{bchapter}
\endbibitem

\bibitem[\protect\citeauthoryear{Charayron et~al.}{2023}]{Charayron2023}
\begin{bchapter}
\bauthor{\bsnm{Charayron}, \binits{R.}},
\bauthor{\bsnm{Lefebvre}, \binits{T.}},
\bauthor{\bsnm{Bartoli}, \binits{N.}},
\bauthor{\bsnm{Morlier}, \binits{J.}}:
\bctitle{Multi-fidelity {B}ayesian optimization strategy applied to overall drone design}.
In: \bbtitle{{AIAA} {SCITECH} 2023 Forum}.
\bpublisher{American Institute of Aeronautics and Astronautics},
\blocation{National Harbor, MD, USA}
(\byear{2023}).
\doiurl{10.2514/6.2023-2366}
\end{bchapter}
\endbibitem

\bibitem[\protect\citeauthoryear{Bouhlel et~al.}{2019}]{Bouhlel2019}
\begin{barticle}
\bauthor{\bsnm{Bouhlel}, \binits{M.A.}},
\bauthor{\bsnm{Hwang}, \binits{J.T.}},
\bauthor{\bsnm{Bartoli}, \binits{N.}},
\bauthor{\bsnm{Lafage}, \binits{R.}},
\bauthor{\bsnm{Morlier}, \binits{J.}},
\bauthor{\bsnm{Martins}, \binits{J.R.R.A.}}:
\batitle{A {P}ython surrogate modeling framework with derivatives}.
\bjtitle{Advances in Engineering Software}
\bvolume{135},
\bfpage{102662}
(\byear{2019})
\doiurl{10.1016/j.advengsoft.2019.03.005}
\end{barticle}
\endbibitem

\bibitem[\protect\citeauthoryear{Lyu et~al.}{2018}]{Lyu2018}
\begin{bchapter}
\bauthor{\bsnm{Lyu}, \binits{W.}},
\bauthor{\bsnm{Yang}, \binits{F.}},
\bauthor{\bsnm{Yan}, \binits{C.}},
\bauthor{\bsnm{Zhou}, \binits{D.}},
\bauthor{\bsnm{Zeng}, \binits{X.}}:
\bctitle{Batch {B}ayesian optimization via multi-objective acquisition ensemble for automated analog circuit design}.
In: \bbtitle{Proceedings of the 35th International Conference on Machine Learning}.
\bpublisher{PMLR},
\blocation{Stockholm, SE}
(\byear{2018})
\end{bchapter}
\endbibitem

\bibitem[\protect\citeauthoryear{Cowen-Rivers et~al.}{2020}]{CowenRivers2020}
\begin{botherref}
\oauthor{\bsnm{Cowen-Rivers}, \binits{A.I.}},
\oauthor{\bsnm{Lyu}, \binits{W.}},
\oauthor{\bsnm{Tutunov}, \binits{R.}},
\oauthor{\bsnm{Wang}, \binits{Z.}},
\oauthor{\bsnm{Grosnit}, \binits{A.}},
\oauthor{\bsnm{Griffiths}, \binits{R.R.}},
\oauthor{\bsnm{Maraval}, \binits{A.M.}},
\oauthor{\bsnm{Jianye}, \binits{H.}},
\oauthor{\bsnm{Wang}, \binits{J.}},
\oauthor{\bsnm{Peters}, \binits{J.}},
\oauthor{\bsnm{Ammar}, \binits{H.B.}}:
{HEBO}: Pushing the limits of sample-efficient hyperparameter optimisation
(2020)
\doiurl{10.48550/ARXIV.2012.03826}
\end{botherref}
\endbibitem

\bibitem[\protect\citeauthoryear{Ginsbourger et~al.}{2010}]{Ginsbourger2010}
\begin{bchapter}
\bauthor{\bsnm{Ginsbourger}, \binits{D.}},
\bauthor{\bsnm{Riche}, \binits{R.L.}},
\bauthor{\bsnm{Carraro}, \binits{L.}}:
\bctitle{Kriging is well-suited to parallelize optimization}.
In: \bbtitle{Computational Intelligence in Expensive Optimization Problems},
pp. \bfpage{131}--\blpage{162}.
\bpublisher{Springer},
\blocation{Berlin Heidelberg}
(\byear{2010}).
\doiurl{10.1007/978-3-642-10701-6_6}
\end{bchapter}
\endbibitem

\bibitem[\protect\citeauthoryear{Cox and John}{1992}]{Cox1992}
\begin{bchapter}
\bauthor{\bsnm{Cox}, \binits{D.D.}},
\bauthor{\bsnm{John}, \binits{S.}}:
\bctitle{A statistical method for global optimization}.
In: \bbtitle{1992 {IEEE} International Conference on Systems, Man, and Cybernetics}.
\bpublisher{{IEEE}},
\blocation{Chicago, IL, USA}
(\byear{1992}).
\doiurl{10.1109/icsmc.1992.271617}
\end{bchapter}
\endbibitem

\bibitem[\protect\citeauthoryear{Hawe and Sykulski}{2007}]{Hawe2007}
\begin{barticle}
\bauthor{\bsnm{Hawe}, \binits{G.I.}},
\bauthor{\bsnm{Sykulski}, \binits{J.K.}}:
\batitle{An enhanced probability of improvement utility function for locating pareto optimal solutions}.
\bjtitle{16th Conference on the Computation of Electromagnetic Fields, COMPUMAG, Aachen, Germany}
\bvolume{3},
\bfpage{965}--\blpage{966}
(\byear{2007})
\end{barticle}
\endbibitem

\bibitem[\protect\citeauthoryear{Rahat et~al.}{2017}]{Rahat2017}
\begin{bchapter}
\bauthor{\bsnm{Rahat}, \binits{A.A.M.}},
\bauthor{\bsnm{Everson}, \binits{R.M.}},
\bauthor{\bsnm{Fieldsend}, \binits{J.E.}}:
\bctitle{{Alternative infill strategies for expensive multi-objective optimisation}}.
In: \bbtitle{Proceedings of the Genetic and Evolutionary Computation Conference on - GECCO '17},
pp. \bfpage{873}--\blpage{880}.
\bpublisher{ACM Press},
\blocation{New York, USA}
(\byear{2017}).
\doiurl{10.1145/3071178.3071276}
\end{bchapter}
\endbibitem

\bibitem[\protect\citeauthoryear{Sohst et~al.}{2021}]{Sohst2021}
\begin{botherref}
\oauthor{\bsnm{Sohst}, \binits{M.}},
\oauthor{\bsnm{Afonso}, \binits{F.}},
\oauthor{\bsnm{Suleman}, \binits{A.}}:
Surrogate-based optimization based on the probability of feasibility.
Structural and Multidisciplinary Optimization
\textbf{65}(1)
(2021)
\doiurl{10.1007/s00158-021-03134-4}
\end{botherref}
\endbibitem

\bibitem[\protect\citeauthoryear{Bussemaker}{2023}]{Bussemaker2023a}
\begin{barticle}
\bauthor{\bsnm{Bussemaker}, \binits{J.H.}}:
\batitle{{SBArchOpt}: Surrogate-based architecture optimization}.
\bjtitle{Journal of Open Source Software}
\bvolume{8}(\bissue{89}),
\bfpage{5564}
(\byear{2023})
\doiurl{10.21105/joss.05564}
\end{barticle}
\endbibitem

\bibitem[\protect\citeauthoryear{Blank and Deb}{2020}]{Blank2020}
\begin{barticle}
\bauthor{\bsnm{Blank}, \binits{J.}},
\bauthor{\bsnm{Deb}, \binits{K.}}:
\batitle{Pymoo: Multi-objective optimization in python}.
\bjtitle{{IEEE} Access}
\bvolume{8},
\bfpage{89497}--\blpage{89509}
(\byear{2020})
\doiurl{10.1109/access.2020.2990567}
\end{barticle}
\endbibitem

\bibitem[\protect\citeauthoryear{Lindauer et~al.}{2019}]{Lindauer2019}
\begin{botherref}
\oauthor{\bsnm{Lindauer}, \binits{M.}},
\oauthor{\bsnm{Eggensperger}, \binits{K.}},
\oauthor{\bsnm{Feurer}, \binits{M.}},
\oauthor{\bsnm{Biedenkapp}, \binits{A.}},
\oauthor{\bsnm{Marben}, \binits{J.}},
\oauthor{\bsnm{Müller}, \binits{P.}},
\oauthor{\bsnm{Hutter}, \binits{F.}}:
{BOAH}: A tool suite for multi-fidelity {B}ayesian optimization \& analysis of hyperparameters
(2019)
\doiurl{10.48550/ARXIV.1908.06756}
\end{botherref}
\endbibitem

\bibitem[\protect\citeauthoryear{Balandat et~al.}{2019}]{Balandat2019}
\begin{barticle}
\bauthor{\bsnm{Balandat}, \binits{M.}},
\bauthor{\bsnm{Karrer}, \binits{B.}},
\bauthor{\bsnm{Jiang}, \binits{D.R.}},
\bauthor{\bsnm{Daulton}, \binits{S.}},
\bauthor{\bsnm{Letham}, \binits{B.}},
\bauthor{\bsnm{Wilson}, \binits{A.G.}},
\bauthor{\bsnm{Bakshy}, \binits{E.}}:
\batitle{{BoTorch}: A framework for efficient monte-carlo {B}ayesian optimization}.
\bjtitle{Advances in Neural Information Processing Systems 33}
(\byear{2019})
\doiurl{10.48550/ARXIV.1910.06403}
\end{barticle}
\endbibitem

\bibitem[\protect\citeauthoryear{Picheny et~al.}{2023}]{Picheny2023}
\begin{botherref}
\oauthor{\bsnm{Picheny}, \binits{V.}},
\oauthor{\bsnm{Berkeley}, \binits{J.}},
\oauthor{\bsnm{Moss}, \binits{H.B.}},
\oauthor{\bsnm{Stojic}, \binits{H.}},
\oauthor{\bsnm{Granta}, \binits{U.}},
\oauthor{\bsnm{Ober}, \binits{S.W.}},
\oauthor{\bsnm{Artemev}, \binits{A.}},
\oauthor{\bsnm{Ghani}, \binits{K.}},
\oauthor{\bsnm{Goodall}, \binits{A.}},
\oauthor{\bsnm{Paleyes}, \binits{A.}},
\oauthor{\bsnm{Vakili}, \binits{S.}},
\oauthor{\bsnm{Pascual-Diaz}, \binits{S.}},
\oauthor{\bsnm{Markou}, \binits{S.}},
\oauthor{\bsnm{Qing}, \binits{J.}},
\oauthor{\bsnm{Loka}, \binits{N.R.B.S.}},
\oauthor{\bsnm{Couckuyt}, \binits{I.}}:
{Trieste}: Efficiently exploring the depths of black-box functions with {TensorFlow}
(2023)
\doiurl{10.48550/ARXIV.2302.08436}
\end{botherref}
\endbibitem

\bibitem[\protect\citeauthoryear{Bartoli et~al.}{2019}]{Bartoli2019}
\begin{barticle}
\bauthor{\bsnm{Bartoli}, \binits{N.}},
\bauthor{\bsnm{Lefebvre}, \binits{T.}},
\bauthor{\bsnm{Dubreuil}, \binits{S.}},
\bauthor{\bsnm{Olivanti}, \binits{R.}},
\bauthor{\bsnm{Priem}, \binits{R.}},
\bauthor{\bsnm{Bons}, \binits{N.}},
\bauthor{\bsnm{Martins}, \binits{J.R.R.A.}},
\bauthor{\bsnm{Morlier}, \binits{J.}}:
\batitle{{Adaptive modeling strategy for constrained global optimization with application to aerodynamic wing design}}.
\bjtitle{Aerospace Science and Technology}
\bvolume{90},
\bfpage{85}--\blpage{102}
(\byear{2019})
\doiurl{10.1016/j.ast.2019.03.041}
\end{barticle}
\endbibitem

\bibitem[\protect\citeauthoryear{Bekemeyer et~al.}{2022}]{Bekemeyer2022}
\begin{bchapter}
\bauthor{\bsnm{Bekemeyer}, \binits{P.}},
\bauthor{\bsnm{Bertram}, \binits{A.}},
\bauthor{\bsnm{Hines~Chaves}, \binits{D.A.}},
\bauthor{\bsnm{Dias~Ribeiro}, \binits{M.}},
\bauthor{\bsnm{Garbo}, \binits{A.}},
\bauthor{\bsnm{Kiener}, \binits{A.}},
\bauthor{\bsnm{Sabater}, \binits{C.}},
\bauthor{\bsnm{Stradtner}, \binits{M.}},
\bauthor{\bsnm{Wassing}, \binits{S.}},
\bauthor{\bsnm{Widhalm}, \binits{M.}},
\bauthor{\bsnm{Goertz}, \binits{S.}},
\bauthor{\bsnm{Jaeckel}, \binits{F.}},
\bauthor{\bsnm{Hoppe}, \binits{R.}},
\bauthor{\bsnm{Hoffmann}, \binits{N.}}:
\bctitle{Data-driven aerodynamic modeling using the {DLR} {SMARTy} toolbox}.
In: \bbtitle{AIAA AVIATION 2022 Forum}.
\bpublisher{American Institute of Aeronautics and Astronautics},
\blocation{Chicago, USA}
(\byear{2022}).
\doiurl{10.2514/6.2022-3899}
\end{bchapter}
\endbibitem

\bibitem[\protect\citeauthoryear{García~Sánchez}{2024}]{GarciaSanchez2024}
\begin{botherref}
\oauthor{\bsnm{García~Sánchez}, \binits{R.}}:
Adaptation of an {MDO} platform for system architecture optimization.
mathesis,
Delft University of Technology,
Delft, NL
(January 2024)
\end{botherref}
\endbibitem

\bibitem[\protect\citeauthoryear{Hallé-Hannan et~al.}{2024}]{HalleHannan2024}
\begin{botherref}
\oauthor{\bsnm{Hallé-Hannan}, \binits{E.}},
\oauthor{\bsnm{Audet}, \binits{C.}},
\oauthor{\bsnm{Diouane}, \binits{Y.}},
\oauthor{\bsnm{Le~Digabel}, \binits{S.}},
\oauthor{\bsnm{Saves}, \binits{P.}}:
A graph-structured distance for heterogeneous datasets with meta variables.
Optimization Online
(2024)
\end{botherref}
\endbibitem

\bibitem[\protect\citeauthoryear{Baraton et~al.}{2023}]{Lucas_b}
\begin{bchapter}
\bauthor{\bsnm{Baraton}, \binits{L.}},
\bauthor{\bsnm{Urbano}, \binits{A.}},
\bauthor{\bsnm{Brevault}, \binits{L.}},
\bauthor{\bsnm{Balesdent}, \binits{M.}}:
\bctitle{Comparative review of multidisciplinary design analysis and optimization architectures for the preliminary design of a liquid rocket engine}.
In: \bbtitle{Aerospace Europe Conference 2023}
(\byear{2023})
\end{bchapter}
\endbibitem

\bibitem[\protect\citeauthoryear{Hung et~al.}{2009}]{HuJoMe2009}
\begin{barticle}
\bauthor{\bsnm{Hung}, \binits{Y.}},
\bauthor{\bsnm{Joseph}, \binits{V.R.}},
\bauthor{\bsnm{Melkote}, \binits{S.N.}}:
\batitle{{Design and Analysis of Computer Experiments With Branching and Nested Factors}}.
\bjtitle{Technometrics}
\bvolume{51}(\bissue{4}),
\bfpage{354}--\blpage{365}
(\byear{2009})
\doiurl{10.1198/TECH.2009.07097}
\end{barticle}
\endbibitem

\bibitem[\protect\citeauthoryear{Golovin et~al.}{2017}]{one-hot}
\begin{bchapter}
\bauthor{\bsnm{Golovin}, \binits{D.}},
\bauthor{\bsnm{Solnik}, \binits{B.}},
\bauthor{\bsnm{Moitra}, \binits{S.}},
\bauthor{\bsnm{Kochanski}, \binits{G.}},
\bauthor{\bsnm{Karro}, \binits{J.}},
\bauthor{\bsnm{Sculley}, \binits{D.}}:
\bctitle{Google vizier: A service for black-box optimization}.
In: \bbtitle{Proceedings of the 23rd ACM SIGKDD International Conference on Knowledge Discovery and Data Mining}
(\byear{2017})
\end{bchapter}
\endbibitem

\bibitem[\protect\citeauthoryear{Saves et~al.}{2022}]{Saves2022}
\begin{bchapter}
\bauthor{\bsnm{Saves}, \binits{P.}},
\bauthor{\bsnm{Diouane}, \binits{Y.}},
\bauthor{\bsnm{Bartoli}, \binits{N.}},
\bauthor{\bsnm{Lefebvre}, \binits{T.}},
\bauthor{\bsnm{Morlier}, \binits{J.}}:
\bctitle{A general square exponential kernel to handle mixed-categorical variables for gaussian process}.
In: \bbtitle{{AIAA} {AVIATION} 2022 Forum}.
\bpublisher{American Institute of Aeronautics and Astronautics},
\blocation{Chicago, IL, USA}
(\byear{2022}).
\doiurl{10.2514/6.2022-3870}
\end{bchapter}
\endbibitem

\bibitem[\protect\citeauthoryear{Renardy et~al.}{2021}]{Renardy2021}
\begin{barticle}
\bauthor{\bsnm{Renardy}, \binits{M.}},
\bauthor{\bsnm{Joslyn}, \binits{L.R.}},
\bauthor{\bsnm{Millar}, \binits{J.A.}},
\bauthor{\bsnm{Kirschner}, \binits{D.E.}}:
\batitle{To sobol or not to sobol? the effects of sampling schemes in systems biology applications}.
\bjtitle{Mathematical Biosciences}
\bvolume{337},
\bfpage{108593}
(\byear{2021})
\doiurl{10.1016/j.mbs.2021.108593}
\end{barticle}
\endbibitem

\bibitem[\protect\citeauthoryear{Kaltenecker et~al.}{2019}]{Kaltenecker2019}
\begin{bchapter}
\bauthor{\bsnm{Kaltenecker}, \binits{C.}},
\bauthor{\bsnm{Grebhahn}, \binits{A.}},
\bauthor{\bsnm{Siegmund}, \binits{N.}},
\bauthor{\bsnm{Guo}, \binits{J.}},
\bauthor{\bsnm{Apel}, \binits{S.}}:
\bctitle{Distance-based sampling of software configuration spaces}.
In: \bbtitle{2019 {IEEE}/{ACM} 41st International Conference on Software Engineering ({ICSE})}.
\bpublisher{{IEEE}},
\blocation{Montreal, CA}
(\byear{2019}).
\doiurl{10.1109/icse.2019.00112}
\end{bchapter}
\endbibitem

\bibitem[\protect\citeauthoryear{Collette et~al.}{2013}]{Collette2013}
\begin{bchapter}
\bauthor{\bsnm{Collette}, \binits{Y.}},
\bauthor{\bsnm{Hansen}, \binits{N.}},
\bauthor{\bsnm{Pujol}, \binits{G.}},
\bauthor{\bsnm{Aponte}, \binits{D.S.}},
\bauthor{\bsnm{Riche}, \binits{R.L.}}:
\bctitle{Object-Oriented Programming of Optimizers {\textendash} Examples in Scilab}.
In: \beditor{\bsnm{Breitkopf}, \binits{P.}},
\beditor{\bsnm{Coelho}, \binits{R.F.}} (eds.)
\bbtitle{Multidisciplinary Design Optimization in Computational Mechanics},
pp. \bfpage{499}--\blpage{538}.
\bpublisher{Wiley},
\blocation{London, UK}
(\byear{2013}).
\doiurl{10.1002/9781118600153.ch14}
\end{bchapter}
\endbibitem

\bibitem[\protect\citeauthoryear{Li et~al.}{2021}]{Li2021}
\begin{bchapter}
\bauthor{\bsnm{Li}, \binits{Y.}},
\bauthor{\bsnm{Shen}, \binits{Y.}},
\bauthor{\bsnm{Zhang}, \binits{W.}},
\bauthor{\bsnm{Chen}, \binits{Y.}},
\bauthor{\bsnm{Jiang}, \binits{H.}},
\bauthor{\bsnm{Liu}, \binits{M.}},
\bauthor{\bsnm{Jiang}, \binits{J.}},
\bauthor{\bsnm{Gao}, \binits{J.}},
\bauthor{\bsnm{Wu}, \binits{W.}},
\bauthor{\bsnm{Yang}, \binits{Z.}},
\bauthor{\bsnm{Zhang}, \binits{C.}},
\bauthor{\bsnm{Cui}, \binits{B.}}:
\bctitle{{OpenBox}: A generalized black-box optimization service}.
In: \bbtitle{Proceedings of the 27th {ACM} {SIGKDD} Conference on Knowledge Discovery and Data Mining}.
\bpublisher{{ACM}},
\blocation{New York, USA}
(\byear{2021}).
\doiurl{10.1145/3447548.3467061}
\end{bchapter}
\endbibitem

\bibitem[\protect\citeauthoryear{Zhang et~al.}{2020}]{Zhang2020}
\begin{bchapter}
\bauthor{\bsnm{Zhang}, \binits{Q.}},
\bauthor{\bsnm{Han}, \binits{Z.}},
\bauthor{\bsnm{Yang}, \binits{F.}},
\bauthor{\bsnm{Zhang}, \binits{Y.}},
\bauthor{\bsnm{Liu}, \binits{Z.}},
\bauthor{\bsnm{Yang}, \binits{M.}},
\bauthor{\bsnm{Zhou}, \binits{L.}}:
\bctitle{Retiarii: A deep learning {Exploratory-Training} framework}.
In: \bbtitle{14th USENIX Symposium on Operating Systems Design and Implementation (OSDI 20)},
pp. \bfpage{919}--\blpage{936}.
\bpublisher{USENIX Association},
\blocation{Virtual}
(\byear{2020}).
\burl{https://www.usenix.org/conference/osdi20/presentation/zhang-quanlu}
\end{bchapter}
\endbibitem

\bibitem[\protect\citeauthoryear{Elsken et~al.}{2019}]{Elsken2019}
\begin{barticle}
\bauthor{\bsnm{Elsken}, \binits{T.}},
\bauthor{\bsnm{Metzen}, \binits{J.H.}},
\bauthor{\bsnm{Hutter}, \binits{F.}}:
\batitle{Neural architecture search: A survey}.
\bjtitle{Journal of Machine Learning Research}
\bvolume{20}(\bissue{55}),
\bfpage{1}--\blpage{21}
(\byear{2019})
\end{barticle}
\endbibitem

\bibitem[\protect\citeauthoryear{Gray et~al.}{2019}]{Gray2019}
\begin{barticle}
\bauthor{\bsnm{Gray}, \binits{J.S.}},
\bauthor{\bsnm{Hwang}, \binits{J.T.}},
\bauthor{\bsnm{Martins}, \binits{J.R.R.A.}},
\bauthor{\bsnm{Moore}, \binits{K.T.}},
\bauthor{\bsnm{Naylor}, \binits{B.A.}}:
\batitle{{OpenMDAO}: an open-source framework for multidisciplinary design, analysis, and optimization}.
\bjtitle{Structural and Multidisciplinary Optimization}
\bvolume{59}(\bissue{4}),
\bfpage{1075}--\blpage{1104}
(\byear{2019})
\doiurl{10.1007/s00158-019-02211-z}
\end{barticle}
\endbibitem

\bibitem[\protect\citeauthoryear{Hendricks and Gray}{2019}]{Hendricks2019}
\begin{barticle}
\bauthor{\bsnm{Hendricks}, \binits{E.S.}},
\bauthor{\bsnm{Gray}, \binits{J.S.}}:
\batitle{{pyCycle}: A tool for efficient optimization of gas turbine engine cycles}.
\bjtitle{Aerospace}
\bvolume{6}(\bissue{8}),
\bfpage{87}
(\byear{2019})
\doiurl{10.3390/aerospace6080087}
\end{barticle}
\endbibitem

\bibitem[\protect\citeauthoryear{Donelli et~al.}{2023}]{Donelli2023}
\begin{botherref}
\oauthor{\bsnm{Donelli}, \binits{G.}},
\oauthor{\bsnm{Ciampa}, \binits{P.D.}},
\oauthor{\bsnm{Mello}, \binits{J.M.G.}},
\oauthor{\bsnm{Odaguil}, \binits{F.I.K.}},
\oauthor{\bsnm{Cuco}, \binits{A.P.C.}},
\oauthor{\bsnm{Laan}, \binits{T.}}:
A value-driven concurrent approach for aircraft design-manufacturing-supply chain.
Production and Manufacturing Research
\textbf{11}(1)
(2023)
\doiurl{10.1080/21693277.2023.2279709}
\end{botherref}
\endbibitem

\bibitem[\protect\citeauthoryear{Nikolentzos et~al.}{2021}]{Nikolentzos2021}
\begin{barticle}
\bauthor{\bsnm{Nikolentzos}, \binits{G.}},
\bauthor{\bsnm{Siglidis}, \binits{G.}},
\bauthor{\bsnm{Vazirgiannis}, \binits{M.}}:
\batitle{Graph kernels: A survey}.
\bjtitle{Journal of Artificial Intelligence Research}
\bvolume{72},
\bfpage{943}--\blpage{1027}
(\byear{2021})
\doiurl{10.1613/jair.1.13225}
\end{barticle}
\endbibitem

\bibitem[\protect\citeauthoryear{Sirico and Herber}{2023}]{Sirico2023}
\begin{botherref}
\oauthor{\bsnm{Sirico}, \binits{A.}},
\oauthor{\bsnm{Herber}, \binits{D.R.}}:
On the use of geometric deep learning for the iterative classification and down-selection of analog electric circuits.
Journal of Mechanical Design
\textbf{146}(5)
(2023)
\doiurl{10.1115/1.4063659}
\end{botherref}
\endbibitem

\end{thebibliography}
